\pdfoutput=1
\RequirePackage{silence}
\documentclass[a4paper,11pt]{amsart}
\usepackage[hmarginratio={1:1},vmarginratio={1:1},lmargin=55.0pt,tmargin=55.0pt]{geometry}

\usepackage[numbers]{natbib}

\allowdisplaybreaks

\usepackage[utf8]{inputenc}

\usepackage{amssymb,amsmath,amsthm,amsfonts,mathrsfs,bbm}
\usepackage[table]{xcolor}
\usepackage{graphicx}
\usepackage{mathtools}

\definecolor{spinach}{RGB}{46,139,87}
\definecolor{tomato}{RGB}{255,99,71}
\definecolor{pumpkin}{RGB}{224,180,80}

\definecolor{orchid}{RGB}{143,40,194}
\definecolor{lava}{RGB}{207,16,32}
\definecolor{mydarkblue}{RGB}{10,10,150}

\usepackage[yyyymmdd,hhmmss]{datetime}
\newcommand\GL[2][]{\todo[color=lava!20,inline,#1]{GL: #2}}
\newcommand\DT[2][]{\todo[color=orchid!20,inline,#1]{DT: #2}}

\def\changed#1{#1}
\def\ochanged#1{#1}

\usepackage{xparse}
\def\citeform#1{\textbf{#1}}

\usepackage{enumitem}
\setlist[enumerate]{itemsep=0.15cm,label=\emph{\upshape(\alph*)}}
\setlist[enumerate,2]{itemsep=0.15cm,label=\emph{\upshape(\roman*)}}

\let\emph\relax
\DeclareTextFontCommand{\emph}{\em}

\usepackage{array}
\newcolumntype{C}{>{$}c<{$}}

\usepackage[vcentermath]{youngtab}

\usepackage[all]{xy}
\usepackage{tikz}
\tikzstyle{densely dotted}=[dash pattern=on \pgflinewidth off .5pt]
\tikzset{anchorbase/.style={baseline={([yshift=-0.5ex]current bounding box.center)}},
tinynodes/.style={font=\tiny, text height=0.25ex, text depth=0.05ex},
smallnodes/.style={font=\scriptsize, text height=0.75ex, text depth=0.15ex},
usual/.style={line width=.9,color=black},
crossline/.style={preaction={draw=white,line width=5.0pt,-},preaction={draw=black,line width=.9pt,-}},
mor/.style={line width=0.75,color=black,fill=yellow!50},
}
\usetikzlibrary{cd}
\usetikzlibrary{decorations}
\usetikzlibrary{decorations.markings}
\usetikzlibrary{decorations.pathreplacing}
\usetikzlibrary{decorations.pathmorphing}
\usetikzlibrary{arrows.meta,shapes,positioning,matrix,calc}
\usetikzlibrary{shapes.callouts}
\usetikzlibrary{shapes.geometric}
\tikzstyle directed=[postaction={decorate,decoration={markings,mark=at position #1 with {\arrow[line width=0.3mm, black]{>}}}}]
\tikzstyle rdirected=[postaction={decorate,decoration={markings,mark=at position #1 with {\arrow[line width=0.3mm, black]{<}}}}]

\newcommand{\qbin}[2]{\genfrac{[}{]}{0pt}{}{#1}{#2}}

\usepackage{aliascnt}
\newtheorem{Theorem}{Theorem}[section]
\newaliascnt{Proposition}{Theorem}
\newtheorem{Proposition}[Proposition]{Proposition}
\aliascntresetthe{Proposition}
\expandafter\def\csname Propositionautorefname\endcsname{Proposition}
\newaliascnt{Corollary}{Theorem}
\newtheorem{Corollary}[Corollary]{Corollary}
\aliascntresetthe{Corollary}
\expandafter\def\csname Corollaryautorefname\endcsname{Corollary}
\newaliascnt{Lemma}{Theorem}
\newtheorem{Lemma}[Lemma]{Lemma}
\aliascntresetthe{Lemma}
\expandafter\def\csname Lemmaautorefname\endcsname{Lemma}
\theoremstyle{definition}
\newaliascnt{Definition}{Theorem}
\newtheorem{Definition}[Definition]{Definition}
\aliascntresetthe{Definition}
\expandafter\def\csname Definitionautorefname\endcsname{Definition}
\def\Theoremautorefname{Theorem}
\def\equationautorefname~#1\null{(#1)\null}

\numberwithin{equation}{section}
\renewcommand{\thesubsection}{\thesection.\arabic{subsection}}
\renewcommand{\theequation}{\thesection.\arabic{equation}}
\newaliascnt{Notation}{Theorem}

\aliascntresetthe{Notation}
\expandafter\def\csname Notationautorefname\endcsname{Notation}
\newaliascnt{Example}{Theorem}
\newtheorem{Example}[Example]{Example}
\aliascntresetthe{Example}
\expandafter\def\csname Exampleautorefname\endcsname{Example}
\newaliascnt{Examples}{Theorem}

\aliascntresetthe{Examples}
\expandafter\def\csname Examplesautorefname\endcsname{Examples}
\def\Examplesautorefname{Example}
\theoremstyle{remark}
\newaliascnt{Remark}{Theorem}
\newtheorem{Remark}[Remark]{Remark}
\aliascntresetthe{Remark}
\expandafter\def\csname Remarkautorefname\endcsname{Remark}
\newaliascnt{Assumption}{Theorem}

\aliascntresetthe{Assumption}
\expandafter\def\csname Assumptionautorefname\endcsname{Assumption}
\newaliascnt{Construction}{Theorem}

\aliascntresetthe{Construction}
\expandafter\def\csname Constructionautorefname\endcsname{Construction}

\usepackage[hypertexnames=false]{hyperref}
\usepackage{bookmark}
\hypersetup{
pdftoolbar=true,
pdfmenubar=true,
pdffitwindow=false,
pdfstartview={FitH},
pdftitle={Minimal presentations of \texorpdfstring{$\mathfrak{gl}_n$}{gln}-web categories},
pdfauthor={Genta Latifi and Daniel Tubbenhauer},
pdfsubject={},
pdfcreator={Genta Latifi and Daniel Tubbenhauer},
pdfproducer={Minimal presentations of \texorpdfstring{$\mathfrak{gl}_n$}{gln}-web categories},
pdfkeywords={},
pdfnewwindow=true,
colorlinks=true,
linkcolor=mydarkblue,
citecolor=teal,
filecolor=magenta,
urlcolor=orchid,
linkbordercolor=lava,
citebordercolor=teal,
urlbordercolor=orchid,
linktocpage=true
}

\def\makeautorefname#1#2{\csdef{#1autorefname}{#2}}
\makeautorefname{section}{Section}%
\makeautorefname{subsection}{Section}%
\makeautorefname{subsubsection}{Section}%

\begin{document}
\title[Minimal presentations of gln-web categories]
{Minimal presentations of gln-web categories}
\author[G. Latifi and D. Tubbenhauer]{Genta Latifi and Daniel Tubbenhauer}

\address{G.L.: Universit{\"a}t Z{\"u}rich, Institut f{\"u}r Mathematik, Winterthurerstrasse 190, CH-8057, Z{\"u}rich, \newline
Switzerland}
\email{genta.latifi@math.uzh.ch}

\address{D.T.: The University of Sydney, School of Mathematics and Statistics F07, Office Carslaw 827, NSW 2006, Australia, \href{http://www.dtubbenhauer.com}{www.dtubbenhauer.com}, \href{https://orcid.org/0000-0001-7265-5047}{ORCID 0000-0001-7265-5047}}
\email{daniel.tubbenhauer@sydney.edu.au}

\begin{abstract}
In this paper we study categories of $\mathfrak{gl}_n$-webs 
which describe associated representation categories of the quantum group $\text{U}_{q}(\mathfrak{gl}_n)$. We give a 
minimal presentation of the category of $\mathfrak{gl}_n$-webs over a field with generic quantum parameters. We additionally describe an integral presentation which differs from others in the literature 
because it is ``as coefficient-free as possible''.
\end{abstract}


\maketitle

\tableofcontents

\arrayrulewidth=0.5mm
\setlength{\arrayrulewidth}{0.5mm}

\section{Introduction}

\setcounter{subsection}{1}

In this paper we study a category of planar graphs coloured by non-negative integers, called a web category. \changed{Web categories, in various incarnations, have been around for donkey's years.}

One of the classical examples of diagrammatic categories is the Temperley--Lieb category, denoted by $\mathcal{TL}$. It was independently introduced by Rumer--Teller--Weyl in their study of invariant tensors and by Temperley--Lieb \ochanged{in} statistical mechanics. 

Let $q$ be a generic formal parameter and let $\mathbb{C}(q)$ be the field of rational functions in $q$. The category $\mathcal{TL}$ is the $\mathbb{C}(q)$-linear category with objects being non-negative integers $k \in\mathbb{Z}_{\geq 0}$, and the morphisms are generated by (taking tensor products $\otimes$ and compositions $\circ$ and $\mathbb{C}(q)$-linear \ochanged{combinations}) of the non-intersecting planar diagrams
of the form
\[
\xy
(0,2)*{
\begin{tikzpicture} [anchorbase,scale=1]
\draw[very thick] (0,0) to (0,1);
\end{tikzpicture}
}
\endxy\quad , \quad
\xy
(0,2)*{
\begin{tikzpicture} [anchorbase,scale=1]
\draw[very thick] (0,0) to [out=90,in=180] (.25,.5) to [out=0,in=90] (.5,0);
\node at (.25,.9) {};
\end{tikzpicture}
}
\endxy\quad , \quad
\xy
(0,-3)*{
\begin{tikzpicture} [anchorbase,scale=1]
\draw[very thick] (0,1) to [out=270,in=180] (.25,.5) to [out=0,in=270] (.5,1);
\node at (.25,0.1) {};
\end{tikzpicture}
}
\endxy 
\]
subject to the relation 
\[
\xy
(0,0)*{
\begin{tikzpicture} [anchorbase,scale=1]
\draw[very thick] (-1,0) circle (0.3cm);
\end{tikzpicture}}
\endxy=-[2],
\]
where $[2]=q+q^{-1}$. The compositions and tensor products of diagrams are illustrated below:
\[
\begin{tikzpicture} [anchorbase,scale=1]
\draw[very thick] (0,0) to [out=90,in=180] (.25,.5) to [out=0,in=90] (.5,0);
\node at (.25,.9) {};
\end{tikzpicture}
\ \circ \ 
\begin{tikzpicture} [anchorbase,scale=1]
\draw[very thick] (0,1) to [out=270,in=180] (.25,.5) to [out=0,in=270] (.5,1);
\node at (.25,0.1) {};
\end{tikzpicture}
=
\begin{tikzpicture} [anchorbase,scale=1]
\draw[very thick] (-1,0) circle (0.3cm);
\end{tikzpicture}
\quad , \quad
\begin{tikzpicture} [anchorbase,scale=1]
\draw[very thick] (0,0) to (0,1);
\end{tikzpicture}
\quad \otimes \
\begin{tikzpicture} [anchorbase,scale=1]
\draw[very thick] (0,0) to [out=90,in=180] (.25,.5) to [out=0,in=90] (.5,0);
\node at (.25,.9) {};
\end{tikzpicture}
=\quad
\begin{tikzpicture} [anchorbase,scale=1]
\draw[very thick] (0,0) to (0,1);
\end{tikzpicture}
\quad
\begin{tikzpicture} [anchorbase,scale=1]
\draw[very thick] (0,0) to [out=90,in=180] (.25,.5) to [out=0,in=90] (.5,0);
\node at (.25,.9) {};
\end{tikzpicture} \ .
\] 
Let $\text{U}_q(\mathfrak{sl}_2)$ be the quantum special linear group. \ochanged{(Without further notice, we consider only left modules of type $1$, the latter whenever applicable.)} Let $\text{U}_q(\mathfrak{sl}_2)$-$\mathbf{fdMod}_\wedge$ $=$ $\text{U}_q(\mathfrak{sl}_2)$-$\mathbf{fdMod}_S$ be the full subcategory of finite-dimensional $\text{U}_q(\mathfrak{sl}_2)$-modules generated by tensor products of the two-dimensional vector representation $\mathbb{C}^2_q = \mathbb{C}(q)^2$ of $\text{U}_q(\mathfrak{sl}_2)$. From the work of Rumer, Teller, and Weyl in \cite{RuTeWe-sl2}, it follows that $\mathcal{TL}$ describes $\text{U}_q(\mathfrak{sl}_2)$-$\mathbf{fdMod}_\wedge$. More precisely, we have the following.

\begin{Theorem}(Rumer--Teller--Weyl \cite{RuTeWe-sl2}.)
The category $\mathcal{TL}$ and $\text{U}_q(\mathfrak{sl}_2)$-$\mathbf{fdMod}_\wedge$ are equivalent as pivotal categories.
\end{Theorem}

It follows from this theorem that the additive Karoubi closure of $\mathcal{TL}$ recovers 
all of $\text{U}_q(\mathfrak{sl}_2)$-$\mathbf{fdMod}$.

Generalizing these ideas, Kuperberg \cite{Ku-spiders-rank-2} gave a diagrammatic description of $\text{U}_q(\mathfrak{sl}_3)$-$\mathbf{fdMod}_{\wedge}$, the full subcategory of finite-dimensional $\text{U}_q(\mathfrak{sl}_3)$-modules whose objects are finite tensor products of $\bigwedge_q^k\mathbb{C}_q^3$, the fundamental $\text{U}_q(\mathfrak{sl}_3)$-modules. 
As for the Temperley--Lieb case, taking the additive Karoubi closure of Kuperberg's web category describes the entire category $\text{U}_q(\mathfrak{sl}_3)$-$\mathbf{fdMod}$.

\smallskip

Generalising Kuperberg's result to $\text{U}_q(\mathfrak{sl}_n)$-$\mathbf{fdMod}_{\wedge}$ \changed{was} an open problem for quite a while, but some partial answers were known before. 
The question was then 
fully answered by Cautis, Kamnitzer, \ochanged{and} Morrison in \cite{CaKaMo-webs-skew-howe}.
\changed{Similar questions were then addressed in many works, e.g. in \cite{RoTu-symmetric-howe}, which studies symmetric powers instead of exterior 
powers.}

In this paper we study the category  $\textbf{SWeb}_{\uparrow,\downarrow}(\mathfrak{gl}_n)$ of $\mathfrak{gl}_n$-webs, whose objects are words from $\{k_\uparrow,k_\downarrow \ | \ k \in \mathbb{Z}_{\geq 0} \}$. We say that the object $k_\uparrow$ is upward oriented and $k_\downarrow$ is downward oriented. A (finite) word in these symbols is \changed{decorated by a vector arrow, e.g.} $\vec{k}$. Then a morphism from $\vec{k}$ to $\vec{l}$, called a web, is a planar diagram with bottom boundary $\vec{k}$ and top boundary $\vec{l}$, which is obtained by piecing together any of the generating diagrams:
\begin{equation} \label{webs intro}
\begin{tikzpicture}[anchorbase,scale=.3]
\draw [very thick, directed=.55] (0, .75) to (0,2.5);
\draw [very thick, directed=.45] (1,-1) to [out=90,in=330] (0,.75);
\draw [very thick, directed=.45] (-1,-1) to [out=90,in=210] (0,.75); 
\node at (0, 3) {\tiny $k{+}l$};
\node at (-1,-1.5) {\tiny $k$};
\node at (1,-1.5) {\tiny $l$};
\end{tikzpicture}
, \quad
\begin{tikzpicture}[anchorbase,scale=.3]
\draw [very thick, directed=.55] (0,-1) to (0,.75);
\draw [very thick, directed=.65] (0,.75) to [out=30,in=270] (1,2.5);
\draw [very thick, directed=.65] (0,.75) to [out=150,in=270] (-1,2.5); 
\node at (0, -1.5) {\tiny $k{+}l$};
\node at (-1,3) {\tiny $k$};
\node at (1,3) {\tiny $l$};
\end{tikzpicture}
, \quad
\begin{tikzpicture}[anchorbase,scale=.3]
\draw [very thick, rdirected=.55] (0, .75) to (0,2.5);
\draw [very thick, rdirected=.45] (1,-1) to [out=90,in=330] (0,.75);
\draw [very thick, rdirected=.45] (-1,-1) to [out=90,in=210] (0,.75); 
\node at (0, 3) {\tiny $k{+}l$};
\node at (-1,-1.5) {\tiny $k$};
\node at (1,-1.5) {\tiny $l$};
\end{tikzpicture}
, \quad
\begin{tikzpicture}[anchorbase,scale=.3]
\draw [very thick, rdirected=.55] (0,-1) to (0,.75);
\draw [very thick, rdirected=.65] (0,.75) to [out=30,in=270] (1,2.5);
\draw [very thick, rdirected=.65] (0,.75) to [out=150,in=270] (-1,2.5); 
\node at (0, -1.5) {\tiny $k{+}l$};
\node at (-1,3) {\tiny $k$};
\node at (1,3) {\tiny $l$};
\end{tikzpicture}
, \quad
\begin{tikzpicture} [anchorbase,scale=1]
\draw[very thick,rdirected=-.95] (0,1) to [out=270,in=180] (.25,.5) to [out=0,in=270] (.5,1);
\node at (.25,.25) {};
\node at (0,1.25) {\tiny $k$};
\node at (.5,1.25) {\tiny $k$};
\end{tikzpicture}
, \quad 
\begin{tikzpicture} [anchorbase,scale=1]
\draw[very thick,rdirected=-.95] (0,0) to [out=90,in=180] (.25,.5) to [out=0,in=90] (.5,0);
\node at (.25,.75) {};
\node at (0,-.25) {\tiny $k$};
\node at (.5,-.25) {\tiny $k$};
\end{tikzpicture}, \quad 
\begin{tikzpicture} [anchorbase,scale=1]
\draw[very thick,directed=.95] (0,1) to [out=270,in=180] (.25,.5) to [out=0,in=270] (.5,1);
\node at (.25,.25) {};
\node at (0,1.25) {\tiny $k$};
\node at (.5,1.25) {\tiny $k$};
\end{tikzpicture}
, \quad 
\begin{tikzpicture} [anchorbase,scale=1]
\draw[very thick,directed=.95] (0,0) to [out=90,in=180] (.25,.5) to [out=0,in=90] (.5,0);
\node at (.25,.75) {};
\node at (0,-.25) {\tiny $k$};
\node at (.5,-.25) {\tiny $k$};
\end{tikzpicture}
.
\end{equation}
\changed{General morphisms are linear combinations of such webs.}

Here, ``piecing together'' means glueing a diagram $w_2$ on top of another diagram $w_1$ (which will correspond to a composition of diagrams $w_2 \circ w_1$) and putting a diagram $w_2$ to the right of a diagram $w_1$ (which will correspond to a tensor product $w_1 \otimes w_2$). For example, 
\[
\begin{tikzpicture}[anchorbase,scale=.3]
\draw [very thick, directed=.6] (0,-1) to (0,.75);
\draw [very thick, directed=.75] (0,.75) to [out=30,in=270] (1,2.5);
\draw [very thick, directed =.75] (0,.75) to [out=150,in=270] (-1,2.5); 
\draw [very thick, directed=-.5] (1,-2.75) to [out=90,in=330] (0,-1);
\draw [very thick, directed=-.5] (-1,-2.75) to [out=90,in=210] (0,-1);
\draw [very thick] (1,2.5) to [in=180, out=90] (2,3.5) to [in=90, out=0](3,2.5);
\draw [very thick,rdirected=.5] (3,-2.75) to (3,2.5);
\draw [very thick,] (-1,2.5) to (-1,3.5);
\node at (-1,4) {\tiny $1$};
\node at (3,-3.15) {\tiny $2$};
\node at (-1,-3.15) {\tiny $2$};
\node at (1,-3.15) {\tiny $1$};
\node at (-1.2,0) {\tiny $3$};
\end{tikzpicture}
= \left(
\begin{tikzpicture}[anchorbase,anchorbase,scale=0.75]
\draw [very thick,directed=.5] (0,0) to (0,2);
\node at (0,-.2) {\tiny $1$};
\node at (0,2.2) {\tiny $1$};
\end{tikzpicture}
\otimes
\begin{tikzpicture}[anchorbase,scale=0.75]
\draw [very thick,directed=.95] (0,0)
to[out=90,in=180] (0.5,0.5) to[out=0,in=90] (1,0);
\node at (0,-.2) {\tiny $2$};
\node at (1,-.2) {\tiny $2$};
\end{tikzpicture}
\right)
\circ
\left(
\begin{tikzpicture}[anchorbase,anchorbase,scale=0.75]
\draw [very thick,directed=.5] (0.5,1) to[out=150,in=270] (0,2);
\draw [very thick,directed=.5] (0.5,1) to[out=30,in=270] (1,2);
\draw [very thick,directed=.5] (0.5,0) to (0.5,1);
\node at (0.5,-.2) {\tiny $3$};
\node at (0,2.2) {\tiny $1$};
\node at (1,2.2) {\tiny $2$};
\end{tikzpicture}
\otimes
\begin{tikzpicture}[anchorbase,scale=0.75]
\draw [very thick,directed=.5] (0,2) to (0,0);
\node at (0,-.2) {\tiny $2$};
\node at (0,2.2) {\tiny $2$};
\end{tikzpicture}
\right)
\circ
\left(
\begin{tikzpicture}[anchorbase,scale=0.75]
\draw [very thick,directed=.5] (0,0) to[out=90,in=210] (0.5,1);
\draw [very thick,directed=.5] (1,0) to[out=90,in=330] (0.5,1);
\draw [very thick,directed=.5] (0.5,1) to (0.5,2);
\node at (0,-.2) {\tiny $2$};
\node at (1,-.2) {\tiny $1$};
\node at (0.5,2.2) {\tiny $3$};
\end{tikzpicture}
\otimes
\begin{tikzpicture}[anchorbase,scale=0.75]
\draw [very thick,directed=.5] (0,2) to (0,0);
\node at (0,-.2) {\tiny $2$};
\node at (0,2.2) {\tiny $2$};
\end{tikzpicture}
\right)
\]
is a morphism from $\vec{k}=2_\uparrow \otimes 1_\uparrow \otimes 2_\downarrow$ to $\vec{l}= 1_\uparrow$.
As we will explain, this category of $\mathfrak{gl}_n$-webs describes in a ``minimal'' way (following Turaev's minimal presentation of the tangle categories \cite{Tu-cat-tangles}; see also \cite{BrDaKu-qwebs-typep} which served as a motivation for us), the entire category of finite-dimensional $\text{U}_q(\mathfrak{gl}_n)$-modules. 

The paper is organized as follows.
We first start \ochanged{with} the full subcategory $\text{U}_{q}(\mathfrak{gl}_n) \text{-} \mathrm{\mathbf{fdMod}}_{S,S^*}$ generated by tensor products of symmetric powers \ochanged{$\text{Sym}_q^k\mathbb{C}_q^n$} ($k \in \mathbb{N})$ of the standard $\text{U}_q(\mathfrak{gl}_n)$-representation and their duals. We point out here that in contrast to \cite{CaKaMo-webs-skew-howe}, where $(\bigwedge_q^k\mathbb{C}_q^n)^* \cong \bigwedge_q^{n-k}\mathbb{C}_q^n$, we do not have such an isomorphism in the case of symmetric powers. Moreover, we explain that we have described the category of representations in a minimal way, namely that we have used a minimal set of generators (only those webs in \autoref{webs intro} with $k=1$ or $l=1$) and a minimal set of defining relations for the web category to describe the category of modules. 

Finally, we give an integral presentation of $\mathfrak{gl}_n$-webs which is different from the ones in the literature. Whether this integral 
presentation also describes $\text{U}_{q}(\mathfrak{gl}_n) \text{-} \mathrm{\mathbf{fdMod}}_{S,S^*}$ is an open problem, but we show that our web category surjects onto this category.

\ochanged{(In fact, we work in slightly greater generality, allowing e.g. any field with a transcendental parameter $q$. For example, working over a finite field instead of $\mathbb{C}$ and a generic variable $q$ is allowed.)}

\changed{For the reader familiar with \cite{TuVaWe-super-howe}, 
conceptually, our upward symmetric category is the red only part of the red-green calculus of that paper, but here we recast it as a standalone minimal presentation.}
\medskip

\noindent\textbf{Acknowledgements.}
We would like to thank Anna Beliakova, Sasha Kleshchev, \ochanged{the referee}, Catharina Stroppel and Emmanuel Wagner for
freely sharing ideas, patiently answering questions and helpful discussions.

G.L. would like to \ochanged{thank} MSRI for sponsoring a research stay during the program Quantum Symmetries where parts of the mathematics in this paper \ochanged{were} discovered.

This paper is part of the first \ochanged{author's} Ph.D. thesis 
at the Universit{\"a}t Z{\"u}rich, and the support 
of the Universit{\"a}t Z{\"u}rich \ochanged{over the years} is gratefully 
acknowledged. None of us \ochanged{is} at the Universit{\"a}t Z{\"u}rich anymore, but it was 
a great time!

\section{Web categories}

We use the following conventions for \ochanged{\textit{quantum numbers}}. 
Let $s\in\mathbb{Z}$ and $t\in\mathbb{Z}_{\geq 0}$.
By convention, $[0]!=1=\begin{bsmallmatrix} s\\ 0 \end{bsmallmatrix}$ and 
\[
[s]=\frac{q^{s}-q^{-s}}{q-q^{-1}},\quad
[t]!=[t][t-1]\dots[1],\quad
\begin{bmatrix} s\\ t \end{bmatrix}
=
\frac{[s][s-1]\dots[s-t+1]}{[t][t-1]\dots[1]}.
\]
\ochanged{We also define the}
\textit{quantum trinomial coefficient} via 
\begin{center}
$\!\begin{bmatrix} a \\ b,c \end{bmatrix}\!:=\frac{[a][a-1] \cdots [a-b-c+1]}{[b]![c]!} . $
\end{center}
Fix a non-negative integer $n$. This integer will be fixed here and throughout. Let 
$\Bbbk$ be any field. Let $\Bbbk(q^\frac{1}{n})$ be the field of rational functions over $q^\frac{1}{n}$, where $q$ is chosen such 
that all quantum numbers are invertible in $\Bbbk(q^\frac{1}{n})$. 
\ochanged{(For example, $q$ could be a generic parameter and $\Bbbk$ arbitrary, or $q=1$ and $\Bbbk$ is of characteristic zero.)}

\begin{Remark}
The $n$th root of $q$ is only needed to match the braidings on webs with 
the braiding coming from the $R$-matrix, see \autoref{up-vs-down}. 
The $n$th root does not play 
any other role.
\end{Remark}

\subsection{Definition of web categories}

\begin{Definition}\label{free sym up webs} 
The \textit{free thin symmetric upward pointing category} of $\mathfrak{gl}_n$-webs $\mathbf{fSWeb}_{\uparrow}(\mathfrak{gl}_n)$ is the $\Bbbk(q^\frac{1}{n})$-linear category monoidally generated by\\
$\bullet $ objects $\{k_{\uparrow} \ | \ k\in \mathbb{Z}_{\geq 0} \}$,\\
$\bullet$ morphisms 
\[
\xy
(0,0)*{
\begin{tikzpicture}[scale=.3]
\draw [very thick, directed=.55] (0, .75) to (0,2.5);
\draw [very thick, directed=.45] (1,-1) to [out=90,in=330] (0,.75);
\draw [very thick, directed=.45] (-1,-1) to [out=90,in=210] (0,.75); 
\node at (0, 3) {\tiny $k{+}1$};
\node at (-1,-1.5) {\tiny $k$};
\node at (1,-1.5) {\tiny $1$};
\end{tikzpicture}  
};
\endxy : k_{\uparrow} \otimes 1_{\uparrow} \rightarrow (k+1)_{\uparrow} , \quad
\xy
(0,0)*{
\begin{tikzpicture}[scale=.3]
\draw [very thick, directed=.55] (0,-1) to (0,.75);
\draw [very thick, directed=.65] (0,.75) to [out=30,in=270] (1,2.5);
\draw [very thick, directed=.65] (0,.75) to [out=150,in=270] (-1,2.5); 
\node at (0, -1.5) {\tiny $k{+}1$};
\node at (-1,3) {\tiny $k$};
\node at (1,3) {\tiny $1$};
\end{tikzpicture}
};
\endxy :(k+1)_{\uparrow} \rightarrow k_{\uparrow} \otimes 1_{\uparrow} ,
\]
called thin merges (or $(k,1)$-merges) and thin splits (or $(k,1)$-splits), respectively.
\end{Definition}

Here and throughout, we will use the conventions of \cite{RoTu-symmetric-howe} for reading the diagrams, meaning that $v \circ  u$ is obtained by glueing $v$ on top of $u$ and $ u \otimes v$ is obtained by putting $v$ to the right of $u$, for $u, v \in \text{Hom}_{\mathbf{fSWeb}_{ \uparrow}(\mathfrak{gl}_n)}(\vec{k}, \vec{l}).$ For example, 

\begin{center}

\[
\xy
(0,0)*{
\begin{tikzpicture} [scale=1]
\draw[very thick] (0,0) to [out=90,in=180] (.25,.5) to [out=0,in=90] (.5,0);
\node at (.25,.75) {};
\node at (0,-.15) {\tiny $k$};
\node at (.5,-.15) {\tiny $k$};
\end{tikzpicture}
}
\endxy
\circ
\xy
(0,0)*{
\begin{tikzpicture} [scale=1]
\draw[very thick] (0,1) to [out=270,in=180] (.25,.5) to [out=0,in=270] (.5,1);
\node at (.25,.25) {};
\node at (0,1.15) {\tiny $k$};
\node at (.5,1.15) {\tiny $k$};
\end{tikzpicture}
}
\endxy
=\xy
(0,0)*{
\begin{tikzpicture} [scale=1]
\draw[very thick] (-1,0) circle (0.3cm);
\node at (-0.5,0) {\tiny $k$};
\end{tikzpicture}}
\endxy ,
\xy
(0,0)*{
\begin{tikzpicture} [scale=1]
\draw[very thick] (0,1) to [out=270,in=180] (.25,.5) to [out=0,in=270] (.5,1);
\node at (.25,.25) {};
\node at (0,1.15) {\tiny $k$};
\node at (.5,1.15) {\tiny $k$};
\end{tikzpicture}
}
\endxy
\circ
\xy
(0,0)*{
\begin{tikzpicture} [scale=1]
\draw[very thick] (0,0) to [out=90,in=180] (.25,.5) to [out=0,in=90] (.5,0);
\node at (.25,.75) {};
\node at (0,-.15) {\tiny $k$};
\node at (.5,-.15) {\tiny $k$};
\end{tikzpicture}
}
\endxy
=
\xy
(0,0)*{
\begin{tikzpicture}[scale=.3]
\draw [very thick] (-1,2.5) to [out=270,in=180] (0,0.5) to [out=0,in=270] (1,2.5);
\draw [very thick] (-1,-2.75) to [out=90,in=180] (0,-0.75) to [out=0,in=90] (1,-2.75);
\node at (-1,3) {\tiny $k$};
\node at (1,3) {\tiny $k$};
\node at (-1,-3.15) {\tiny $k$};
\node at (1,-3.15) {\tiny $k$};
\end{tikzpicture}
};
\endxy,
\xy
(0,0)*{
\begin{tikzpicture}[scale=.3]
\draw [very thick] (0,-1) to (0,.75);
\draw [very thick] (0,.75) to [out=30,in=270] (1,2.5);
\draw [very thick] (0,.75) to [out=150,in=270] (-1,2.5); 
\draw [very thick] (1,-2.75) to [out=90,in=330] (0,-1);
\draw [very thick] (-1,-2.75) to [out=90,in=210] (0,-1);
\node at (-1,3) {\tiny $l_1$};
\node at (1,3) {\tiny $l_2$};
\node at (-1,-3.15) {\tiny $k_1$};
\node at (1,-3.15) {\tiny $k_2$};
\node at (-1.625,0) {\tiny $k_1{+}k_2$};
\end{tikzpicture}
};
\endxy
\otimes
\xy
(0,0)*{
\begin{tikzpicture} [scale=.3]
\draw[very thick] (0,-2.75) to (0,2.5);
\node at (0,-3.15) {\tiny $k_3$};
\node at (0,3) {\tiny $k_3$};
\end{tikzpicture}
}
\endxy
=
\xy
(-5,0)*{
\begin{tikzpicture}[scale=.3]
\draw [very thick] (0,-1) to (0,.75);
\draw [very thick] (0,.75) to [out=30,in=270] (1,2.5);
\draw [very thick] (0,.75) to [out=150,in=270] (-1,2.5); 
\draw [very thick] (1,-2.75) to [out=90,in=330] (0,-1);
\draw [very thick] (-1,-2.75) to [out=90,in=210] (0,-1);
\node at (-1,3) {\tiny $l_1$};
\node at (1,3) {\tiny $l_2$};
\node at (-1,-3.15) {\tiny $k_1$};
\node at (1,-3.15) {\tiny $k_2$};
\node at (-1.625,0) {\tiny $k_1{+}k_2$};
\end{tikzpicture}
};
(5,0)*{
\begin{tikzpicture} [scale=.3]
\draw[very thick] (0,-2.75) to (0,2.5);
\node at (0,-3.15) {\tiny $k_3$};
\node at (0,3) {\tiny $k_3$};
\end{tikzpicture}
}
\endxy
\]

\end{center}
where in the final equation $k_1+k_2=l_1+l_2$. The cup and cap diagrams will be defined later on. Moreover, here we have omitted the arrows, meaning that these conventions hold for all orientations (we elaborate on orientations below).\\

We want to consider a certain quotient of $\mathbf{fSWeb}_{\uparrow}(\mathfrak{gl}_n)$. We therefore need to define some other webs. We will later impose a (co)associativity relation on the webs, see \autoref{(co)assoc}. For now we proceed as follows. 

\begin{Definition}
The \textit{thin over-crossing} and \textit{thin under-crossing} are defined respectively as
\begin{gather}\label{thin-braid}
\begin{gathered}
\xy
(0,0)*{
\begin{tikzpicture}[scale=.3]
\draw [very thick, ->] (-1,-1) to (1,1);
\draw [very thick] (1,-1) to (0.25,-0.25);
\draw [very thick, ->] (-0.25,0.25) to (-1,1);
\node at (-1,-1.5) {\tiny $1$};
\node at (1,-1.5) {\tiny $1$};
\node at (1,1.45) {\tiny $1$};
\node at (-1,1.45) {\tiny $1$};
\end{tikzpicture}
};
\endxy  = -q^{-1-\frac{1}{n}} \Bigg(  
\xy
(0,0)*{
\begin{tikzpicture}[scale=.3]
\draw [very thick, directed=.55] (-1, -2.5) to (-1,2.5);
\draw [very thick, directed=.55] (1, -2.5) to (1,2.5);
\node at (-1,3) {\tiny $1$};
\node at (1,3) {\tiny $1$};
\node at (-1,-3) {\tiny $1$};
\node at (1,-3) {\tiny $1$};
\end{tikzpicture}
};
\endxy - q
\xy
(0,0)*{
\begin{tikzpicture}[scale=.3]
\draw [very thick, directed=.55] (0,-1) to (0,.75);
\draw [very thick, directed=.65] (0,.75) to [out=30,in=270] (1,2.5);
\draw [very thick, directed=.65] (0,.75) to [out=150,in=270] (-1,2.5); 
\draw [very thick, directed=.65] (-1, -2.5) to [out=90, in=210] (0,-1);
\draw [very thick, directed=.65] (1,-2.5) to [out =90, in=330] (0,-1);
\node at (-1,-3) {\tiny $1$};
\node at (1,-3) {\tiny $1$};
\node at (0.75,0) {\tiny $2$};
\node at (-1,3) {\tiny $1$};
\node at (1,3) {\tiny $1$};
\end{tikzpicture}
};
\endxy \Bigg),
\\
\xy
(0,0)*{
\begin{tikzpicture}[scale=.3]
\draw [very thick, ->] (1,-1) to (-1,1);
\draw [very thick] (-1,-1) to (-0.25,-0.25);
\draw [very thick, ->] (0.25,0.25) to (1,1);
\node at (-1,-1.5) {\tiny $1$};
\node at (1,-1.5) {\tiny $1$};
\node at (1,1.45) {\tiny $1$};
\node at (-1,1.45) {\tiny $1$};
\end{tikzpicture}
};
\endxy  = -q^{1+\frac{1}{n}} \Bigg( 
\xy
(0,0)*{
\begin{tikzpicture}[scale=.3]
\draw [very thick, directed=.55] (-1, -2.5) to (-1,2.5);
\draw [very thick, directed=.55] (1, -2.5) to (1,2.5);
\node at (-1,3) {\tiny $1$};
\node at (1,3) {\tiny $1$};
\node at (-1,-3) {\tiny $1$};
\node at (1,-3) {\tiny $1$};
\end{tikzpicture}
};
\endxy- q^{-1} 
\xy
(0,0)*{
\begin{tikzpicture}[scale=.3]
\draw [very thick, directed=.55] (0,-1) to (0,.75);
\draw [very thick, directed=.65] (0,.75) to [out=30,in=270] (1,2.5);
\draw [very thick, directed=.65] (0,.75) to [out=150,in=270] (-1,2.5); 
\draw [very thick, directed=.65] (-1, -2.5) to [out=90, in=210] (0,-1);
\draw [very thick, directed=.65] (1,-2.5) to [out =90, in=330] (0,-1);
\node at (-1,-3) {\tiny $1$};
\node at (1,-3) {\tiny $1$};
\node at (0.75,0) {\tiny $2$};
\node at (-1,3) {\tiny $1$};
\node at (1,3) {\tiny $1$};
\end{tikzpicture}
};
\endxy
\Bigg).
\end{gathered}
\end{gather}
\end{Definition}

\begin{Definition}\label{thick-braid} 
The \textit{thick over-crossing of type $(k,l)$} is defined by ``exploding'' edges:
\begin{equation*}
\xy
(0,0)*{
\begin{tikzpicture}[scale=.3]
\draw [very thick, ->] (2,-2) to (-2,2);
\draw [very thick, ->,crossline] (-2,-2) to (2,2);
\node at (-2,-2.5) {\tiny $k$};
\node at (2,-2.5) {\tiny $l$};
\node at (-2,2.5) {\tiny $l$};
\node at (2,2.5) {\tiny $k$};
\end{tikzpicture}
};
\endxy= \frac{1}{[k]!} \frac{1}{[l]!}
\xy
(0,0)*{
\begin{tikzpicture}[scale=.75]
\draw [very thick] (1,-1) to (.5,-.5);
\draw [very thick, ->] (-.5,.5) to (-1,1);
\draw [very thick] (.5,-.5) to [out=-165,in=-105] (-.5,.5);
\draw [very thick] (.5,-.5) to [out=-285,in=15] (-.5,.5);
\draw [very thick,crossline] (-1,-1) to (-.5,-.5);
\draw [very thick, ->,crossline] (.5,.5) to (1,1);
\draw [very thick,crossline] (-.5,-.5) to [out=105, in=165] (.5,.5);
\draw [very thick,crossline] (-.5,-.5) to [out=-15, in=285] (.5,.5);
\node at (-1,-1.25) {\tiny $k$};
\node at (1,-1.25) {\tiny $l$};
\node at (-.7,0) {\tiny $1$};
\node at (.7,0) {\tiny $1$};
\node at (-.15,-.15) {\tiny \rotatebox{-45}{$$}};
\node at (.15,.15) {\tiny \rotatebox{-45}{$$}};
\node at (-.2,.7) {\tiny $1$};
\node at (.2,-.8) {\tiny $1$};
\node at (0,0) {\tiny $\cdots$};
\node at (-1,1.25) {\tiny $l$};
\node at (1,1.25) {\tiny $k$};
\end{tikzpicture}
};
\endxy.
\end{equation*}
Similarly we define the thick under-crossing. By definition, if $k=0$ or $l=0$ then the crossing is the corresponding identity. 
\end{Definition}

\begin{Definition}\label{(1,k)-merge} 
We define the $(1,k)$-merges as 
\[\begin{tikzpicture}[anchorbase,scale=0.75]
\draw[very thick,directed=0.5] (0,0)
node[below]{$1$}to[out=90,in=205] (0.5,0.75);
\draw[very thick,directed=0.5] (1,0)
node[below]{$k$}to[out=90,in=-25] (0.5,0.75);
\draw[very thick,directed=0.5] (0.5,0.75)
to (0.5,1.5)node[above]{$k+1$};
\end{tikzpicture}
=
q^{\frac{k}{n}-k}
\begin{tikzpicture}[anchorbase,scale=0.75]
\draw[very thick,directed=0.5] (0,0)
to[out=90,in=205] (0.5,0.75);
\draw[very thick,directed=0.5] (1,0)
to[out=90,in=-25] (0.5,0.75);
\draw[very thick,directed=0.5] (0.5,0.75)
to (0.5,1.5)node[above]{$k+1$};
\draw[very thick] (1,-1) node[below]{$k$} to (0,0);
\draw[very thick,crossline] (0,-1)
node[below]{$1$} to (1,0);
\end{tikzpicture},\]
and similarly the $(1,k)$-splits.
\end{Definition}

The scalar in the above definition is needed for \autoref{splits&merges comp crossings} to work. Further, we will also need the following.

\begin{Definition}\label{thick merges} 
The \textit{thick merge} of type $(k,l)$ (or the $(k,l)$-merge) is defined as below
\[
\xy
(0,0)*{
\begin{tikzpicture}[scale=.5]
\draw [very thick, directed=.55] (-1,3) to (-1,6);
\draw [very thick, rdirected=.55] (-1,3) to (-3,0);
\draw [very thick,  rdirected=.55] (-1,3) to (1,0);
\node at (-1,6.25) {\tiny $k{+}l$};
\node at (-3,-.25) {\tiny $k$};
\node at (1,-.25) {\tiny $l$};
\end{tikzpicture}
};
\endxy=\frac{1}{[l]!}
\xy
(0,0)*{
\begin{tikzpicture}[scale=.5, xscale=-1]
\draw [very thick, directed=.55] (-1,3) to (-1,6);
\draw [very thick, rdirected=.55] (-1,3) to (-3,0);
\draw [very thick,  rdirected=.55] (-1,3) to (1,0);
\node at (-1,6.25) {\tiny $k{+}l$};
\node at (-3,-.25) {\tiny $l$};
\node at (1,-.25) {\tiny $k$};
\draw [very thick, directed=.55] (-2.55,.67) to[out=0,in=270] 
(-1.2,1.25) to[out=90,in=220] (-.8,2.73);
\node at (-1.5,1.5) {\tiny $\cdots$};
\node at (-1,1.5) {\tiny $1$};
\node at (-2.5,1.5) {\tiny $1$};
\end{tikzpicture}
};
\endxy
\]
and similarly the $(k,l)$-split.
\end{Definition}

The meticulous reader observes that \autoref{thick-braid} and \autoref{thick merges} are not well-defined since we have not specified the precise form of the webs. However, in \autoref{invertibility of digon} we will explain that this is not a problem, so the concerned reader might choose any form of these webs. Similarly for the other web categories that we define later.

\begin{Definition}\label{ladders} 
The $F^{(j)}$ and $E^{(j)}$-\textit{ladders} are given by
\begin{equation*}
\changed{F^{(j)}=
\xy
(0,0)*{
\begin{tikzpicture}[scale=.3]
\draw [very thick, directed=.55] (-2,-2) to (-2,0);
\draw [very thick, directed=.55] (-2,0) to (-2,2);
\draw [very thick, directed=.55] (2,-2) to (2,0);
\draw [very thick, directed=.55] (2,0) to (2,2);
\draw [very thick, directed=.55] (-2,0) to (2,0);
\node at (-2,-2.5) {\tiny $k$};
\node at (2,-2.5) {\tiny $l$};
\node at (-2,2.5) {\tiny $k{-}j$};
\node at (2,2.5) {\tiny $l{+}j$};
\node at (0,0.75) {\tiny $j$};
\end{tikzpicture}
};
\endxy=
\xy
(0,0)*{
\begin{tikzpicture}[scale=.3]
\draw [very thick, directed=.55] (1,.75) to (1,2);
\draw [very thick] (2,-1) to [out=90,in=320] (1,.75);
\draw [very thick, rdirected=0.65] (2,-1) to (2,-2);
\draw [very thick, directed=.55] (-1,-2) to (-1,-.25);
\draw [very thick] (-1,-0.25) to [out=150,in=270] (-2,1.5);
\draw [very thick, directed=0.35] (-2,1.5) to (-2,2);
\draw [very thick, directed=.55] (-1,-.25) to (1,.75);
\node at (-1,-2.5) {\tiny $k$};
\node at (2,-2.5) {\tiny $l$};
\node at (-2,2.5) {\tiny $k{-}j$};
\node at (1,2.5) {\tiny $l{+}j$};
\node at (0,.75) {\tiny $j$};
\end{tikzpicture}
};
\endxy\quad\text{and}\quad
E^{(j)}=
\xy
(0,0)*{
\begin{tikzpicture}[scale=.3]
\draw [very thick, directed=.55] (-2,-2) to (-2,0);
\draw [very thick, directed=.55] (-2,0) to (-2,2);
\draw [very thick, directed=.55] (2,-2) to (2,0);
\draw [very thick, directed=.55] (2,0) to (2,2);
\draw [very thick, rdirected=.55] (-2,0) to (2,0);
\node at (-2,-2.5) {\tiny $k$};
\node at (2,-2.5) {\tiny $l$};
\node at (-2,2.5) {\tiny $k{+}j$};
\node at (2,2.5) {\tiny $l{-}j$};
\node at (0,0.75) {\tiny $j$};
\end{tikzpicture}
};
\endxy=
\xy
(0,0)*{\reflectbox{
\begin{tikzpicture}[scale=.3]
\draw [very thick, directed=.55] (1,.75) to (1,2);
\draw [very thick] (2,-1) to [out=90,in=320] (1,.75);
\draw [very thick, rdirected=0.65] (2,-1) to (2,-2);
\draw [very thick, directed=.55] (-1,-2) to (-1,-.25);
\draw [very thick] (-1,-0.25) to [out=150,in=270] (-2,1.5);
\draw [very thick, directed=0.35] (-2,1.5) to (-2,2);
\draw [very thick, directed=.55] (-1,-.25) to (1,.75);
\node at (-1,-2.5) {\reflectbox{\tiny $l$}};
\node at (2,-2.5) {\reflectbox{\tiny $k$}};
\node at (-2,2.5) {\reflectbox{\tiny $l{-}j$}};
\node at (1,2.5) {\reflectbox{\tiny $k{+}j$}};
\node at (-0.25,.75) {\reflectbox{\tiny $j$}};
\end{tikzpicture}
}};
\endxy.}
\end{equation*}
\end{Definition}

Now we are ready to define the aforementioned quotient of $\mathbf{fSWeb}_{\uparrow}(\mathfrak{gl}_n)$.

\begin{Definition}\label{up-webs} 
The \textit{thin symmetric upward pointing category} of $\mathfrak{gl}_n$-webs $\mathbf{SWeb}_{\uparrow}(\mathfrak{gl}_n)$ is obtained from $\mathbf{fSWeb}_{\uparrow}(\mathfrak{gl}_n)$ by imposing the following \changed{three} relations on the morphisms:
\end{Definition}

\begin{flushleft}
Thin associativity and coassociativity
\end{flushleft} 
\begin{equation}\label{(co)assoc}
\xy,
(0,0)*{
\begin{tikzpicture}[scale=.3]
\draw [very thick, directed=.45] (0,.75) to [out=90,in=220] (1,2.5);
\draw [very thick, directed=.45] (1,-1) to [out=90,in=330] (0,.75);
\draw [very thick, directed=.45] (-1,-1) to [out=90,in=210] (0,.75);
\draw [very thick, directed=.45] (3,-1) to [out=90,in=330] (1,2.5);
\draw [very thick, directed=.45] (1,2.5) to (1,4.25);
\node at (-1,-1.5) {\tiny $h$};
\node at (1,-1.5) {\tiny $1$};
\node at (-1.375,1.5) {\tiny $h{+}1$};
\node at (3,-1.5) {\tiny $l$};
\node at (1,4.75) {\tiny $h{+}1{+}l$};
\end{tikzpicture}
};
\endxy=\xy
(0,0)*{
\begin{tikzpicture}[scale=.3]
\draw [very thick, directed=.45] (0,.75) to [out=90,in=340] (-1,2.5);
\draw [very thick, directed=.45] (-1,-1) to [out=90,in=210] (0,.75);
\draw [very thick, directed=.45] (1,-1) to [out=90,in=330] (0,.75);
\draw [very thick, directed=.45] (-3,-1) to [out=90,in=220] (-1,2.5);
\draw [very thick, directed=.45] (-1,2.5) to (-1,4.25);
\node at (1,-1.5) {\tiny $l$};
\node at (-1,-1.5) {\tiny $1$};
\node at (1.25,1.5) {\tiny $1{+}l$};
\node at (-3,-1.5) {\tiny $h$};
\node at (-1,4.75) {\tiny $h{+}1{+}l$};
\end{tikzpicture}
};
\endxy
\quad\text{and}\quad
\xy
(0,0)*{\rotatebox{180}{
\begin{tikzpicture}[scale=.3]
\draw [very thick, rdirected=.55] (0,.75) to [out=90,in=220] (1,2.5);
\draw [very thick, rdirected=.55] (1,-1) to [out=90,in=330] (0,.75);
\draw [very thick, rdirected=.55] (-1,-1) to [out=90,in=210] (0,.75);
\draw [very thick, rdirected=.55] (3,-1) to [out=90,in=330] (1,2.5);
\draw [very thick, rdirected=.55] (1,2.5) to (1,4.25);
\node at (-1,-1.5) {\rotatebox{180}{\tiny $l$}};
\node at (1,-1.5) {\rotatebox{180}{\tiny $1$}};
\node at (-1.325,1.5) {\rotatebox{180}{\tiny $1\! +\! l$}};
\node at (3,-1.5) {\rotatebox{180}{\tiny $h$}};
\node at (1,4.75) {\rotatebox{180}{\tiny $h\! +\! 1\! +\! l$}};
\end{tikzpicture}
}};
\endxy=
\xy
(0,0)*{\reflectbox{\rotatebox{180}{
\begin{tikzpicture}[scale=.3]
\draw [very thick, rdirected=.55] (0,.75) to [out=90,in=220] (1,2.5);
\draw [very thick, rdirected=.55] (1,-1) to [out=90,in=330] (0,.75);
\draw [very thick, rdirected=.55] (-1,-1) to [out=90,in=210] (0,.75);
\draw [very thick, rdirected=.55] (3,-1) to [out=90,in=330] (1,2.5);
\draw [very thick, rdirected=.55] (1,2.5) to (1,4.25);
\node at (-1,-1.5) {\reflectbox{\rotatebox{180}{\tiny $h$}}};
\node at (1,-1.5) {\reflectbox{\rotatebox{180}{\tiny $1$}}};
\node at (-1.325,1.5) {\reflectbox{\rotatebox{180}{\tiny $h\! +\! 1$}}};
\node at (3,-1.5) {\reflectbox{\rotatebox{180}{\tiny $l$}}};
\node at (1,4.75) {\reflectbox{\rotatebox{180}{\tiny $h\! +\! 1\! +\! l$}}};
\end{tikzpicture}
}}};
\endxy.
\end{equation}
\begin{flushleft}The thin square switch 
\end{flushleft}
\begin{equation}\label{thinsqure-sw}
\xy
(0,0)*{
\begin{tikzpicture}[scale=.3]
\draw [very thick, directed=.55] (-2,-4) to (-2,-2);
\draw [very thick, directed=1] (-2,-2) to (-2,0.25);
\draw [very thick, directed=.55] (2,-4) to (2,-2);
\draw [very thick, directed=1] (2,-2) to (2,0.25);
\draw [very thick, directed=.55] (-2,-2) to (2,-2);
\draw [very thick] (-2,0.25) to (-2,2);
\draw [very thick, directed=.55] (-2,2) to (-2,4);
\draw [very thick] (2,0.25) to (2,2);
\draw [very thick, directed=.55] (2,2) to (2,4);
\draw [very thick, rdirected=.55] (-2,2) to (2,2);
\node at (-2,-4.5) {\tiny $k$};
\node at (2,-4.5) {\tiny $l$};
\node at (-2,4.5) {\tiny $k$};
\node at (2,4.5) {\tiny $l$};
\node at (-3.5,0) {\tiny $k{-}1$};
\node at (3.5,0) {\tiny $l{+}1$};
\node at (0,-1.25) {\tiny $1$};
\node at (0,2.75) {\tiny $1$};
\end{tikzpicture}
};
\endxy=\xy
(0,0)*{
\begin{tikzpicture}[scale=.3]
\draw [very thick, directed=.55] (-2,-4) to (-2,-2);
\draw [very thick, directed=1] (-2,-2) to (-2,0.25);
\draw [very thick, directed=.55] (2,-4) to (2,-2);
\draw [very thick, directed=1] (2,-2) to (2,0.25);
\draw [very thick, rdirected=.55] (-2,-2) to (2,-2);
\draw [very thick] (-2,0.25) to (-2,2);
\draw [very thick, directed=.55] (-2,2) to (-2,4);
\draw [very thick] (2,0.25) to (2,2);
\draw [very thick, directed=.55] (2,2) to (2,4);
\draw [very thick, directed=.55] (-2,2) to (2,2);
\node at (-2,-4.5) {\tiny $k$};
\node at (2,-4.5) {\tiny $l$};
\node at (-2,4.5) {\tiny $k$};
\node at (2,4.5) {\tiny $l$};
\node at (-3.5,0) {\tiny $k{+}1$};
\node at (3.5,0) {\tiny $l{-}1$};
\node at (0,-1.25) {\tiny $1$};
\node at (0,2.75) {\tiny $1$};
\end{tikzpicture}
};
\endxy+[k-l]\xy
(0,0)*{
\begin{tikzpicture}[scale=.3]
\draw [very thick, directed=.55] (-2,-4) to (-2,4);
\draw [very thick, directed=.55] (2,-4) to (2,4);
\node at (-2,-4.5) {\tiny $k$};
\node at (2,-4.5) {\tiny $l$};
\node at (-2,4.5) {\tiny $k$};
\node at (2,4.5) {\tiny $l$};
\end{tikzpicture}
};
\endxy \ ,
\end{equation}
where in both defining relations we used the thick merges, splits and ladders from above.

\begin{flushleft}
\changed{Dumbbell relation,}
\end{flushleft}
\changed{which we will state around \autoref{Eq:RedGreen} below (it does not play a role until \autoref{S:Basis}), and which depends on $n$.}

\begin{Definition} 
The two types of compositions of merges and splits
\begin{equation*}
\xy
(-5,0)*{
\begin{tikzpicture}[scale=.3]
\draw [very thick, directed=.6] (0,-1) to (0,.75);
\draw [very thick, directed=.75] (0,.75) to [out=30,in=270] (1,2.5);
\draw [very thick, directed =.75] (0,.75) to [out=150,in=270] (-1,2.5); 
\draw [very thick, directed=-.5] (1,-2.75) to [out=90,in=330] (0,-1);
\draw [very thick, directed=-.5] (-1,-2.75) to [out=90,in=210] (0,-1);
\node at (-1,3) {\tiny $k$};
\node at (1,3) {\tiny $l$};
\node at (-1,-3.15) {\tiny $k$};
\node at (1,-3.15) {\tiny $l$};
\node at (-1.625,0) {\tiny $k{+}l$};
\end{tikzpicture}
};
\endxy \quad \text{and}  \quad 
\xy
(0,0)*{
\begin{tikzpicture}[scale=.3]
\draw [very thick, directed=.55] (0,.75) to (0,2.5);
\draw [very thick, directed=.55] (0,-2.75) to [out=30,in=330] (0,.75);
\draw [very thick, directed=.55] (0,-2.75) to [out=150,in=210] (0,.75);
\draw [very thick, directed=.55] (0,-4.5) to (0,-2.75);
\node at (0,-5) {\tiny $k{+}l$};
\node at (0,3) {\tiny $k{+}l$};
\node at (-2,-1) {\tiny $k$};
\node at (2,-1) {\tiny $l$};
\end{tikzpicture}
};
\endxy
\end{equation*}
are called the \textit{$(k,l)$-dumbbell} (or \textit{dumbbell of thickness $k+l$}) and the \textit{$(k,l)$-digon}, respectively.
\end{Definition}

Note that we also illustrate webs with potential zero or negative edges. In these cases, by convention, if a diagram has a 0-labelled edge then we remove that edge and if a diagram has a negative-labelled edge then the entire diagram is 0, e.g. the so-called \textit{thin digon removal} 
\begin{equation*}
\xy
(0,0)*{
\begin{tikzpicture}[scale=.3]
\draw [very thick, directed=.55] (0,.75) to (0,2.5);
\draw [very thick, directed=.55] (0,-2.75) to [out=30,in=330] (0,.75);
\draw [very thick, directed=.55] (0,-2.75) to [out=150,in=210] (0,.75);
\draw [very thick, directed=.55] (0,-4.5) to (0,-2.75);
\node at (0,-5) {\tiny $k$};
\node at (0,3) {\tiny $k$};
\node at (-2,-1) {\tiny $k{-}1$};
\node at (2,-1) {\tiny $1$};
\end{tikzpicture}
};
\endxy=[k] \xy
(0,0)*{
\begin{tikzpicture}[scale=.3]
\draw [very thick, directed=.55] (0,-4.5) to (0,2.5);
\node at (0,-5) {\tiny $k$};
\node at (0,3) {\tiny $k$};
\end{tikzpicture}
};
\endxy
\end{equation*}
is a special case of the thin square switches \autoref{thinsqure-sw}, namely in is the case for $l=0$. By \autoref{(co)assoc} the thin digon removal implies
\begin{equation*}
\xy
(0,0)*{
\begin{tikzpicture}[scale=.3]
\draw [very thick, directed=.55] (0,.75) to (0,2.5);
\draw [very thick, directed=.55] (0,-2.75) to [out=30,in=330] (0,.75);
\draw [very thick, directed=.55] (0,-2.75) to [out=150,in=210] (0,.75);
\draw [very thick, directed=.55] (0,-4.5) to (0,-2.75);
\node at (0,-5) {\tiny $k$};
\node at (0,3) {\tiny $k$};
\node at (-2,-1) {\tiny $1$};
\node at (2,-1) {\tiny $1$};
\node at (0,-1) {\tiny $\cdots$};
\end{tikzpicture}
};
\endxy=[k]! \xy
(0,0)*{
\begin{tikzpicture}[scale=.3]
\draw [very thick, directed=.55] (0,-4.5) to (0,2.5);
\node at (0,-5) {\tiny $k$};
\node at (0,3) {\tiny $k$};
\end{tikzpicture}
};
\endxy.
\end{equation*}

\begin{Remark}\label{invertibility of digon}
Note that since \ochanged{the quantum numbers are invertible in our ground field is $\Bbbk(q^\frac{1}{n})$}, the relation above is an invertible relation, namely we have 
\[\xy
(0,0)*{
\begin{tikzpicture}[scale=.3]
\draw [very thick, directed=.55] (0,-4.5) to (0,2.5);
\node at (0,-5) {\tiny $k$};
\node at (0,3) {\tiny $k$};
\end{tikzpicture}
};
\endxy = \frac{1}{[k]!} \xy
(0,0)*{
\begin{tikzpicture}[scale=.3]
\draw [very thick, directed=.55] (0,.75) to (0,2.5);
\draw [very thick, directed=.55] (0,-2.75) to [out=30,in=330] (0,.75);
\draw [very thick, directed=.55] (0,-2.75) to [out=150,in=210] (0,.75);
\draw [very thick, directed=.55] (0,-4.5) to (0,-2.75);
\node at (0,-5) {\tiny $k$};
\node at (0,3) {\tiny $k$};
\node at (-2,-1) {\tiny $1$};
\node at (2,-1) {\tiny $1$};
\node at (0,-1) {\tiny $\cdots$};
\end{tikzpicture}
};
\endxy ,\]
which justifies \autoref{thick-braid}.
Moreover, together with the associativity consequence from below
\[
\xy
(0,0)*{
\begin{tikzpicture}[scale=.5, xscale=-1]
\draw [very thick, directed=.55] (-1,3) to (-1,6);
\draw [very thick, rdirected=.55] (-1,3) to (-3,0);
\draw [very thick,  rdirected=.55] (-1,3) to (1,0);
\node at (-1,6.25) {\tiny $k{+}l$};
\node at (-3,-.25) {\tiny $l$};
\node at (1,-.25) {\tiny $k$};
\draw [very thick, directed=.55] (-2.55,.67) to[out=0,in=270] 
(-1.2,1.25) to[out=90,in=220] (-.8,2.73);
\node at (-1.5,1.5) {\tiny $\cdots$};
\node at (-1,1.5) {\tiny $1$};
\node at (-2.5,1.5) {\tiny $1$};
\end{tikzpicture}
};
\endxy=
\xy
(0,0)*{
\begin{tikzpicture}[scale=.6]
\draw [very thick,  directed=.65] (-1,3) to (-1,6);
\draw [very thick] (-1,3) to (-1.45,2.333);
\draw [very thick,directed=.55] (-3,0) to (-1,3);
\draw [very thick, ] (-1,3) to (-.8,2.73);
\draw [very thick, rdirected=.65] (-.8,2.73) to (-.55,2.333);
\draw [very thick, rdirected=.55] (0.55,.67) to (1,0);
\draw [very thick, directed=.55] (0.55,.67) to [out=187.7,in=243.7] 
(-.55,2.333);
\draw [very thick, directed=.55] (0.55,.67) to [out=63.7,in=363.7] 
(-.55,2.333);
\node at (-1,6.25) {\tiny $k{+}l$};
\node at (-3,-.25) {\tiny $k$};
\node at (1,-.25) {\tiny $l$};
\node at (-.4,2.8) {\tiny $l$};
\node at (0,1.5) {\tiny \rotatebox{33.7}{$\cdots$}};
\node at (-0.8,1.15) {\tiny $1$};
\node at (0.8,1.85) {\tiny $1$};
\end{tikzpicture}
};
\endxy \]  
it also justifies \autoref{thick merges}.

\changed{In more details, any way to interpret \autoref{thick-braid} and \autoref{thick merges} implies \autoref{thick (co)assoc}, which in turn implies that all ways to read \autoref{thick-braid} and \autoref{thick merges} are the same webs in the quotient $\mathbf{SWeb}_{\uparrow}(\mathfrak{gl}_n)$ of $\mathbf{fSWeb}_{\uparrow}(\mathfrak{gl}_n)$: different representatives are related by the already imposed local relations, hence define the same morphism in the quotient.}
\end{Remark}

We wanted to define our web category in a way that it is as thin as possible, i.e. with as few $k$-labelled edges as possible, where $k>1$.
We now show that in this thin-defined web category we can deduce thick versions of relations \autoref{(co)assoc} and \autoref{thinsqure-sw}.

\begin{Lemma}\label{thick (co)assoc} 
Thick (co)associativity 
\[
\xy,
(0,0)*{
\begin{tikzpicture}[scale=.3]
\draw [very thick, directed=.45] (0,.75) to [out=90,in=220] (1,2.5);
\draw [very thick, directed=.45] (1,-1) to [out=90,in=330] (0,.75);
\draw [very thick, directed=.45] (-1,-1) to [out=90,in=210] (0,.75);
\draw [very thick, directed=.45] (3,-1) to [out=90,in=330] (1,2.5);
\draw [very thick, directed=.45] (1,2.5) to (1,4.25);
\node at (-1,-1.5) {\tiny $h$};
\node at (1,-1.5) {\tiny $k$};
\node at (-1.375,1.5) {\tiny $h{+}k$};
\node at (3,-1.5) {\tiny $l$};
\node at (1,4.75) {\tiny $h{+}k{+}l$};
\end{tikzpicture}
};
\endxy=\xy
(0,0)*{
\begin{tikzpicture}[scale=.3]
\draw [very thick, directed=.45] (0,.75) to [out=90,in=340] (-1,2.5);
\draw [very thick, directed=.45] (-1,-1) to [out=90,in=210] (0,.75);
\draw [very thick, directed=.45] (1,-1) to [out=90,in=330] (0,.75);
\draw [very thick, directed=.45] (-3,-1) to [out=90,in=220] (-1,2.5);
\draw [very thick, directed=.45] (-1,2.5) to (-1,4.25);
\node at (1,-1.5) {\tiny $l$};
\node at (-1,-1.5) {\tiny $k$};
\node at (1.25,1.5) {\tiny $k{+}l$};
\node at (-3,-1.5) {\tiny $h$};
\node at (-1,4.75) {\tiny $h{+}k{+}l$};
\end{tikzpicture}
};
\endxy
\quad\text{and}\quad
\xy
(0,0)*{\rotatebox{180}{
\begin{tikzpicture}[scale=.3]
\draw [very thick, rdirected=.55] (0,.75) to [out=90,in=220] (1,2.5);
\draw [very thick, rdirected=.55] (1,-1) to [out=90,in=330] (0,.75);
\draw [very thick, rdirected=.55] (-1,-1) to [out=90,in=210] (0,.75);
\draw [very thick, rdirected=.55] (3,-1) to [out=90,in=330] (1,2.5);
\draw [very thick, rdirected=.55] (1,2.5) to (1,4.25);
\node at (-1,-1.5) {\rotatebox{180}{\tiny $l$}};
\node at (1,-1.5) {\rotatebox{180}{\tiny $k$}};
\node at (-1.325,1.5) {\rotatebox{180}{\tiny $k\! +\! l$}};
\node at (3,-1.5) {\rotatebox{180}{\tiny $h$}};
\node at (1,4.75) {\rotatebox{180}{\tiny $h\! +\! k\! +\! l$}};
\end{tikzpicture}
}};
\endxy=
\xy
(0,0)*{\reflectbox{\rotatebox{180}{
\begin{tikzpicture}[scale=.3]
\draw [very thick, rdirected=.55] (0,.75) to [out=90,in=220] (1,2.5);
\draw [very thick, rdirected=.55] (1,-1) to [out=90,in=330] (0,.75);
\draw [very thick, rdirected=.55] (-1,-1) to [out=90,in=210] (0,.75);
\draw [very thick, rdirected=.55] (3,-1) to [out=90,in=330] (1,2.5);
\draw [very thick, rdirected=.55] (1,2.5) to (1,4.25);
\node at (-1,-1.5) {\reflectbox{\rotatebox{180}{\tiny $h$}}};
\node at (1,-1.5) {\reflectbox{\rotatebox{180}{\tiny $k$}}};
\node at (-1.325,1.5) {\reflectbox{\rotatebox{180}{\tiny $h\! +\! k$}}};
\node at (3,-1.5) {\reflectbox{\rotatebox{180}{\tiny $l$}}};
\node at (1,4.75) {\reflectbox{\rotatebox{180}{\tiny $h\! +\! k\! +\! l$}}};
\end{tikzpicture}
}}};
\endxy
\]
follows from thin (co)associativity \autoref{(co)assoc} and thin square switch relation \autoref{thinsqure-sw}.
\end{Lemma}

\begin{proof} 
From the discussion above, we can explode the $k$-labelled edge and proceed like below
\[
\xy,
(0,0)*{
\begin{tikzpicture}[scale=.5]
\draw [very thick, directed=.45] (0,.75) to [out=90,in=220] (1,2.5);
\draw [very thick, directed=.45] (1,-1) to [out=90,in=330] (0,.75);
\draw [very thick, directed=.45] (-1,-1) to [out=90,in=210] (0,.75);
\draw [very thick, directed=.45] (3,-1) to [out=90,in=330] (1,2.5);
\draw [very thick, directed=.45] (1,2.5) to (1,4.25);
\node at (-1,-1.5) {\tiny $h$};
\node at (1,-1.5) {\tiny $k$};
\node at (-1.375,1.5) {\tiny $h{+}k$};
\node at (3,-1.5) {\tiny $l$};
\node at (1,4.75) {\tiny $h{+}k{+}l$};
\end{tikzpicture}
};
\endxy= \frac{1}{[k]!}
\xy,
(0,0)*{
\begin{tikzpicture}[scale=.5]
\draw [very thick, directed=.65] (-1,-1) to [out=90,in=210] (1,2.5);
\draw [very thick, directed=.65] (3,-1) to [out=90,in=330] (1,2.5);
\draw [very thick, directed=.45] (1,2.5) to (1,4.25);
\node at (-1,-1.5) {\tiny $h$};
\node at (1,-1.5) {\tiny $k$};
\node at (0,0) {\tiny $\cdots$};
\node at (-.5,2) {\tiny $h{+}k$};
\node at (3,-1.5) {\tiny $l$};
\node at (1,4.75) {\tiny $h{+}k{+}l$};
\node at (-.6,-.3) {\tiny $1$};
\node at (.7,0.2) {\tiny $1$};
\draw [very thick] (.5,-.5) to [out=-165,in=-105] (-.5,.5);
\draw [very thick] (.5,-.5) to [out=-285,in=15] (-.5,.5);
\draw [very thick] (1,-1) to (.5,-.5);
\draw [very thick] (-.5,.5) to (-.7,.7);
\end{tikzpicture}
};
\endxy= \frac{1}{[k]!}
\xy,
(0,0)*{
\begin{tikzpicture}[scale=.5]
\draw [very thick, directed=.65] (-1,-1) to [out=90,in=210] (1,2.5);
\draw [very thick, directed=.65] (3,-1) to [out=90,in=330] (1,2.5);
\draw [very thick, directed=.45] (1,2.5) to (1,4.25);
\draw [very thick] (1,-1) to (.5,-.5);
\node at (-1,-1.5) {\tiny $h$};
\node at (.9,-1.5) {\tiny $k$};
\node at (-.1,2.5) {\tiny $h{+}k$};
\node at (3,-1.5) {\tiny $l$};
\node at (1,4.75) {\tiny $h{+}k{+}l$};
\node at (0,-.8) {\tiny $1$};
\node at (.8,0) {\tiny $1$};
\node at (0,0) {\tiny $\cdots$};
\draw [very thick] (.5,-.5) to [out=225, in=-45] (-.2,-.4);
\draw [very thick] (.5,-.5) to [out=45, in=-45] (.4,.2);
\draw [very thick] (-.2,-.4) to (-.85,0.2);
\draw [very thick] (.4,.2) to (-.55,1);
\end{tikzpicture}
};
\endxy 
\]
\[
=\frac{1}{[k]!}
\xy,
(0,0)*{ 
\begin{tikzpicture}[scale=.5, xscale=-1]
\draw [very thick, directed=.65] (-1,-1) to [out=90,in=210] (1,2.5);
\draw [very thick, directed=.65] (3,-1) to [out=90,in=330] (1,2.5);
\draw [very thick, directed=.45] (1,2.5) to (1,4.25);
\draw [very thick] (1,-1) to (.5,-.5);
\node at (-1,-1.5) {\tiny $l$};
\node at (.9,-1.5) {\tiny $k$};
\node at (-.1,2.5) {\tiny $h{+}k$};
\node at (3,-1.5) {\tiny $h$};
\node at (1,4.75) {\tiny $h{+}k{+}l$};
\node at (0,-.9) {\tiny $1$};
\node at (.8,0) {\tiny $1$};
\node at (0,0) {\tiny $\cdots$};
\draw [very thick] (.5,-.5) to [out=225, in=-45] (-.2,-.4);
\draw [very thick] (.5,-.5) to [out=45, in=-45] (.4,.2);
\draw [very thick] (-.2,-.4) to (-.85,0.2);
\draw [very thick] (.4,.2) to (-.55,1);
\end{tikzpicture}
};
\endxy =
\frac{1}{[k]!}
\xy,
(0,0)*{
\begin{tikzpicture}[scale=.5, xscale=-1]
\draw [very thick, directed=.65] (-1,-1) to [out=90,in=210] (1,2.5);
\draw [very thick, directed=.65] (3,-1) to [out=90,in=330] (1,2.5);
\draw [very thick, directed=.45] (1,2.5) to (1,4.25);
\node at (-1,-1.5) {\tiny $l$};
\node at (1,-1.5) {\tiny $k$};
\node at (0,0) {\tiny $\cdots$};
\node at (-.5,2) {\tiny $h{+}k$};
\node at (3,-1.5) {\tiny $h$};
\node at (1,4.75) {\tiny $h{+}k{+}l$};
\node at (-.6,-.3) {\tiny $1$};
\node at (.7,0.2) {\tiny $1$};
\draw [very thick] (.5,-.5) to [out=-165,in=-105] (-.5,.5);
\draw [very thick] (.5,-.5) to [out=-285,in=15] (-.5,.5);
\draw [very thick] (1,-1) to (.5,-.5);
\draw [very thick] (-.5,.5) to (-.7,.7);
\end{tikzpicture}
};
\endxy=
\xy
(0,0)*{
\begin{tikzpicture}[scale=.5]
\draw [very thick, directed=.45] (0,.75) to [out=90,in=340] (-1,2.5);
\draw [very thick, directed=.45] (-1,-1) to [out=90,in=210] (0,.75);
\draw [very thick, directed=.45] (1,-1) to [out=90,in=330] (0,.75);
\draw [very thick, directed=.45] (-3,-1) to [out=90,in=220] (-1,2.5);
\draw [very thick, directed=.45] (-1,2.5) to (-1,4.25);
\node at (1,-1.5) {\tiny $l$};
\node at (-1,-1.5) {\tiny $k$};
\node at (1.25,1.5) {\tiny $k{+}l$};
\node at (-3,-1.5) {\tiny $h$};
\node at (-1,4.75) {\tiny $h{+}k{+}l$};
\end{tikzpicture}
};
\endxy,
\]
where in the second, third and fourth equality we used relation \autoref{(co)assoc} to pull thin edges along thick ones from one side to the other. Thick coassociativity is proved similarly.
\end{proof}

\begin{Lemma}\label{thick sq switch and div powers} 
The \textit{thick square switches} 
\begin{center}

\begin{equation*}\
\xy
(0,0)*{
\begin{tikzpicture}[scale=.3]
\draw [very thick, directed=.55] (-2,-4) to (-2,-2);
\draw [very thick, directed=1] (-2,-2) to (-2,0.25);
\draw [very thick, directed=.55] (2,-4) to (2,-2);
\draw [very thick, directed=1] (2,-2) to (2,0.25);
\draw [very thick, directed=.55] (-2,-2) to (2,-2);
\draw [very thick] (-2,0.25) to (-2,2);
\draw [very thick, directed=.55] (-2,2) to (-2,4);
\draw [very thick] (2,0.25) to (2,2);
\draw [very thick, directed=.55] (2,2) to (2,4);
\draw [very thick, rdirected=.55] (-2,2) to (2,2);
\node at (-2,-4.5) {\tiny $k$};
\node at (2,-4.5) {\tiny $l$};
\node at (-2,4.5) {\tiny $k{-}h{+}g$};
\node at (2,4.5) {\tiny $l{+}h{-}g$};
\node at (-3.5,0) {\tiny $k{-}h$};
\node at (3.5,0) {\tiny $l{+}h$};
\node at (0,-1.25) {\tiny $h$};
\node at (0,2.75) {\tiny $g$};
\end{tikzpicture}
};
\endxy = \sum_{i} \begin{bmatrix} k-l+g -h\\ i \end{bmatrix}
\xy
(0,0)*{
\begin{tikzpicture}[scale=.3]
\draw [very thick, directed=.55] (-2,-4) to (-2,-2);
\draw [very thick, directed=1] (-2,-2) to (-2,0.25);
\draw [very thick, directed=.55] (2,-4) to (2,-2);
\draw [very thick, directed=1] (2,-2) to (2,0.25);
\draw [very thick, rdirected=.55] (-2,-2) to (2,-2);
\draw [very thick] (-2,0.25) to (-2,2);
\draw [very thick, directed=.55] (-2,2) to (-2,4);
\draw [very thick] (2,0.25) to (2,2);
\draw [very thick, directed=.55] (2,2) to (2,4);
\draw [very thick, directed=.55] (-2,2) to (2,2);
\node at (-2,-4.5) {\tiny $k$};
\node at (2,-4.5) {\tiny $l$};
\node at (-2,4.5) {\tiny $k{-}h{+}g$};
\node at (2,4.5) {\tiny $l{-}g{+}h$};
\node at (-4,0) {\tiny $k{+}g{-}i$};
\node at (4,0) {\tiny $l{-}g{+}i$};
\node at (0,-1.25) {\tiny $g{-}i$};
\node at (0,2.75) {\tiny $h{-}i$};
\end{tikzpicture}
};
\endxy, 
\end{equation*}
\begin{equation*}
\xy
(0,0)*{
\begin{tikzpicture}[scale=.3]
\draw [very thick, directed=.55] (-2,-4) to (-2,-2);
\draw [very thick, directed=1] (-2,-2) to (-2,0.25);
\draw [very thick, directed=.55] (2,-4) to (2,-2);
\draw [very thick, directed=1] (2,-2) to (2,0.25);
\draw [very thick, rdirected=.55] (-2,-2) to (2,-2);
\draw [very thick] (-2,0.25) to (-2,2);
\draw [very thick, directed=.55] (-2,2) to (-2,4);
\draw [very thick] (2,0.25) to (2,2);
\draw [very thick, directed=.55] (2,2) to (2,4);
\draw [very thick, directed=.55] (-2,2) to (2,2);
\node at (-2,-4.5) {\tiny $k$};
\node at (2,-4.5) {\tiny $l$};
\node at (-2,4.5) {\tiny $k{+}h{-}g$};
\node at (2,4.5) {\tiny $l{-}h{+}g$};
\node at (-4,0) {\tiny $k{+}h$};
\node at (3.75,0) {\tiny $l{-}h$};
\node at (0,-1.25) {\tiny $h$};
\node at (0,2.75) {\tiny $g$};
\end{tikzpicture}
};
\endxy= \sum_{i} \begin{bmatrix} l-k+g-h \\ i \end{bmatrix}
\xy
(0,0)*{
\begin{tikzpicture}[scale=.3]
\draw [very thick, directed=.55] (-2,-4) to (-2,-2);
\draw [very thick, directed=1] (-2,-2) to (-2,0.25);
\draw [very thick, directed=.55] (2,-4) to (2,-2);
\draw [very thick, directed=1] (2,-2) to (2,0.25);
\draw [very thick, directed=.55] (-2,-2) to (2,-2);
\draw [very thick] (-2,0.25) to (-2,2);
\draw [very thick, directed=.55] (-2,2) to (-2,4);
\draw [very thick] (2,0.25) to (2,2);
\draw [very thick, directed=.55] (2,2) to (2,4);
\draw [very thick, rdirected=.55] (-2,2) to (2,2);
\node at (-2,-4.5) {\tiny $k$};
\node at (2,-4.5) {\tiny $l$};
\node at (-2,4.5) {\tiny $k{+}h{-}g$};
\node at (2,4.5) {\tiny $l{-}h{+}g$};
\node at (-4,0) {\tiny $k{-}g{+}i$};
\node at (4,0) {\tiny $l{+}g{-}i$};
\node at (0,-1.25) {\tiny $g{-}i$};
\node at (0,2.75) {\tiny $h{-}i$};
\end{tikzpicture}
};
\endxy ,
\end{equation*}

\end{center}
and \textit{divided powers collapsing} 
\begin{equation*}\
\xy
(0,0)*{
\begin{tikzpicture}[scale=.3]
\draw [very thick, directed=.55] (-2,-4) to (-2,-2);
\draw [very thick, directed=1] (-2,-2) to (-2,0.25);
\draw [very thick, directed=.55] (2,-4) to (2,-2);
\draw [very thick, directed=1] (2,-2) to (2,0.25);
\draw [very thick, rdirected=.55] (-2,-2) to (2,-2);
\draw [very thick] (-2,0.25) to (-2,2);
\draw [very thick, directed=.55] (-2,2) to (-2,4);
\draw [very thick] (2,0.25) to (2,2);
\draw [very thick, directed=.55] (2,2) to (2,4);
\draw [very thick, rdirected=.55] (-2,2) to (2,2);
\node at (-2,-4.5) {\tiny $k$};
\node at (2,-4.5) {\tiny $l$};
\node at (-2,4.5) {\tiny $k{+}h{+}g$};
\node at (2,4.5) {\tiny $l{-}g{-}h$};
\node at (-3.5,0) {\tiny $k{+}h$};
\node at (3.5,0) {\tiny $l{-}h$};
\node at (0,-1.25) {\tiny $h$};
\node at (0,2.75) {\tiny $g$};
\end{tikzpicture}
};
\endxy = \begin{bmatrix} h+g \\ g \end{bmatrix}
\xy
(0,0)*{
\begin{tikzpicture}[scale=.3]
\draw [very thick, directed=.55] (-2,-4) to (-2,-2);
\draw [very thick] (-2,-2) to (-2,0.25);
\draw [very thick, directed=.55] (2,-4) to (2,-2);
\draw [very thick] (2,-2) to (2,0.25);
\draw [very thick] (-2,0.25) to (-2,2);
\draw [very thick, directed=.55] (-2,2) to (-2,4);
\draw [very thick] (2,0.25) to (2,2);
\draw [very thick, directed=.55] (2,2) to (2,4);
\draw [very thick, rdirected=.55] (-2,0) to (2,0);
\node at (-2,-4.5) {\tiny $k$};
\node at (2,-4.5) {\tiny $l$};
\node at (-2,4.5) {\tiny $k{+}h{+}g$};
\node at (2,4.5) {\tiny $l{-}g{-}h$};
\node at (0,-1.25) {\tiny $h{+}g$};
\end{tikzpicture}
};
\endxy, 
\end{equation*}
\begin{equation*}\
\xy
(0,0)*{
\begin{tikzpicture}[scale=.3]
\draw [very thick, directed=.55] (-2,-4) to (-2,-2);
\draw [very thick, directed=1] (-2,-2) to (-2,0.25);
\draw [very thick, directed=.55] (2,-4) to (2,-2);
\draw [very thick, directed=1] (2,-2) to (2,0.25);
\draw [very thick, directed=.55] (-2,-2) to (2,-2);
\draw [very thick] (-2,0.25) to (-2,2);
\draw [very thick, directed=.55] (-2,2) to (-2,4);
\draw [very thick] (2,0.25) to (2,2);
\draw [very thick, directed=.55] (2,2) to (2,4);
\draw [very thick, directed=.55] (-2,2) to (2,2);
\node at (-2,-4.5) {\tiny $k$};
\node at (2,-4.5) {\tiny $l$};
\node at (-2,4.5) {\tiny $k{-}h{-}g$};
\node at (2,4.5) {\tiny $l{+}h{+}g$};
\node at (-3.5,0) {\tiny $k{-}h$};
\node at (3.5,0) {\tiny $l{+}h$};
\node at (0,-1.25) {\tiny $h$};
\node at (0,2.75) {\tiny $g$};
\end{tikzpicture}
};
\endxy = \begin{bmatrix} h+g \\ g \end{bmatrix}
\xy
(0,0)*{
\begin{tikzpicture}[scale=.3]
\draw [very thick, directed=.55] (-2,-4) to (-2,-2);
\draw [very thick] (-2,-2) to (-2,0.25);
\draw [very thick, directed=.55] (2,-4) to (2,-2);
\draw [very thick] (2,-2) to (2,0.25);
\draw [very thick] (-2,0.25) to (-2,2);
\draw [very thick, directed=.55] (-2,2) to (-2,4);
\draw [very thick] (2,0.25) to (2,2);
\draw [very thick, directed=.55] (2,2) to (2,4);
\draw [very thick, directed=.55] (-2,0) to (2,0);
\node at (-2,-4.5) {\tiny $k$};
\node at (2,-4.5) {\tiny $l$};
\node at (-2,4.5) {\tiny $k{-}h{-}g$};
\node at (2,4.5) {\tiny $l{+}g{+}h$};
\node at (0,-1.25) {\tiny $h{+}g$};
\end{tikzpicture}
};
\endxy ,
\end{equation*}
hold in $\mathbf{SWeb}_{\uparrow}(\mathfrak{gl}_n)$.
\end{Lemma}

\begin{proof} 
This is \cite[Remark 2.5]{SaTu-bcd-webs} and can be shown by using the definition of $\mathbf{SWeb}_{\uparrow}(\mathfrak{gl}_n)$ and \autoref{thick (co)assoc}. 
\end{proof}

As a special case of the lemma above, also the \textit{thick digon removal} holds, i.e.
\[
\xy
(0,0)*{
\begin{tikzpicture}[scale=.3]
\draw [very thick, directed=.55] (0,.75) to (0,2.5);
\draw [very thick, directed=.55] (0,-2.75) to [out=30,in=330] (0,.75);
\draw [very thick, directed=.55] (0,-2.75) to [out=150,in=210] (0,.75);
\draw [very thick, directed=.55] (0,-4.5) to (0,-2.75);
\node at (0,-5) {\tiny $k{+}l$};
\node at (0,3) {\tiny $k{+}l$};
\node at (-2,-1) {\tiny $k$};
\node at (2,-1) {\tiny $l$};
\end{tikzpicture}
};
\endxy=\begin{bmatrix}
k+l \\l
\end{bmatrix}
\xy
(0,0)*{
\begin{tikzpicture}[scale=.3]
\draw [very thick, directed=.55] (0,-4.5) to (0,2.5);
\node at (0,-5) {\tiny $k{+}l$};
\node at (0,3) {\tiny $k{+}l$};
\end{tikzpicture}
};
\endxy.
\]

\begin{Corollary}\label{serre relat for webs} 
Symmetric webs satisfy the so-called Serre relation
\[
\xy
(0,0)*{
\begin{tikzpicture}[scale=.3]
\draw [very thick, directed=.55] (-2,-4) to (-2,-2);
\draw [very thick, ] (-2,-2) to (-2,0.25);
\draw [very thick, directed=.55] (2,-4) to (2,-2);
\draw [very thick, ] (2,-2) to (2,0.25);
\draw [very thick, directed=.55] (2,-2) to (6,-2);
\draw [very thick] (-2,0.25) to (-2,2);
\draw [very thick, directed=.55] (-2,2) to (-2,4);
\draw [very thick] (2,0.25) to (2,2);
\draw [very thick, directed=.55] (2,2) to (2,4);
\draw [very thick, directed=.6] (6,-4) to (6,4);
\draw [very thick, directed=.55] (-2,0) to (2,0);
\draw [very thick, directed=.55] (-2,2) to (2,2);
\node at (-2,-4.5) {\tiny $k$};
\node at (2,-4.5) {\tiny $l$};
\node at (-2,4.5) {\tiny $k{-}2$};
\node at (2,4.5) {\tiny $l{+}1$};
\node at (0,.75) {\tiny $1$};
\node at (0,2.75) {\tiny $1$};
\node at (6,-4.5) {\tiny $m$};
\node at (6,4.5) {\tiny $m{+}1$};
\node at (4,-1.25) {\tiny $1$};
\end{tikzpicture}
};
\endxy -[2]
\xy
(0,0)*{
\begin{tikzpicture}[scale=.3]
\draw [very thick, directed=.55] (-2,-4) to (-2,-2);
\draw [very thick, directed=1] (-2,-2) to (-2,0.25);
\draw [very thick, directed=.55] (2,-4) to (2,-2);
\draw [very thick, ] (2,-2) to (2,0.25);
\draw [very thick, directed=.55] (-2,-2) to (2,-2);
\draw [very thick] (-2,0.25) to (-2,2);
\draw [very thick, directed=.55] (-2,2) to (-2,4);
\draw [very thick] (2,0.25) to (2,2);
\draw [very thick, directed=.55] (2,2) to (2,4);
\draw [very thick, directed=.6] (6,-4) to (6,4);
\draw [very thick, directed=.55] (2,0) to (6,0);
\draw [very thick, directed=.55] (-2,2) to (2,2);
\node at (-2,-4.5) {\tiny $k$};
\node at (2,-4.5) {\tiny $l$};
\node at (-2,4.5) {\tiny $k{-}2$};
\node at (2,4.5) {\tiny $l{+}1$};
\node at (0,-1.25) {\tiny $1$};
\node at (0,2.75) {\tiny $1$};
\node at (4,.75) {\tiny $1$};
\node at (6,-4.5) {\tiny $m$};
\node at (6,4.5) {\tiny $m{+}1$};
\end{tikzpicture}
};
\endxy +
\xy
(0,0)*{
\begin{tikzpicture}[scale=.3]
\draw [very thick, directed=.55] (-2,-4) to (-2,-2);
\draw [very thick, ] (-2,-2) to (-2,0.25);
\draw [very thick, directed=.55] (2,-4) to (2,-2);
\draw [very thick, ] (2,-2) to (2,0.25);
\draw [very thick, directed=.55] (-2,-2) to (2,-2);
\draw [very thick] (-2,0.25) to (-2,2);
\draw [very thick, directed=.55] (-2,2) to (-2,4);
\draw [very thick] (2,0.25) to (2,2);
\draw [very thick, directed=.55] (2,2) to (2,4);
\draw [very thick, directed=.6] (6,-4) to (6,4);
\draw [very thick, directed=.55] (-2,0) to (2,0);
\draw [very thick, directed=.55] (2,2) to (6,2);
\node at (-2,-4.5) {\tiny $k$};
\node at (2,-4.5) {\tiny $l$};
\node at (-2,4.5) {\tiny $k{-}2$};
\node at (2,4.5) {\tiny $l{+}1$};
\node at (0,-1.25) {\tiny $1$};
\node at (0,0.75) {\tiny $1$};
\node at (4,2.75) {\tiny $1$};
\node at (6,-4.5) {\tiny $m$};
\node at (6,4.5) {\tiny $m{+}1$};
\end{tikzpicture}
};
\endxy=0.
\]
\end{Corollary}

\begin{proof}
The argument is as in \cite[Lemma 2.2.1]{CaKaMo-webs-skew-howe} and follows from applying associativity and the thick square switches to the middle web. Note that the divided powers from the leftmost and rightmost webs collapse, using the previous lemma.
\end{proof}

\begin{Remark}
The so-called higher order Serre relations 
as \ochanged{e.g. in} \cite[Chapter 7]{Lu-quantum-book}
hold in symmetric webs as well.
\end{Remark}

Further, we want to emphasize that our thin crossings from \autoref{thin-braid} satisfy the braiding relations.

\begin{Lemma}\label{thin R2&R3}  
The (thin) Reidemeister 2 and Reidemeister 3 moves hold in $\mathbf{SWeb}_{\uparrow}(\mathfrak{gl}_n)$, i.e.
\[
\xy
(0,0)*{
\begin{tikzpicture}[scale=.3]
\draw [very thick, ->] (-1,-1) to (1,1);
\draw [very thick] (1,-1) to (0.25,-0.25);
\draw [very thick, ->] (-0.25,0.25) to (-1,1);
\node at (-1,-3.5) {\tiny $1$};
\node at (1,-3.5) {\tiny $1$};
\node at (1,1.45) {\tiny $1$};
\node at (-1,1.45) {\tiny $1$};
\draw [very thick] (1,-3) to (-1,-1);
\draw [very thick] (-1,-3) to (-.25,-2.25);
\draw [very thick] (0.25,-1.75) to (1,-1);
\end{tikzpicture}
};
\endxy =
\xy
(0,0)*{
\begin{tikzpicture}[scale=.3]
\draw [very thick, directed=.55] (-1, -2.25) to (-1,2.25);
\draw [very thick, directed=.55] (1, -2.25) to (1,2.25);
\node at (-1,2.75) {\tiny $1$};
\node at (1,2.75) {\tiny $1$};
\node at (-1,-2.75) {\tiny $1$};
\node at (1,-2.75) {\tiny $1$};
\end{tikzpicture}
};
\endxy=
\xy
(0,0)*{
\begin{tikzpicture}[scale=.3]
\draw [very thick, ->] (1,-1) to (-1,1);
\draw [very thick] (-1,-1) to (-0.25,-0.25);
\draw [very thick, ->] (0.25,0.25) to (1,1);
\node at (-1,-3.5) {\tiny $1$};
\node at (1,-3.5) {\tiny $1$};
\node at (1,1.45) {\tiny $1$};
\node at (-1,1.45) {\tiny $1$};
\draw [very thick] (-1,-3) to (1,-1);
\draw [very thick] (1,-3) to (0.25,-2.25);
\draw [very thick] (-0.25,-1.75) to (-1,-1);
\end{tikzpicture}
};
\endxy,
\]
\[
\xy
(0,0)*{
\begin{tikzpicture}[scale=.3]
\draw [very thick, ->] (-1,-1) to (1,1);
\draw [very thick] (1,-1) to (0.25,-0.25);
\draw [very thick, ->] (-0.25,0.25) to (-1,1);
\draw [very thick,->] (3,-1) to (3,1);
\draw [very thick] (-1,-3) to (-1,-1);
\draw [very thick] (1,-1) to (3,-3);
\draw [very thick, crossline] (1,-3) to (3,-1);
\draw [very thick] (-1,-3) to (1,-5);
\draw [very thick,crossline] (-1,-5) to (1,-3);
\draw [very thick] (3,-5) to (3,-3);
\node at (-1,-5.5) {\tiny $1$};
\node at (1,-5.5) {\tiny $1$};
\node at (3,-5.5) {\tiny $1$};
\node at (-1,1.5) {\tiny $1$};
\node at (1,1.5) {\tiny $1$};
\node at (3,1.5) {\tiny $1$};
\end{tikzpicture}
};
\endxy=
\xy
(0,0)*{
\begin{tikzpicture}[scale=.3]
\draw [very thick, ->] (-1,-1) to (1,1);
\draw [very thick] (1,-1) to (0.25,-0.25);
\draw [very thick, ->] (-0.25,0.25) to (-1,1);
\draw [very thick,->] (-3,-1) to (-3,1);
\draw [very thick] (-1,-3) to (-3,-1);
\draw [very thick,crossline] (-3,-3) to (-1,-1);
\draw [very thick] (1,-3) to (1,-1);
\draw [very thick] (-3,-5) to (-3,-3);
\draw [very thick] (1,-5) to (-1,-3);
\draw [very thick,crossline] (-1,-5) to (1,-3);
\node at (-3,1.5) {\tiny $1$};
\node at (-1,1.5) {\tiny $1$};
\node at (1,1.5) {\tiny $1$};
\node at (-3,-5.5) {\tiny $1$};
\node at (1,-5.5) {\tiny $1$};
\node at (-1,-5.5) {\tiny $1$};
\end{tikzpicture}
};
\endxy.
\]
\end{Lemma}

\begin{proof}
This can be verified via a straightforward calculation. 
\end{proof}

We will see later that also the thick version of the previous lemma holds, meaning that $\mathbf{SWeb}_{\uparrow}(\mathfrak{gl}_n)$ is a braided monoidal category. For now, note that as a direct consequence of the \autoref{thin R2&R3} we have.

\begin{Corollary}
Thin over and under-crossings are inverses of each other.\qed
\end{Corollary}

Since in $\text{U}_{q}(\mathfrak{gl}_n)$-\textbf{fdMod} every object has a dual and we have evaluation and coevaluation maps, we want incorporate the analogous notions into our diagrammatic category. Therefore, we define the following.

\begin{Definition}\label{free symm-webs} 
The \textit{free thin symmetric (upward-downward pointing) category} of $\mathfrak{gl}_n$-webs, denote it 
by $\mathbf{fSWeb}_{\uparrow,\downarrow}(\mathfrak{gl}_n)$, is the $\Bbbk(q^{\frac{1}{n}})$-linear category monoidally generated by\\
$\bullet$ objects $ \{k_{\uparrow} , \  k_{\downarrow}  \ | \ k\in \mathbb{Z}_{\geq 0} \} $ with $0_\uparrow=0_\downarrow=0$,\\
$\bullet$ morphisms

\[
\xy
(0,0)*{
\begin{tikzpicture}[scale=.3]
\draw [very thick, directed=.55] (0, .75) to (0,2.5);
\draw [very thick, directed=.45] (1,-1) to [out=90,in=330] (0,.75);
\draw [very thick, directed=.45] (-1,-1) to [out=90,in=210] (0,.75); 
\node at (0, 3) {\tiny $k{+}1$};
\node at (-1,-1.5) {\tiny $k$};
\node at (1,-1.5) {\tiny $1$};
\end{tikzpicture}
};
\endxy : k_{\uparrow} \otimes 1_{\uparrow} \rightarrow (k+1)_{\uparrow}, \quad
\xy
(0,0)*{
\begin{tikzpicture}[scale=.3]
\draw [very thick, directed=.55] (0,-1) to (0,.75);
\draw [very thick, directed=.65] (0,.75) to [out=30,in=270] (1,2.5);
\draw [very thick, directed=.65] (0,.75) to [out=150,in=270] (-1,2.5); 
\node at (0, -1.5) {\tiny $k{+}1$};
\node at (-1,3) {\tiny $k$};
\node at (1,3) {\tiny $1$};
\end{tikzpicture}
};
\endxy: (k+1)_{\uparrow} \rightarrow k_{\uparrow} \otimes 1_{\uparrow}, \quad
\]
\[
\xy
(0,0)*{
\begin{tikzpicture}[scale=.3]
\draw [very thick, rdirected=.55] (0, .75) to (0,2.5);
\draw [very thick, rdirected=.45] (1,-1) to [out=90,in=330] (0,.75);
\draw [very thick, rdirected=.45] (-1,-1) to [out=90,in=210] (0,.75); 
\node at (0, 3) {\tiny $k{+}1$};
\node at (-1,-1.5) {\tiny $k$};
\node at (1,-1.5) {\tiny $1$};
\end{tikzpicture}
};
\endxy : k_{\downarrow} \otimes 1_{\downarrow} \rightarrow (k+1)_{\downarrow}, \quad
\xy
(0,0)*{
\begin{tikzpicture}[scale=.3]
\draw [very thick, rdirected=.55] (0,-1) to (0,.75);
\draw [very thick, rdirected=.65] (0,.75) to [out=30,in=270] (1,2.5);
\draw [very thick, rdirected=.65] (0,.75) to [out=150,in=270] (-1,2.5); 
\node at (0, -1.5) {\tiny $k{+}1$};
\node at (-1,3) {\tiny $k$};
\node at (1,3) {\tiny $1$};
\end{tikzpicture}
};
\endxy: (k+1)_{\downarrow} \rightarrow k_{\downarrow} \otimes 1_{\downarrow}, \quad
\]
\begin{equation*}
\xy
(0,0)*{
\begin{tikzpicture} [scale=1]
\draw[very thick,rdirected=-.95] (0,1) to [out=270,in=180] (.25,.5) to [out=0,in=270] (.5,1);
\node at (.25,.25) {};
\node at (0,1.25) {\tiny $1$};
\node at (.5,1.25) {\tiny $1$};
\end{tikzpicture}
}
\endxy : 0 \rightarrow 1_{\uparrow} \otimes 1_{\downarrow}, \quad 
\xy
(0,0)*{
\begin{tikzpicture} [scale=1]
\draw[very thick,rdirected=-.95] (0,0) to [out=90,in=180] (.25,.5) to [out=0,in=90] (.5,0);
\node at (.25,.75) {};
\node at (0,-.25) {\tiny $1$};
\node at (.5,-.25) {\tiny $1$};
\end{tikzpicture}
}
\endxy : 1_{\downarrow} \otimes 1_{\uparrow} \rightarrow 0 . 
\]
\end{Definition}
As before, we want to consider a quotient of $\mathbf{fSWeb}_{\uparrow,\downarrow}(\mathfrak{gl}_n)$ by some relations on the morphisms and in order to give these relations we need to define some other webs. First, we define upward pointing thin crossings as in \autoref{thin-braid}, thick upward crossings as in \autoref{thick-braid}, the $(1,k)$- merges and splits as in \autoref{(1,k)-merge} and also the thick merges and splits as in \autoref{thick merges}. Then we continue as follows.

\begin{Definition}\label{mates} 
We define the \changed{following}:\\
\begin{itemize}

\item \textit{Leftward (thin) crossings} 
\[
\xy
(0,0)*{
\begin{tikzpicture}[scale=.3]
\draw [very thick, ->] (1,-1) to (-1,1);
\draw [very thick] (1,1) to (0.25,0.25);
\draw [very thick, ->] (-0.25,-0.25) to (-1,-1);
\node at (-1,-1.5) {\tiny $1$};
\node at (1,-1.5) {\tiny $1$};
\node at (-1,1.5) {\tiny $1$};
\node at (1,1.5) {\tiny $1$};
\end{tikzpicture}
};
\endxy : =
\xy   
(0,0)*{
\begin{tikzpicture}[scale=.3]
\draw[very thick] (-2.5,1) to [out=90,in=180] (-1.75,2.5) to [out=0,in=90] (-1,1);
\draw[very thick] (1,-1) to [out=270,in=180] (1.75,-2.5) to [out=0,in=270] (2.5,-1);
\draw [very thick] (-1,-1) to (1,1);
\draw [very thick] (1,-1) to (0.25,-0.25);
\draw [very thick] (-0.25, 0.25) to (-1,1);
\draw[very thick] (-1,-1) to (-1,-2.5);
\draw[very thick, ->] (1,1) to (1,2.5);
\draw[very thick] (2.5,-1) to (2.5, 2.5);
\draw[very thick, ->] (-2.5, 1) to (-2.5,-2.5);
\node at (-2.5,-3) {\tiny $1$};
\node at (2.5,3) {\tiny $1$};
\node at (-1,-3) {\tiny $1$};
\node at (1,3) {\tiny $1$};
\end{tikzpicture}
};
\endxy \quad \text{and} \quad
\xy
(0,0)*{
\begin{tikzpicture}[scale=.3]
\draw [very thick, ->] (1,1) to (-1,-1);
\draw [very thick, ->] (-0.25,0.25) to (-1,1);
\draw [very thick] (0.25,-0.25) to (1,-1);
\node at (-1,-1.5) {\tiny $1$};
\node at (1,-1.5) {\tiny $1$};
\node at (-1,1.5) {\tiny $1$};
\node at (1,1.5) {\tiny $1$};
\end{tikzpicture}
};
\endxy := \xy   
(0,0)*{
\begin{tikzpicture}[scale=.3]
\draw[very thick] (-2.5,1) to [out=90,in=180] (-1.75,2.5) to [out=0,in=90] (-1,1);
\draw[very thick] (1,-1) to [out=270,in=180] (1.75,-2.5) to [out=0,in=270] (2.5,-1);
\draw [very thick] (1,-1) to (-1,1);
\draw [very thick] (-1,-1) to (-0.25,-0.25);
\draw [very thick] (0.25, 0.25) to (1,1);
\draw[very thick] (-1,-1) to (-1,-2.5);
\draw[very thick, ->] (1,1) to (1,2.5);
\draw[very thick] (2.5,-1) to (2.5, 2.5);
\draw[very thick, ->] (-2.5, 1) to (-2.5,-2.5);
\node at (-2.5,-3) {\tiny $1$};
\node at (2.5,3) {\tiny $1$};
\node at (1,3) {\tiny $1$};
\node at (-1,-3) {\tiny $1$};
\end{tikzpicture}
};
\endxy ;
\]

\item \textit{downward pointing (thin) over and under-crossings} 
\[
\xy
(0,0)*{
\begin{tikzpicture} [scale=1]
\draw [very thick, ->]  (.5,1.) to (0,.5);
\draw [very thick] (0,1) to (.2,.8);
\draw [very thick, ->] (.31,.7) to (.5,.5);
\node at (0,.25) {\tiny $1$};
\node at (0.5,.25) {\tiny $1$};
\node at (0,1.2) {\tiny $1$};
\node at (0.5,1.2) {\tiny $1$};
\end{tikzpicture}
}
\endxy:=
\xy
(0,0)*{
\begin{tikzpicture} [scale=1]
\draw [very thick]  (.5,1.) to (0,.5);
\draw [very thick] (0,1) to (.2,.8);
\draw [very thick] (.31,.7) to (.5,.5);
\draw [very thick] (0,1) to [in=0, out=90] (-.25,1.25) to [in=90, out=180] (-.5,1);
\draw [very thick] (.5,1) to [in=0, out=90] (-0.25,1.75) to [in=90, out=180] (-1,1);
\draw [very thick] (1,.5) to [in=0, out=270] (.75,.25) to [in=270, out=180] (.5,.5);
\draw [very thick] (1.5,.5) to [in=0, out=270] (.75,-.25) to [in=270, out=180] (0,.5);
\draw [very thick, ->](-.5,1) to (-.5,-.25);
\draw [very thick, ->](-1,1) to (-1,-.25);
\draw [very thick](1,2) to (1,.5);
\draw [very thick](1.5,2) to (1.5,.5);
\node at (-1,-.5) {\tiny $1$};
\node at (-.5,-.5) {\tiny $1$};
\node at (1,2.25) {\tiny $1$};
\node at (1.5,2.25) {\tiny $1$};
\end{tikzpicture}
}
\endxy \quad , \quad 
\xy
(0,0)*{
\begin{tikzpicture} [scale=1]
\draw [very thick]  (.5,1.) to (0.3,.8);
\draw [very thick,->] (0,1) to (.5,.5);
\draw [very thick, ->] (.2,.7) to (0,.5);
\node at (0,1.2) {\tiny $1$};
\node at (0.45,1.2) {\tiny $1$};
\node at (0,.25) {\tiny $1$};
\node at (0.5,.25) {\tiny $1$};
\end{tikzpicture}
}
\endxy:=
\xy
(0,0)*{
\begin{tikzpicture} [scale=1]
\draw [very thick]  (.5,1.) to (0.3,.8);
\draw [very thick] (0,1) to (.5,.5);
\draw [very thick] (.2,.7) to (0,.5);
\draw [very thick] (0,1) to [in=0, out=90] (-.25,1.25) to [in=90, out=180] (-.5,1);
\draw [very thick] (.5,1) to [in=0, out=90] (-0.25,1.75) to [in=90, out=180] (-1,1);
\draw [very thick] (1,.5) to [in=0, out=270] (.75,.25) to [in=270, out=180] (.5,.5);
\draw [very thick] (1.5,.5) to [in=0, out=270] (.75,-.25) to [in=270, out=180] (0,.5);
\draw [very thick, ->](-.5,1) to (-.5,-.25);
\draw [very thick, ->](-1,1) to (-1,-.25);
\draw [very thick](1,2) to (1,.5);
\draw [very thick](1.5,2) to (1.5,.5);
\node at (-1,-.5) {\tiny $1$};
\node at (-.5,-.5) {\tiny $1$};
\node at (1,2.25) {\tiny $1$};
\node at (1.5,2.25) {\tiny $1$};
\end{tikzpicture}
}
\endxy.
\]

\end{itemize}
\end{Definition}

\begin{Definition}\label{thick cap} 
The \textit{$k$-labelled caps} are defined via explosion
\[
\xy
(0,0)*{
\begin{tikzpicture} [scale=1]
\draw[very thick, rdirected=-.9] (0,0) to [out=90,in=180] (.25,.5) to [out=0,in=90] (.5,0);
\node at (.25,.75) {};
\node at (0,-.15) {\tiny $k$};
\node at (.5,-.15) {\tiny $k$};
\end{tikzpicture}
}
\endxy
=
\frac{1}{[k]!}
\xy
(0,0)*{
\begin{tikzpicture}[scale=.3]
\draw [very thick, rdirected=-.5] (0,-1) to (0,.75);
\draw [very thick] (0,.75) to [out=30,in=270] (1,2.5);
\draw [very thick] (0,.75) to [out=150,in=270] (-1,2.5); 
\draw [very thick,rdirected=-1] (1,2.5) to [out=90,in=180] (2,3.5) to [out=0,in=90] (3,2.5);
\draw [very thick,rdirected=-1] (-1,2.5) to [out=90,in=180] (2,5.5) to [out=0,in=90] (5,2.5);
\draw [very thick] (4,-1) to (4,.75);
\draw [very thick] (4,-1) to (4,.75);
\draw [very thick] (4,.75) to [out=30,in=270] (5,2.5);
\draw [very thick] (4,.75) to [out=150,in=270] (3,2.5);
\node at (0,-1.35) {\tiny $k$};
\node at (4,-1.35) {\tiny $k$};
\node at (2,3) {\tiny $1$};
\node at (2,6) {\tiny $1$};
\node at (2,4.75) {\tiny $\vdots$};
\end{tikzpicture}
};
\endxy
\]
and similarly the \textit{$k$-labelled cups}.
\end{Definition}

We will later impose that the (thin) leftward crossings are invertible and we will draw their inverses as (thin) \textit{rightward crossings}, i.e.
\begin{equation*}
{\bigg(\xy
(0,0)*{
\begin{tikzpicture}[scale=.3]
\draw [very thick, ->] (1,-1) to (-1,1);
\draw [very thick] (1,1) to (0.25,0.25);
\draw [very thick, ->] (-0.25,-0.25) to (-1,-1);
\node at (-1,-1.5) {\tiny $1$};
\node at (1,-1.5) {\tiny $1$};
\node at (-1,1.5) {\tiny $1$};
\node at (1,1.5) {\tiny $1$};
\end{tikzpicture}
};
\endxy\bigg)}^{-1} : = 
\xy
(0,0)*{
\begin{tikzpicture}[scale=.3]
\draw [very thick, ->] (-1,-1) to (1,1);
\draw [very thick] (-1,1) to (-0.25,0.25);
\draw [very thick, ->] (0.25,-0.25) to (1,-1);
\node at (-1,-1.5) {\tiny $1$};
\node at (1,-1.5) {\tiny $1$};
\node at (-1,1.5) {\tiny $1$};
\node at (1,1.5) {\tiny $1$};
\end{tikzpicture}
};
\endxy  \quad \text{and} \quad
{\bigg(\xy
(0,0)*{
\begin{tikzpicture}[scale=.3]
\draw [very thick, ->] (1,1) to (-1,-1);
\draw [very thick, ->] (-0.25,0.25) to (-1,1);
\draw [very thick] (0.25,-0.25) to (1,-1);
\node at (-1,-1.5) {\tiny $1$};
\node at (1,-1.5) {\tiny $1$};
\node at (-1,1.5) {\tiny $1$};
\node at (1,1.5) {\tiny $1$};
\end{tikzpicture}
};
\endxy \bigg)}^{-1} :=
\xy
(0,0)*{
\begin{tikzpicture}[scale=.3]
\draw [very thick, ->] (-1,1) to (1,-1);
\draw [very thick, ->] (0.25,0.25) to (1,1);
\draw [very thick] (-1,-1) to (-0.25,-.25);
\node at (-1,-1.5) {\tiny $1$};
\node at (1,-1.5) {\tiny $1$};
\node at (-1,1.5) {\tiny $1$};
\node at (1,1.5) {\tiny $1$};
\end{tikzpicture}
};
\endxy.
\end{equation*}

\begin{Definition}\label{thick leftward cross} 
The \textit{thick leftward} crossings are defined by explosion
\begin{equation*}
\xy
(0,0)*{
\begin{tikzpicture}[scale=.3]
\draw [very thick, ->] (2,-2) to (-2,2);
\draw [very thick, rdirected=0.05,crossline] (-2,-2) to (2,2);
\node at (-2,-2.5) {\tiny $k$};
\node at (2,-2.5) {\tiny $l$};
\node at (-2,2.5) {\tiny $l$};
\node at (2,2.5) {\tiny $k$};
\end{tikzpicture}
};
\endxy= \frac{1}{[k]!} \frac{1}{[l]!}
\xy
(0,0)*{
\begin{tikzpicture}[scale=.75]
\draw [very thick] (1,-1) to (.5,-.5);
\draw [very thick, ->] (-.5,.5) to (-1,1);
\draw [very thick] (.5,-.5) to [out=-165,in=-105] (-.5,.5);
\draw [very thick] (.5,-.5) to [out=-285,in=15] (-.5,.5);
\draw [very thick,crossline,rdirected=0.14] (-1,-1) to (-.5,-.5);
\draw [very thick, ,crossline] (.5,.5) to (1,1);
\draw [very thick,crossline] (-.5,-.5) to [out=105, in=165] (.5,.5);
\draw [very thick,crossline] (-.5,-.5) to [out=-15, in=285] (.5,.5);
\node at (-1,-1.25) {\tiny $k$};
\node at (1,-1.25) {\tiny $l$};
\node at (-.7,0) {\tiny $1$};
\node at (.7,0) {\tiny $1$};
\node at (-.15,-.15) {\tiny \rotatebox{-45}{$$}};
\node at (.15,.15) {\tiny \rotatebox{-45}{$$}};
\node at (-.2,.7) {\tiny $1$};
\node at (.2,-.8) {\tiny $1$};
\node at (0,0) {\tiny $\cdots$};
\node at (-1,1.25) {\tiny $l$};
\node at (1,1.25) {\tiny $k$};
\end{tikzpicture}
};
\endxy.
\end{equation*}
and similarly the \textit{thick rightward} crossings.
\end{Definition}

\begin{Definition}\label{right cap&cup} 
The \textit{rightward cap and cup} are defined as
\[
\xy
(0,0)*{
\begin{tikzpicture} [scale=1]
\draw[very thick,directed=1] (0,.75) to [out=90,in=180] (.25,1.25) to [out=0,in=90] (.5,.75);
\node at (0, .5) {\tiny $1$};
\node at (.5, .5) {\tiny $1$};
\end{tikzpicture}
}
\endxy := q^{n-\frac{1}{n}} 
\xy
(0,0)*{
\begin{tikzpicture} [scale=1]
\draw [very thick, ->] (0,.5) to (.5,1.);
\draw [very thick] (0,1) to (.2,.8);
\draw [very thick, ->] (.3,.65) to (.5,.5);
\draw [very thick, directed = 1] (.5,1) to [in=0, out=90] (.25,1.5) to [in=90, out=180] (0,1);
\node at (.25,.25) {};
\node at (0, .25) {\tiny $1$};
\node at (0.5, .25) {\tiny $1$};
\end{tikzpicture}
}
\endxy, \quad 
\xy
(0,0)*{
\begin{tikzpicture} [scale=1]
\draw[very thick,directed=1] (0,1) to [out=270,in=180] (.25,.5) to [out=0,in=270] (.5,1);
\node at (.25,.25) {};
\node at (0,1.25) {\tiny $1$};
\node at (.5,1.25) {\tiny $1$};
\end{tikzpicture}
}
\endxy := q^{\frac{1}{n}-n} 
\xy
(0,0)*{
\begin{tikzpicture} [scale=1]
\draw [very thick, ->] (0,1.5) to (.5,1);
\draw [very thick, ->] (.35,1.35) to (.5,1.5);
\draw [very thick] (.22,1.2) to (0,1);
\draw[very thick,rdirected=-.95] (0,1) to [out=270,in=180] (.25,.5) to [out=0,in=270] (.5,1);
\node at (.25,.25) {};
\node at (0,1.65) {\tiny $1$};
\node at (.6,1.65) {\tiny $1$};
\end{tikzpicture}
}
\endxy, \quad 
\]
and their thick versions are obtained by explosion.
\end{Definition}

\begin{Definition}\label{symm-webs} 
The \textit{thin symmetric (upward-downward pointing) category} of $\mathfrak{gl}_n$-webs, denote it 
by $\mathbf{SWeb}_{\uparrow,\downarrow}(\mathfrak{gl}_n)$, is the quotient category obtained from $\mathbf{fSWeb}_{\uparrow,\downarrow}(\mathfrak{gl}_n)$ by imposing relations \autoref{(co)assoc} and \autoref{thinsqure-sw} together with their downward pointing versions, and relations \autoref{inv-br}-\autoref{circ-ev} from below.\\
Invertibility of thin leftward crossings  
\begin{equation} \label{inv-br} 
\xy
(0,0)*{
\begin{tikzpicture}[scale=.3]
\draw [very thick, ->] (1,-1) to (-1,1);
\draw [very thick] (1,1) to (0.25,0.25);
\draw [very thick, ->] (-0.25,-0.25) to (-1,-1);
\node at (-1,-1.5) {\tiny $1$};
\node at (1,-1.5) {\tiny $1$};
\node at (-1,1.5) {\tiny $1$};
\node at (1,1.5) {\tiny $1$};
\end{tikzpicture}
};
\endxy  \quad 
\text{and} \quad
\xy
(0,0)*{
\begin{tikzpicture}[scale=.3]
\draw [very thick, ->] (1,1) to (-1,-1);
\draw [very thick, ->] (-0.25,0.25) to (-1,1);
\draw [very thick] (0.25,-0.25) to (1,-1);
\node at (-1,-1.5) {\tiny $1$};
\node at (1,-1.5) {\tiny $1$};
\node at (-1,1.5) {\tiny $1$};
\node at (1,1.5) {\tiny $1$};
\end{tikzpicture}
};
\endxy.
\end{equation}
The thin zigzag relation
\begin{equation}\label{eq-snake} 
\xy
(0,0)*{\reflectbox{
\begin{tikzpicture} [scale=1]
\draw[very thick] (0,0) to [out=90,in=180] (.25,.5) to [out=0,in=90] (.5,0);
\draw[very thick] (-.5,0) to [out=270,in=180] (-.25,-.5) to [out=0,in=270] (0,0);
\draw[very thick, directed=1] (0.5,0) to (0.5,-.75);
\draw[very thick] (-.5,0) to (-.5,.75);
\node at (0.5,-.9) {\reflectbox{\tiny $1$}};
\node at (-0.5,.9) {\reflectbox{\tiny $1$}};
\end{tikzpicture}
}};
\endxy=
\xy
(0,0)*{
\begin{tikzpicture} [scale=1]
\draw[very thick, rdirected=-.9] (0,0) to (0,1.5);
\node at (0,-.15) {\tiny $1$};
\node at (0,1.65) {\tiny $1$};
\end{tikzpicture}
}
\endxy
\quad , \quad
\xy
(0,0)*{
\begin{tikzpicture} [scale=1]
\draw[very thick] (0,0) to [out=90,in=180] (.25,.5) to [out=0,in=90] (.5,0);
\draw[very thick] (-.5,0) to [out=270,in=180] (-.25,-.5) to [out=0,in=270] (0,0);
\draw[very thick] (0.5,0) to (0.5,-.75);
\draw[very thick, directed=1] (-.5,0) to (-.5,.75);
\node at (0.5,-.9) {\tiny $1$};
\node at (-0.5,.9) {\tiny $1$};
\end{tikzpicture}
};
\endxy =
\xy
(0,0)*{
\begin{tikzpicture} [scale=1]
\draw[very thick, directed=1] (0,0) to (0,1.5);
\node at (0,-.15) {\tiny $1$};
\node at (0,1.65) {\tiny $1$};
\end{tikzpicture}
}
\endxy.
\end{equation}
The merge-split slides
\begin{equation}\label{eq-vertexslide}
\xy
(0,0)*{
\begin{tikzpicture}[scale=.3]
\draw [very thick, directed=.5] (0,-1) to (0,.75);
\draw [very thick] (0,.75) to [out=30,in=270] (1,2.5);
\draw [very thick] (0,.75) to [out=150,in=270] (-1,2.5); 
\draw [very thick] (1,2.5) to [out=90,in=180] (2,3.5) to [out=0,in=90] (3,2.5);
\draw [very thick] (-1,2.5) to [out=90,in=180] (2,4.5) to [out=0,in=90] (5,2.5);
\draw [very thick,rdirected=-.75] (3,-1) to (3,2.5);
\draw [very thick,rdirected=-.75] (5,-1) to (5,2.5);
\node at (0,-1.35) {\tiny $k{+}1$};
\node at (3,-1.35) {\tiny $1$};
\node at (5,-1.35) {\tiny $k$};
\end{tikzpicture}
};
\endxy=
\xy
(0,0)*{
\begin{tikzpicture}[scale=.3]
\draw [very thick] (4, .75) to (4,2.5);
\draw [very thick,rdirected=-.75] (5,-1) to [out=90,in=330] (4,.75);
\draw [very thick,rdirected=-.75] (3,-1) to [out=90,in=210] (4,.75);
\draw [very thick, directed=.5] (0, -1) to (0,2.5);
\draw [very thick] (0,2.5) to [out=90,in=180] (2,4.5) to [out=0,in=90] (4,2.5); 
\node at (0,-1.35) {\tiny $k{+}1$};
\node at (3,-1.35) {\tiny $1$};
\node at (5,-1.35) {\tiny $k$};
\end{tikzpicture}
};
\endxy\quad , \quad
\xy
(0,0)*{\reflectbox{
\begin{tikzpicture}[scale=.3]
\draw [very thick,directed=.5] (0,-1) to (0,.75);
\draw [very thick] (0,.75) to [out=30,in=270] (1,2.5);
\draw [very thick] (0,.75) to [out=150,in=270] (-1,2.5); 
\draw [very thick] (1,2.5) to [out=90,in=180] (2,3.5) to [out=0,in=90] (3,2.5);
\draw [very thick] (-1,2.5) to [out=90,in=180] (2,4.5) to [out=0,in=90] (5,2.5);
\draw [very thick, rdirected=-.75] (3,-1) to (3,2.5);
\draw [very thick,rdirected=-.75] (5,-1) to (5,2.5);
\node at (0,-1.35) {\reflectbox{\tiny $k{+}1$}};
\node at (3,-1.35) {\reflectbox{\tiny $1$}};
\node at (5,-1.35) {\reflectbox{\tiny $k$}};
\end{tikzpicture}}
};
\endxy=
\xy
(0,0)*{\reflectbox{
\begin{tikzpicture}[scale=.3]
\draw [very thick] (4, .75) to (4,2.5);
\draw [very thick,rdirected=-.75] (5,-1) to [out=90,in=330] (4,.75);
\draw [very thick,rdirected=-.75] (3,-1) to [out=90,in=210] (4,.75);
\draw [very thick, directed=.5] (0, -1) to (0,2.5);
\draw [very thick] (0,2.5) to [out=90,in=180] (2,4.5) to [out=0,in=90] (4,2.5); 
\node at (0,-1.35) {\reflectbox{\tiny $k{+}1$}};
\node at (3,-1.35) {\reflectbox{\tiny $1$}};
\node at (5,-1.35) {\reflectbox{\tiny $k$}};
\end{tikzpicture}
}};
\endxy, \ k\geq 1
\end{equation}
together with their cup and opposite orientation versions.\\
The 1-labelled circle removal
\begin{equation}\label{circ-ev}
\xy
(0,0)*{
\begin{tikzpicture} [scale=1]
\draw[very thick, directed=.53] (-1,0) circle (0.3cm);
\node at (-0.5,0) {\tiny $1$};
\end{tikzpicture}}
\endxy = [n].
\end{equation}
\end{Definition}

Note that $1_\uparrow$ and $1_\downarrow$ are duals to each other by the thin zigzag relation \autoref{eq-snake}. Later we will show that for every $k \geq 1$, $k_\uparrow$ and $k_\downarrow$ are duals to each other meaning that 
$\textbf{SWeb}_{\uparrow,\downarrow}(\mathfrak{gl}_n)$ is endowed with a rigid structure.

\begin{Remark}\label{down merges&splits} 
We could also define $\textbf{SWeb}_{\uparrow,\downarrow}(\mathfrak{gl}_n)$ using only upward pointing merges and splits but then we would have to include caps and cups of \textit{any} thickness as generators. In that case we would define the downward pointing merges and splits like below
\[
\xy
(0,0)*{

};
\endxy  
\]
(see also \autoref{downard thick merge and split}).
\end{Remark}

\subsection{Properties of webs} 

We now give the first consequences of the defining web relations.

\begin{Lemma}\label{thick-snake}  
\ochanged{The relations} \autoref{eq-snake} and \autoref{eq-vertexslide} imply the thick zigzag relations
\begin{equation*}
\xy
(0,0)*{\reflectbox{
%
}
\endxy.
\end{equation*}
\end{Lemma}

\begin{proof}
Using \autoref{thick (co)assoc} we see that we can adapt the proof of \cite[Lemma 2.21]{RoTu-symmetric-howe}. Namely, after explosion, using thin zigzag and merge-split slides we have
\[
\xy
(0,0)*{\reflectbox{
%
}
\endxy
\]
and similarly for the other thick zigzag relation.
\end{proof}

A direct consequence of this lemma is the following.

\begin{Corollary} \label{web cat is pivotal} 
The category $\mathbf{SWeb}_{\uparrow,\downarrow}(\mathfrak{gl}_n)$ can be endowed with a pivotal structure, where $k_\uparrow$ and $k_\downarrow$ are dual to each other.\qed
\end{Corollary}

We will use this pivotal structure throughout this paper.

\begin{Lemma}\label{thick merge and split slides} 
The thick versions of relation \autoref{eq-vertexslide} hold, i.e. 
\[\xy
(0,0)*{
%
}};
\endxy, \ k,l \geq 1
\]
as well as their cup and opposite orientation versions.
\end{Lemma}

\begin{proof} 
We prove the first relation via induction on $k+l$. The case $k+l=1$ follows trivially. Assume it holds for thickness smaller or equal to $k+l-1$. We have 
\[\xy
(0,0)*{
%
};
\endxy,
\]
where in the third equality we used the induction hypothesis to slide the $(k,l-1)$ split and in the last equality we used relation \autoref{eq-vertexslide} for the $(k+l-1,1)$ split, then associativity and the digon removal.
\end{proof}

\begin{Corollary}\label{downard thick merge and split} 
We have the following
\[
\xy
(0,0)*{
%
};
\endxy.  
\]
\end{Corollary}

\begin{proof}
This is a direct consequence of the previous lemma and thick zigzag relation. 
\end{proof}

The next lemma shows that in $\mathbf{SWeb}_{\uparrow,\downarrow}(\mathfrak{gl}_n)$ we have a relation akin to the usual planar isotopies, namely a relation of Reidemeister 1 type.

\begin{Lemma}\label{thin Reidemeister}
The following holds 
\begin{equation*} 
\xy
(0,0)*{
%
 };
\endxy
\end{equation*}
as well as its downward pointing version and 
\begin{equation*} 
\xy
(0,0)*{
%
 };
\endxy
\end{equation*}
as well as its downward pointing version.
\end{Lemma}

\begin{proof}
From the definition of the rightward cap
\begin{equation*}
q^{n- \frac{1}{n}}
\xy 
(0,0)*{
%
}
\endxy  
\end{equation*}
after adding a crossing and using the fact that the right over-crossing is the inverse of the left under-crossing we obtain 
\begin{equation*} q^{n- \frac{1}{n}}
\xy 
(0,0)*{
%
}
\endxy.
\end{equation*}
Further, by adding caps on both sides we have
\begin{equation*}q^{n- \frac{1}{n}}
\xy 
(0,0)*{
%
}
\endxy
\end{equation*}
and thus by straightening the left-hand side
\begin{equation*} q^{n- \frac{1}{n}}
\xy 
(0,0)*{
%
 };
\endxy
\end{equation*}
where in the first equality we used the definition of the left crossing and in the second we used the zigzag relation. This means that 
\begin{equation*}
q^{n- \frac{1}{n}}
\xy 
(0,0)*{
%
 };
\endxy.
\end{equation*} 
All other equalities can be proved verbatim.
\end{proof}

\begin{Corollary} 
The 1-circle evaluation with reversed orientation also holds, i.e.
\begin{gather}\label{reverse circ-ev}
\xy
(0,0)*{
%
}
\endxy = [n].
\end{gather}
\end{Corollary}
\begin{proof}
This follows from \autoref{right cap&cup}, \autoref{mates}, the previous lemma and \ochanged{the} relation \autoref{circ-ev}.
\end{proof}

\begin{Corollary}\label{sym dumb removal} 
The \textit{thin sideways dumbbell removal} holds, i.e.
\[
\xy
(-5,0)*{
%
};
\endxy,
\]
together with the downward pointing version.
\end{Corollary}

\begin{proof}
This follows from \autoref{thin-braid} and \autoref{thin Reidemeister} and \autoref{reverse circ-ev} as below
\[
\xy
(-5,0)*{
%
};
\endxy
\]
and similarly the downward pointing version.
\end{proof}

\begin{Lemma}\label{thick circle ev} 
We have the following thick circle evaluation formula 
\begin{equation*}
\xy
(0,0)*{
%
.
\end{equation*}
\end{Lemma}

\begin{proof}
We prove this by induction on $k$. The case $k=1$ is relation \autoref{circ-ev}. Using the induction hypothesis we have
\[
\xy
(0,0)*{
%
};
\endxy
\]
\[
\stackrel{\autoref{sym dumb removal}} {=}
\frac{1}{[k]} [n+1]
\xy
(0,0)*{
%
,
\]
where the last equality is a straightforward computation.
\end{proof}

\begin{Lemma}\label{the (k,l) dumbbell removal} 
The $(k,l)$-dumbbell removal, i.e.
\[
\xy
(-5,0)*{
%
};
\endxy,
\]
as well as its downward pointing version hold.
\end{Lemma}

\begin{proof}
We prove the second equality via induction on $l$. For $l=1$ we have
\[
\xy
(-5,0)*{
%
};
\endxy 
\]
\[
\stackrel{\autoref{sym dumb removal}} {=} [n+1] \xy
(0,0)*{
%
};
\endxy.
\]
Now, assume the relation holds for $(l_1,k)$-dumbbells, where $l_1 \leq l-1$. For the $(l,k)$-dumbbell we have 
\[
\xy
(-5,0)*{
%
};
\endxy.
\]
Finally, from the induction hypothesis and the case $l=1$ from above, the \ochanged{latter} is equal to
\[ 
\frac{1}{[l]} %
};
\endxy
\]
completing the proof.
\end{proof}

\begin{Lemma}\label{op-zigzag} 
The other versions of (thick) zigzag relations also hold:
\begin{equation*}
\xy
(0,0)*{\reflectbox{
%
}
\endxy.
\end{equation*}
\end{Lemma}

\begin{proof}
From the previous two lemmas we have
\begin{equation*}
\xy
(0,0)*{
%
};
\endxy.
\end{equation*}
Now we use this together with the digon removal and the relations \autoref{eq-snake} to obtain the thick version
\begin{equation*}
\xy
(0,0)*{
%
};
\endxy.
\end{equation*}
The other equation is proved similarly.
\end{proof}

\begin{Lemma}\label{pitchfork}
The so-called pitchfork relations hold
\[
\raisebox{-.01cm}{\xy
(0,0)*{
%
};
\endxy}}
\]
together with the merge and downwards pointing versions.
\end{Lemma}

\begin{proof} We prove the first relation via induction of $k+l$. For $k+l=1$, this follows trivially. Assume the relation holds for labels whose sum is smaller than $k+l$. We explode the $k$-edge as below
\[ \raisebox{-.01cm}{\xy
(0,0)*{
%
};
\endxy}.
\]
Now, using the induction hypothesis, we slide the edges labelled by 1 and $k-1$ through the $(l_1,l_2)$-split and after removing the digon we obtain 
\[
\raisebox{-.01cm}{\xy
(0,0)*{
%
};
\endxy}.
\]
The second relation is shown verbatim.
\end{proof}

\begin{Corollary}\label{dumbbell-cross slide} 
Combining the result above for merges and splits, we have that also the sliding of edges through dumbbells hold.\qed
\end{Corollary}

\begin{Corollary}\label{thickR2&R3} 
The thick versions of Reidemeister 2 and 3 from \autoref{thin R2&R3} hold. Namely,
\[
\xy
(0,0)*{
%
};
\endxy.
\]
\end{Corollary}

\begin{proof} 
We prove only the first equality in Reidemeister 2. From \autoref{thick-braid}, \autoref{dumbbell-cross slide}, \autoref{thin R2&R3} and the digon removal we have
\[
\xy
(0,0)*{
%
};
\endxy,
\]
where in the second equation we used \autoref{dumbbell-cross slide}  to move the central dumbbells up and then we used the fact that thin over and under-crossings are inverses to each other.
\end{proof}

Note that \ochanged{it follows from the previous lemma} that $\mathbf{SWeb}_{\uparrow,\downarrow}(\mathfrak{gl}_n)$ is a braided category \ochanged{(using the $(k,l)$-crossings, which we will always use for the braiding)}. In particular we have that the $(l,k)$ over-crossing and the $(k,l)$ under-crossing are inverses to each other.
Another direct consequence is the following.

\begin{Corollary}\label{thck left cross inv}
The leftward and rightward thick crossings are inverse to each other.\qed
\end{Corollary}

\begin{Proposition}\label{dualising webs} 
For any two objects $\vec{k}$ and $\vec{l}$ of $\mathbf{SWeb}_{\uparrow,\downarrow}(\mathfrak{gl}_n)$ 
there exist objects $\vec{k}^{\prime}$ and $\vec{l}^{\prime}$ (\changed{implicitly} defined in the proof) of $\mathbf{SWeb}_{\uparrow}(\mathfrak{gl}_n)$ such that 
there is an isomorphism of $\Bbbk(q^{\frac{1}{n}})$-vector spaces
\[
\mathrm{Hom}_{\mathbf{SWeb}_{\uparrow,\downarrow}(\mathfrak{gl}_n)}(\vec{k},\vec{l})
\cong
\mathrm{Hom}_{\mathbf{SWeb}_{\uparrow}(\mathfrak{gl}_n)}(\vec{k}^{\prime},\vec{l}^{\prime}).
\]
\end{Proposition}
\begin{proof} Let $\mathrm{w}$ be a web in $ \mathrm{Hom}_{\mathbf{SWeb}_{\uparrow,\downarrow}(\mathfrak{gl}_n)}(\vec{k},\vec{l})$.
Invertibility of crossings (\autoref{thickR2&R3} and \autoref{thck left cross inv}) and the invertibility of zigzag relations \autoref{eq-snake} imply that the usual isomorphism in braided pivotal categories holds. In pictures, this isomorphism looks e.g.
\[ \mathrm{w} \mapsto
\begin{tikzpicture}[anchorbase,scale=1.25]
\draw[very thick] (0.25,1.5) to (0.25,2);
\draw[very thick] (1.25,1.5) to[out=90,in=270] (0.75,2);
\draw[very thick] (2.25,1.5) to[out=90,in=270] (1.25,2);
\draw[very thick,crossline] (1.75,1.5) to[out=90,in=270] (2.25,2);
\draw[very thick,crossline] (0.75,1.5) to[out=90,in=270] (1.75,2);
\draw[very thick] (0.25,-1) to (0.25,-0.5);
\draw[very thick] (1.25,-1) to[out=90,in=270] (0.75,-0.5);
\draw[very thick,crossline] (0.75,-1) to[out=90,in=270] (1.25,-0.5);
\draw[very thick] (1.75,-1) to (1.75,-0.5);
\draw[very thick] (2.25,-1) to (2.25,-0.5);
\draw[very thick] (0,0) rectangle node[pos=0.5]{w} (2.5,1);
\draw[very thick,directed=0.55] (0.25,1.5) to (0.25,1);
\draw[very thick,directed=0.55] (0.75,1) to (0.75,1.5);
\draw[very thick,directed=0.55] (1.25,1.5) to (1.25,1);
\draw[very thick,directed=0.55] (1.75,1) to (1.75,1.5);
\draw[very thick,directed=0.55] (2.25,1.5) to (2.25,1);
\draw[very thick,directed=0.55] (0.25,-0.5) to (0.25,0);
\draw[very thick,directed=0.55] (0.75,0) to (0.75,-0.5);
\draw[very thick,directed=0.55] (1.25,-0.5) to (1.25,0);
\draw[very thick,directed=0.55] (1.75,0) to (1.75,-0.5);
\draw[very thick,directed=0.55] (2.25,0) to (2.25,-0.5);
\draw[very thick] (0.25,2) to[out=90,in=0]
(0,2.25) to[out=180,in=90] (-0.25,2);
\draw[very thick] (0.75,2) to[out=90,in=0]
(0,2.5) to[out=180,in=90] (-0.75,2);
\draw[very thick] (1.25,2) to[out=90,in=0]
(0,2.75) to[out=180,in=90] (-1.25,2);
\draw[very thick] (2.25,-1) to[out=270,in=180]
(2.5,-1.25) to[out=0,in=270] (2.75,-1);
\draw[very thick] (1.75,-1) to[out=270,in=180]
(2.5,-1.5) to[out=0,in=270] (3.25,-1);
\draw[very thick] (1.25,-1) to[out=270,in=180]
(2.5,-1.75) to[out=0,in=270] (3.75,-1);
\draw[very thick] (1.75,2) to (1.75,3)node[above]{$l_1^{\prime}$};
\draw[very thick] (2.25,2) to (2.25,3)node[above]{$l_2^{\prime}$};
\draw[very thick] (2.75,-1) to (2.75,3)node[above]{$l_3^{\prime}$};
\draw[very thick] (3.25,-1) to (3.25,3)node[above]{$l_4^{\prime}$};
\draw[very thick] (3.75,-1) to (3.75,3)node[above]{$l_5^{\prime}$};
\draw[very thick] (-1.25,2) to (-1.25,-2)node[below]{$k_1^{\prime}$};
\draw[very thick] (-0.75,2) to (-0.75,-2)node[below]{$k_2^{\prime}$};
\draw[very thick] (-0.25,2) to (-0.25,-2)node[below]{$k_3^{\prime}$};
\draw[very thick] (0.25,-1) to (0.25,-2)node[below]{$k_4^{\prime}$};
\draw[very thick] (0.75,-1) to (0.75,-2)node[below]{$k_5^{\prime}$};
\end{tikzpicture},
\]
which is evidently an invertible operation. \changed{From this picture it is then easy to write down formulas, depending on $\vec{k}$ and $\vec{l}$, for $\vec{k}^{\prime}$ and $\vec{l}^{\prime}$,}
\end{proof}

\begin{Remark}\label{dualising is possible in all web cats with duals}
Note that \autoref{dualising webs} only uses invertibility of crossings and the zigzag relations. In particular, it also holds for other web categories with duals and invertible crossings which are going to be discussed in this paper, not just 
$\mathbf{SWeb}_{\uparrow,\downarrow}(\mathfrak{gl}_n)$.
\end{Remark}

The splits and merges, when composed with the corresponding crossings, are ``mirrored'' as we will explain in the following lemma.  

\begin{Lemma}\label{splits&merges comp crossings} 
We have
\[
\xy
(0,0)*{
\begin{tikzpicture}[scale=.3]
\draw [very thick,] (-1,4.5) to (1,2.5);
\draw [very thick,crossline] (-1,2.5) to (1,4.5);
\draw [very thick, directed=.55] (0,-1) to (0,.75);
\draw [very thick, directed=.65] (0,.75) to [out=30,in=270] (1,2.5);
\draw [very thick, directed=.65] (0,.75) to [out=150,in=270] (-1,2.5); 
\node at (0, -1.5) {\tiny $k{+}l$};
\node at (-1,3.25) {\tiny $l$};
\node at (1,3.25) {\tiny $k$};
\end{tikzpicture}
};
\endxy= q^{-\frac{kl}{n}+kl}
\xy
(0,0)*{
\begin{tikzpicture}[scale=.3]
\draw [very thick, directed=.55] (0,-1) to (0,.75);
\draw [very thick, directed=.65] (0,.75) to [out=30,in=270] (1,2.5);
\draw [very thick, directed=.65] (0,.75) to [out=150,in=270] (-1,2.5); 
\node at (0, -1.5) {\tiny $k{+}l$};
\node at (-1,3) {\tiny $k$};
\node at (1,3) {\tiny $l$};
\end{tikzpicture}
};
\endxy
\]
\[
\xy
(0,0)*{
\begin{tikzpicture}[scale=.3]
\draw [very thick] (-1,2.5) to (1,4.5);
\draw [very thick,crossline] (-1,4.5) to (1,2.5);
\draw [very thick, directed=.55] (0,-1) to (0,.75);
\draw [very thick, directed=.65] (0,.75) to [out=30,in=270] (1,2.5);
\draw [very thick, directed=.65] (0,.75) to [out=150,in=270] (-1,2.5); 
\node at (0, -1.5) {\tiny $k{+}l$};
\node at (-1,3.25) {\tiny $l$};
\node at (1,3.25) {\tiny $k$};
\end{tikzpicture}
};
\endxy= q^{\frac{kl}{n}+kl}
\xy
(0,0)*{
\begin{tikzpicture}[scale=.3]
\draw [very thick, directed=.55] (0,-1) to (0,.75);
\draw [very thick, directed=.65] (0,.75) to [out=30,in=270] (1,2.5);
\draw [very thick, directed=.65] (0,.75) to [out=150,in=270] (-1,2.5); 
\node at (0, -1.5) {\tiny $k{+}l$};
\node at (-1,3) {\tiny $k$};
\node at (1,3) {\tiny $l$};
\end{tikzpicture}
};
\endxy
\]
and merges have the same behaviour.
\end{Lemma}

\begin{proof}
We prove this via the exploded crossing formula in \autoref{thick-braid}. First, the case $ k = l = 1$ can be verified via a direct computation using \autoref{thin-braid}. The relevant diagrammatic calculation is now
\[
\begin{tikzpicture}[anchorbase,scale=0.75]
\draw [very thick] (1,2) to (0,3)node[above]{$k$};
\draw [very thick,crossline] (0,2) to (1,3)node[above]{$l$};
\draw [very thick,directed=.5] (0.5,1) to[out=150,in=270] (0,2);
\draw [very thick,directed=.5] (0.5,1) to[out=30,in=270] (1,2);
\draw [very thick,directed=.5] (0.5,0)node[below]{$k+l$} to (0.5,1);
\end{tikzpicture}
=
\frac{1}{[k]![l]!}
\begin{tikzpicture}[anchorbase,scale=0.75]
\draw [very thick] (2,0) to (0,2)node[left]{$1$};
\draw [very thick] (3,0) to (1,2)node[right]{$1$};
\draw [very thick,crossline] (0,0) to (2,2)node[left]{$1$};
\draw [very thick,crossline] (1,0) to (3,2)node[right]{$1$};
\node at (0.5,0){$\cdots$};
\node at (2.5,0){$\cdots$};
\node at (0.5,2){$\cdots$};
\node at (2.5,2){$\cdots$};
\draw [very thick,rdirected=.5] (0.5,2.5) to[out=180,in=90] (0,2);
\draw [very thick,rdirected=.5] (0.5,2.5) to[out=0,in=90] (1,2);
\draw [very thick,rdirected=.5] (0.5,3.5)node[above]{$k$} to(0.5,2.5);
\draw [very thick,rdirected=.5] (2.5,2.5) to[out=180,in=90] (2,2);
\draw [very thick,rdirected=.5] (2.5,2.5) to[out=0,in=90] (3,2);
\draw [very thick,rdirected=.5] (2.5,3.5)node[above]{$l$} to(2.5,2.5);
\draw [very thick,directed=.5] (1.5,-1) to[out=180,in=270] (0.5,-.5);
\draw [very thick,directed=.5] (0.5,-.5) to[out=180,in=270] (0,0);
\draw [very thick,directed=.5] (0.5,-.5) to[out=0,in=270] (1,0);
\draw [very thick,directed=.5] (1.5,-1) to[out=0,in=270] (2.5,-.5);
\draw [very thick,directed=.5] (2.5,-.5) to[out=180,in=270] (2,0);
\draw [very thick,directed=.5] (2.5,-.5) to[out=0,in=270] (3,0);
\draw [very thick,directed=.5] (1.5,-2)node[below]{$k+l$} to (1.5,-1);
\end{tikzpicture}
=
\frac{q^{-\frac{kl}{n}+kl}}{[k]![l]!}
\begin{tikzpicture}[anchorbase,scale=0.75]
\node at (0.5,0){$\cdots$};
\node at (2.5,0){$\cdots$};
\draw [very thick,rdirected=.5] (0.5,0.5) to[out=180,in=90] (0,0);
\draw [very thick,rdirected=.5] (0.5,0.5) to[out=0,in=90] (1,0);
\draw [very thick,rdirected=.5] (0.5,1.5)node[above]{$k$} to(0.5,0.5);
\draw [very thick,rdirected=.5] (2.5,0.5) to[out=180,in=90] (2,0);
\draw [very thick,rdirected=.5] (2.5,0.5) to[out=0,in=90] (3,0);
\draw [very thick,rdirected=.5] (2.5,1.5)node[above]{$l$} to(2.5,0.5);
\draw [very thick,directed=.5] (1.5,-1) to[out=180,in=270] (0.5,-.5);
\draw [very thick,directed=.5] (0.5,-.5) to[out=180,in=270] (0,0);
\draw [very thick,directed=.5] (0.5,-.5) to[out=0,in=270] (1,0);
\draw [very thick,directed=.5] (1.5,-1) to[out=0,in=270] (2.5,-.5);
\draw [very thick,directed=.5] (2.5,-.5) to[out=180,in=270] (2,0);
\draw [very thick,directed=.5] (2.5,-.5) to[out=0,in=270] (3,0);
\draw [very thick,directed=.5] (1.5,-2)node[below]{$k+l$} to (1.5,-1);
\end{tikzpicture}
\]
\[
=
q^{-\frac{kl}{n}+kl}
\begin{tikzpicture}[anchorbase,scale=0.75]
\draw [very thick,directed=.5] (0.5,1)
to[out=150,in=270] (0,2)node[above]{$k$};
\draw [very thick,directed=.5] (0.5,1)
to[out=30,in=270] (1,2)node[above]{$l$};
\draw [very thick,directed=.5] (0.5,0)node[below]{$k+l$} to (0.5,1);
\end{tikzpicture}
.\]
Let us explain the second equality. Successively, using associativity in the bottom bouquet we can always match a crossing with a corresponding split, giving a factor $q^{-\frac{1}{n}+1}$. As there are $kl$ crossings to be absorbed, then the overall factor is $q^{-\frac{kl}{n}+kl}$. All other cases follow mutatis mutandis.
\end{proof}

Further, we have the next two lemmas.

\begin{Lemma} \label{thick Reidemeister} 
The following relation of Reidemeister 1 type 
\begin{equation*} 
\xy
(0,0)*{
\begin{tikzpicture}[scale=.3]
\draw [very thick] (-1,-1) to (1,1);
\draw [very thick] (1,-1) to (.25,-.25);
\draw [very thick] (-1,1) to (-.25,.25);
\draw[very thick] (-1,1) to [in=0, out=90] (-1.5,1.5) to [in=90,out=180] (-2,1);
\draw [very thick] (-2,1) to (-2,-1);
\draw [very thick] (-1,-1) to [in=0, out=270] (-1.5,-1.5) to [in=270, out=180] (-2,-1);
\draw [very thick,->] (1,1) to (1,3);
\draw [very thick] (1,-1) to (1,-2.5);
\node at (1,-3) {\tiny $k$};
\node at (1,3.5) {\tiny $k$};
\end{tikzpicture}
};
\endxy=q^{k(k-\frac{k}{n}+n-1)} 
\xy 
(0,0)*{
\begin{tikzpicture}[scale=.3]
\draw [very thick, ->] (0,-2.5) to (0,3);
\node at (0,-3) {\tiny $k$};
\node at (0,3.5) {\tiny $k$};
\end{tikzpicture}
};
\endxy=
\xy
(0,0)*{
\begin{tikzpicture}[scale=.3]
\draw [very thick] (-1,-1) to (1,1);
\draw [very thick] (1,-1) to (.25,-.25);
\draw [very thick] (-.25,.25) to (-1,1);
\draw [very thick] (1,1) to [in=180, out=90] (1.5,1.5) to [in=90, out=0](2,1);
\draw [very thick] (2,-1) to [in=0, out=270] (1.5,-1.5) to [in=270, out=180](1,-1);
\draw [very thick] (2,-1) to (2,1);
\draw [very thick,->] (-1,1) to (-1,3);
\draw [very thick] (-1,-2.5) to (-1,-1);
\node at (-1,-3) {\tiny $k$};
\node at (-1,3.5) {\tiny $k$};
\end{tikzpicture} };
\endxy
\end{equation*}
as well as its downward pointing version and 
\begin{equation*} 
\xy
(0,0)*{
\begin{tikzpicture}[scale=.3]
\draw [very thick] (1,-1) to (-1,1);
\draw [very thick] (-1,-1) to (-.25,-.25);
\draw [very thick] (1,1) to (.25,.25);
\draw[very thick] (-1,1) to [in=0, out=90] (-1.5,1.5) to [in=90,out=180] (-2,1);
\draw [very thick] (-2,1) to (-2,-1);
\draw [very thick] (-1,-1) to [in=0, out=270] (-1.5,-1.5) to [in=270, out=180] (-2,-1);
\draw [very thick,->] (1,1) to (1,3);
\draw [very thick] (1,-1) to (1,-2.5);
\node at (1,-3) {\tiny $k$};
\node at (1,3.5) {\tiny $k$};
\end{tikzpicture}
};
\endxy=q^{-k(k-\frac{k}{n}+n-1)} 
\xy 
(0,0)*{
\begin{tikzpicture}[scale=.3]
\draw [very thick, ->] (0,-2.5) to (0,3);
\node at (0,-3) {\tiny $k$};
\node at (0,3.5) {\tiny $k$};
\end{tikzpicture}
};
\endxy=
\xy
(0,0)*{
\begin{tikzpicture}[scale=.3]
\draw [very thick] (1,-1) to (-1,1);
\draw [very thick] (-1,-1) to (-.25,-.25);
\draw [very thick] (.25,.25) to (1,1);
\draw [very thick] (1,1) to [in=180, out=90] (1.5,1.5) to [in=90, out=0](2,1);
\draw [very thick] (2,-1) to [in=0, out=270] (1.5,-1.5) to [in=270, out=180](1,-1);
\draw [very thick] (2,-1) to (2,1);
\draw [very thick,->] (-1,1) to (-1,3);
\draw [very thick] (-1,-2.5) to (-1,-1);
\node at (-1,-3) {\tiny $k$};
\node at (-1,3.5) {\tiny $k$};
\end{tikzpicture} };
\endxy
\end{equation*}
as well as its downward pointing version hold.
\end{Lemma}
\begin{proof}
Using explosion, the previous lemma and the various established topological properties, the main thing to observe is that a Reidemeister 1 picture with $k$ parallel thin strands can be topologically rewritten as follows
\[
\begin{tikzpicture}[anchorbase,scale=0.75]
\draw [very thick,directed=0.2] (2,0) to (0,2)node[left]{$1$};
\draw [very thick,directed=0.2] (3,0) to (1,2)node[right]{$1$};
\draw [very thick,crossline,directed=0.2] (0,0) to (2,2)node[left]{$1$};
\draw [very thick,crossline,directed=0.2] (1,0) to (3,2)node[right]{$1$};
\draw [very thick] (3,0) to[out=270,in=180]
(3.5,-0.5) to[out=0,in=270] (4,0) to (4,2)
to[out=90,in=0] (3.5,2.5) to[out=180,in=90] (3,2);
\draw [very thick] (2,0) to[out=270,in=180]
(3.5,-1) to[out=0,in=270] (5,0) to (5,2)
to[out=90,in=0] (3.5,3) to[out=180,in=90] (2,2);
\draw [very thick] (0,2) to (0,3);
\draw [very thick] (1,2) to (1,3);
\draw [very thick] (0,0) to (0,-1);
\draw [very thick] (1,0) to (1,-1);
\node at (0.5,0){$\cdots$};
\node at (2.5,0){$\cdots$};
\node at (0.5,2){$\cdots$};
\node at (2.5,2){$\cdots$};
\end{tikzpicture}
=
\begin{tikzpicture}[anchorbase,scale=0.75]
\draw [very thick] (0.25,1.5) to[out=180,in=270] (0,2.25) to (0,3);
\draw [very thick,crossline] (0,1) to (0,2.25)
to[out=90,in=180] (0.25,2.5) to[out=0,in=90] (0.5,2.25)
to (0.5,1.75) to[out=270,in=0] (0.25,1.5);
\draw [very thick] (1.25,1.5)
to[out=180,in=270] (1,2.25) to (1,3);
\draw [very thick,crossline] (1,1) to (1,2.25)
to[out=90,in=180] (1.25,2.5) to[out=0,in=90] (1.5,2.25)
to (1.5,1.75) to[out=270,in=0] (1.25,1.5);
\draw [very thick] (1,1) to (1,1.5);
\draw [very thick] (0,0) to (0,-1);
\draw [very thick] (1,0) to (1,-1);
\draw [very thick] (-0.5,0) rectangle (1.5,1)node[pos=0.5]{$ft$};
\node at (0.5,-0.25){$\cdots$};
\node at (0.5,1.25){$\cdots$};
\end{tikzpicture}
.\]
Here $ft$ is the full twist on $k$ strands which has $k(k-1)$ crossings. Now we first need to absorb the full twist using the previous lemma giving a factor $q^{-\frac{k(k-1)}{n}+k(k-1)}$. Then we use the thin Reidemeister 1 move from \autoref{thin Reidemeister} to get the factor $q^{k(n- \frac{1}{n})}$. Thus, the overall factor is $q^{k(k-\frac{k}{n}+n-1)} $. All other cases follow mutatis mutandis.
\end{proof}

\begin{Lemma}\label{thick cu(a)p sliding} 
The sliding of thick edges through thick caps holds, together with its cup, reversed orientation and over-crossing versions 
\begin{equation*}
\xy
(0,1)*{
\begin{tikzpicture}[scale=.3]
\draw [very thick, ->] (1.5,0) to [in=0,out=90] (0,1.5) to[in=90,out=180](-1.5,0);
\draw [very thick, ->](-1.3,1.3) to (-1.85,1.85);
\draw [very thick] (0,0) to (-.85,.9);
\node at (-1.5,-.5) {\tiny $k$};
\node at (1.5,-.5) {\tiny $k$};
\node at (0,-.5) {\tiny $l$};
\end{tikzpicture}
}
\endxy=
\xy
(0,1)*{
\begin{tikzpicture}[scale=.3]
\draw [very thick, ->] (1.5,0) to [in=0,out=90] (0,1.5) to[in=90,out=180](-1.5,0);
\draw [very thick] (0.3,0) to (0.75,1);
\draw [very thick, ->] (1,1.5) to (1.25,2);
\node at (-1.5,-.5) {\tiny $k$};
\node at (1.5,-.5) {\tiny $k$};
\node at (0,-.5) {\tiny $l$};
\node at (0,-.5) {\tiny $l$};
\end{tikzpicture}
}
\endxy.
\end{equation*} 
\end{Lemma}
\begin{proof} The case $k=l=1$ follows from a direct check using the definition of the leftward crossings. The general case is: 
\[
\xy
(0,1)*{
\begin{tikzpicture}[scale=.3]
\draw [very thick, ->] (1.5,0) to [in=0,out=90] (0,1.5) to[in=90,out=180](-1.5,0);
\draw [very thick, ->](-1.3,1.3) to (-1.85,1.85);
\draw [very thick] (0,0) to (-.85,.9);
\node at (-1.5,-.5) {\tiny $k$};
\node at (1.5,-.5) {\tiny $k$};
\node at (0,-.5) {\tiny $l$};
\end{tikzpicture}
}
\endxy=
\frac{1}{[k]![l]!}
\xy
(0,0)*{
\begin{tikzpicture}[scale=.2]
\draw [very thick, rdirected=-.5] (0,-1) to (0,.75);
\draw [very thick] (9,-1) to (9,.75);
\draw [very thick] (9,.75) to [out=30,in=270] (10,2.5);
\draw [very thick, ] (9,.75) to [out=150,in=270] (8,2.5);
\draw [very thick] (7,-1) to (6,.5);
\draw [very thick, ->] (1.5,8) to (.5,9.5);
\draw [very thick, ] (6,.5) to [out=-205,in=-95] (1.5,8);
\draw [very thick] (6,.5) to [out=-275,in=05] (1.5,8);
\draw [very thick] (0,.75) to [out=30,in=270] (1,2.5);
\draw [very thick] (0,.75) to [out=150,in=270] (-1,2.5); 
\draw [very thick,rdirected=-1, crossline] (1,2.5) to [out=90,in=180] (4.5,4.5) to [out=0,in=90] (8,2.5);
\draw [very thick,rdirected=-1, crossline] (-1,2.5) to [out=90,in=180] (4.5,6.5) to [out=0,in=90] (10,2.5);
\node at (0,-1.75) {\tiny $k$};
\node at (9,-1.75) {\tiny $k$};
\node at (2,3) {\tiny $1$};
\node at (1.25,5) {\tiny $1$};
\node at (4,5.25) {\tiny $\cdots$};
\node at (7,-1.75) {\tiny $l$};
\node at (.5,10.5) {\tiny $l$};
\node at (6.75,2) {\tiny $1$};
\node at (3,1.5) {\tiny $1$};
\end{tikzpicture}
};
\endxy=
\frac{1}{[k]![l]!}
\xy
(0,0)*{
\begin{tikzpicture}[scale=.2]
\draw [very thick, rdirected=-.5] (0,-1) to (0,.75);
\draw [very thick] (9,-1) to (9,.75);
\draw [very thick] (9,.75) to [out=30,in=270] (10,2.5);
\draw [very thick, ] (9,.75) to [out=150,in=270] (8,2.5);
\draw [very thick] (2,-1) to (3,.5);
\draw [very thick, ->] (8.2,8) to (9.6,9.5);
\draw [very thick, ] (3,.5) to [out=-265,in=-175] (8.2,8);
\draw [very thick] (3,.5) to [out=-335,in=-75] (8.2,8);
\draw [very thick] (0,.75) to [out=30,in=270] (1,2.5);
\draw [very thick] (0,.75) to [out=150,in=270] (-1,2.5); 
\draw [very thick,rdirected=-1, crossline] (1,2.5) to [out=90,in=180] (4.5,4.5) to [out=0,in=90] (8,2.5);
\draw [very thick,rdirected=-1, crossline] (-1,2.5) to [out=90,in=180] (4.5,6.5) to [out=0,in=90] (10,2.5);
\node at (0,-1.75) {\tiny $k$};
\node at (9,-1.75) {\tiny $k$};
\node at (2,3) {\tiny $1$};
\node at (1.25,5) {\tiny $1$};
\node at (2,-1.75) {\tiny $l$};
\node at (9.6,10.5) {\tiny $l$};
\node at (4,.5) {\tiny $1$};
\node at (2.5,1.5) {\tiny $1$};
\node at (5.5,5.25) {\tiny $\cdots$};
\end{tikzpicture}
};
\endxy=
\xy
(0,1)*{
\begin{tikzpicture}[scale=.3]
\draw [very thick, ->] (1.5,0) to [in=0,out=90] (0,1.5) to[in=90,out=180](-1.5,0);
\draw [very thick] (0.3,0) to (0.75,1);
\draw [very thick, ->] (1,1.5) to (1.25,2);
\node at (-1.5,-.5) {\tiny $k$};
\node at (1.5,-.5) {\tiny $k$};
\node at (0,-.5) {\tiny $l$};
\end{tikzpicture}
}
\endxy.
\]
Here we exploded the strings.
\end{proof}

\begin{Definition} \label{defn-symJW}
Let $k\in \mathbb{Z}_{> 0}$. The $k$-th \textit{symmetric $\mathfrak{gl_n}$-projector} $e_k$ is defined via
\[
e_k
=
\frac{1}{[k]!}
\;\;
\xy
(0,0)*{\reflectbox{
\begin{tikzpicture}[scale=.3]
\draw [very thick, directed=.7] (0,-1) to (0,.75);
\draw [very thick] (0,.75) to [out=30,in=270] (1,2.5);
\draw [very thick] (0,.75) to [out=150,in=270] (-1,2.5); 
\draw [very thick] (1,-2.75) to [out=90,in=330] (0,-1);
\draw [very thick] (-1,-2.75) to [out=90,in=210] (0,-1);
\draw [very thick] (1,2.5) to [out=30,in=270] (2,4.25);
\draw [very thick] (1,2.5) to [out=150,in=270] (0,4.25);
\draw [very thick] (2,-4.5) to [out=90,in=330] (1,-2.75);
\draw [very thick] (0,-4.5) to [out=90,in=210] (1,-2.75);
\draw [very thick] (2,4.25) to [out=30,in=270] (3,6.0);
\draw [very thick] (2,4.25) to [out=150,in=270] (1,6.0);
\draw [very thick] (3,-6.25) to [out=90,in=330] (2,-4.5);
\draw [very thick] (1,-6.25) to [out=90,in=210] (2,-4.5);
\draw [very thick] (-1,2.5) to (-1,6.0);
\draw [very thick] (0,4.25) to (0,6.0);
\draw [very thick] (-1,-2.5) to (-1,-6.25);
\draw [very thick] (0,-4.25) to (0,-6.25);
\node at (2,-6.25) {\tiny $\cdots$};
\node at (2,6) {\tiny $\cdots$};
\node at (-1.25,1.75) {\reflectbox{\tiny $1$}};
\node at (2,1.75) {\reflectbox{\tiny $k{-}1$}};
\node at (-1.25,-1.9) {\reflectbox{\tiny $1$}};
\node at (2,-1.9) {\reflectbox{\tiny $k{-}1$}};
\node at (0.8,0) {\reflectbox{\tiny $k$}};
\node at (-0.25,3.5) {\reflectbox{\tiny $1$}};
\node at (3,3.5) {\reflectbox{\tiny $k{-}2$}};
\node at (-0.25,-3.65) {\reflectbox{\tiny $1$}};
\node at (3,-3.65) {\reflectbox{\tiny $k{-}2$}};
\node at (0.75,5.25) {\reflectbox{\tiny $1$}};
\node at (4,5.25) {\reflectbox{\tiny $k{-}3$}};
\node at (0.75,-5.4) {\reflectbox{\tiny $1$}};
\node at (4,-5.4) {\reflectbox{\tiny $k{-}3$}};
\end{tikzpicture}
}};
\endxy
=
\xy
(0,0)*{
\begin{tikzpicture}[scale=.3]
\draw [double, directed=.7] (0,-1) to (0,.75);
\draw [very thick] (0,.75) to [out=30,in=270] (1,2.5);
\draw [very thick] (0,.75) to [out=150,in=270] (-1,2.5); 
\draw [very thick] (1,-2.75) to [out=90,in=330] (0,-1);
\draw [very thick] (-1,-2.75) to [out=90,in=210] (0,-1);
\node at (-1,3) {\tiny $1$};
\node at (0.1,3) {$\cdots$};
\node at (1,3) {\tiny $1$};
\node at (-1,-3.15) {\tiny $1$};
\node at (0.1,-3.15) {$\cdots$};
\node at (1,-3.15) {\tiny $1$};
\node at (-0.7,0) {\tiny $k$};
\end{tikzpicture}
};
\endxy.
\]
\end{Definition}

The reader familiar with the Jones--Wenzl recursion will recognize the following lemma. However, as a word of warning, the scalars 
are not the same as in the classical Jones--Wenzl recursion.

\begin{Lemma}\label{JW-recursion}
We have:
\[
e_k=
\begin{tikzpicture}[anchorbase,scale=.3]
\draw [double, directed=.7] (0,-1) to (0,.75);
\draw [very thick] (0,.75) to [out=30,in=270] (1,2.5);
\draw [very thick] (0,.75) to [out=150,in=270] (-1,2.5); 
\draw [very thick] (1,-2.75) to [out=90,in=330] (0,-1);
\draw [very thick] (-1,-2.75) to [out=90,in=210] (0,-1);
\node at (-1,3) {\tiny $1$};
\node at (0.1,3) {\tiny $\cdots$};
\node at (1,3) {\tiny $1$};
\node at (-1,-3.15) {\tiny $1$};
\node at (0.1,-3.15) {$\cdots$};
\node at (1,-3.15) {\tiny $1$};
\node at (-0.7,0) {\tiny $k$};
\end{tikzpicture}
-\frac{[k-2]}{[k]}
\begin{tikzpicture}[anchorbase,scale=.3]
\draw [double, directed=.7] (0,-1) to (0,.75);
\draw [very thick] (0,.75) to [out=30,in=270] (1,2.5);
\draw [very thick] (0,.75) to [out=150,in=270] (-1,2.5); 
\draw [very thick] (1,-2.75) to [out=90,in=330] (0,-1);
\draw [very thick] (-1,-2.75) to [out=90,in=210] (0,-1);
\draw [very thick] (3,-2.75) to (3,2.5);
\node at (-1,3) {\tiny $1$};
\node at (0.1,3) {$\cdots$};
\node at (1,3) {\tiny $1$};
\node at (-1,-3.15) {\tiny $1$};
\node at (0.1,-3.15) {$\cdots$};
\node at (1,-3.15) {\tiny $1$};
\node at (-1.25,0) {\tiny $k{-}1$};
\node at (3,-3.15) {\tiny $1$};
\node at (3,3) {\tiny $1$};
\end{tikzpicture}
+\frac{[k-1]}{[k]}
\begin{tikzpicture}[anchorbase,scale=.3]
\draw [double, directed=.7] (0,-1) to (0,.75);
\draw [very thick] (0,.75) to [out=30,in=270] (1,2.5);
\draw [very thick] (0,.75) to [out=150,in=270] (-1,2.5); 
\draw [very thick] (1,-2.75) to [out=90,in=330] (0,-1);
\draw [very thick] (3,-2.75) to (3,2.5);
\draw [very thick] (1,-2.75) to [out=270,in=180] (2,-4.5) to [out=0,in=270] (3,-2.75);
\draw [very thick] (1,-6.75) to [out=90,in=180] (2,-5) to [out=0,in=90] (3,-6.75);
\draw [double, directed=.7] (0,-10.25) to (0,-8.5);
\draw [very thick] (0,-8.5) to [out=30,in=270] (1,-6.75);
\draw [very thick] (0,-1) to [out=210,in=150] (0,-8.5); 
\draw [very thick] (-0.75,-1.5) to [out=300,in=60] (-0.75,-8);	
\draw [very thick] (1,-12) to [out=90,in=330] (0,-10.25);
\draw [very thick] (-1,-12) to [out=90,in=210] (0,-10.25); 
\draw [very thick] (3,-12) to (3,-6.75);
\node at (-1,3) {\tiny $1$};
\draw [very thick] (2,-5) to (2,-4.5);
\node at (3,-4.75) {\tiny $2$};
\node at (1,3) {\tiny $1$};
\node at (-1,-12.5) {\tiny $1$};
\node at (1,-12.5) {\tiny $1$};
\node at (-1.3,0) {\tiny $k{-}1$};
\node at (-1.3,-9.5) {\tiny $k{-}1$};
\node at (-1.3,-0.95) {\tiny $k{-}2$};
\node at (-1.3,-8.55) {\tiny $k{-}2$};
\node at (-2.25,-4.75) {\tiny $1$};
\node at (-0.875,-4.75) {$\cdots$};
\node at (0.5,-4.75) {\tiny $1$};
\node at (3,-12.5) {\tiny $1$};
\node at (3,3) {\tiny $1$};
\end{tikzpicture}
\]
for all $k\in\mathbb{Z}_{\geq 3}$.
\end{Lemma}

\begin{proof}
The proof is similar to that of \cite[Lemma 2.13]{RoTu-symmetric-howe}. (The proof of this uses the dumbbell relation, which will state in \autoref{S:Sym} below.)
\end{proof}

\subsection{Exterior and red-green webs}\label{S:Sym}

Recall the exterior $\mathfrak{sl}_{n}$-webs from \cite{CaKaMo-webs-skew-howe}. We point out here that for $\text{U}_{q}(\mathfrak{gl}_n)$ we have $\bigwedge_{q}^n (\Bbbk_{q}^n)\not \cong \Bbbk_q $, whereas for $\text{U}_{q}(\mathfrak{sl}_n)$ we have $\bigwedge_{q}^n (\Bbbk_{q}^n) \cong \Bbbk_q$. ($\bigwedge_{q}^n (\Bbbk_{q}^n)$ is the quantum analog of the determinant representation of the general linear group.) So, for exterior $\mathfrak{sl}_n$-webs the $n$-labelled edges are 0-labelled and thus are removed from the illustrations. \changed{There is also a $\mathfrak{gl}_n$-version of their calculus (with only minor changes), which we assume the reader has seen. Recall furthermore the calculus from \cite{TuVaWe-super-howe}.}

\begin{Remark}
All of the above and the below is formulated for symmetric 
webs, but exterior webs as in \cite{CaKaMo-webs-skew-howe} or red-green webs 
as in \cite{TuVaWe-super-howe} \changed{work} perfectly well using the same strategy.
There is one minor difference in the relations: 
the crossings are defined differently when braiding exterior powers:
\begin{equation*}
\xy
(0,0)*{
\begin{tikzpicture}[scale=.3]
\draw [very thick, ->] (-1,-1) to (1,1);
\draw [very thick] (1,-1) to (0.25,-0.25);
\draw [very thick, ->] (-0.25,0.25) to (-1,1);
\node at (-1,-1.5) {\tiny $1$};
\node at (1,-1.5) {\tiny $1$};
\node at (1,1.45) {\tiny $1$};
\node at (-1,1.45) {\tiny $1$};
\end{tikzpicture}
};
\endxy  = q^{1-\frac{1}{n}} \Bigg(  
\xy
(0,0)*{
\begin{tikzpicture}[scale=.3]
\draw [very thick, directed=.55] (-1, -2.5) to (-1,2.5);
\draw [very thick, directed=.55] (1, -2.5) to (1,2.5);
\node at (-1,3) {\tiny $1$};
\node at (1,3) {\tiny $1$};
\node at (-1,-3) {\tiny $1$};
\node at (1,-3) {\tiny $1$};
\end{tikzpicture}
};
\endxy - q^{-1}
\xy
(0,0)*{
\begin{tikzpicture}[scale=.3]
\draw [very thick, directed=.55] (0,-1) to (0,.75);
\draw [very thick, directed=.65] (0,.75) to [out=30,in=270] (1,2.5);
\draw [very thick, directed=.65] (0,.75) to [out=150,in=270] (-1,2.5); 
\draw [very thick, directed=.65] (-1, -2.5) to [out=90, in=210] (0,-1);
\draw [very thick, directed=.65] (1,-2.5) to [out =90, in=330] (0,-1);
\node at (-1,-3) {\tiny $1$};
\node at (1,-3) {\tiny $1$};
\node at (0.75,0) {\tiny $2$};
\node at (-1,3) {\tiny $1$};
\node at (1,3) {\tiny $1$};
\end{tikzpicture}
};
\endxy \Bigg),
\end{equation*} 
\begin{equation*}
\xy
(0,0)*{
\begin{tikzpicture}[scale=.3]
\draw [very thick, ->] (1,-1) to (-1,1);
\draw [very thick] (-1,-1) to (-0.25,-0.25);
\draw [very thick, ->] (0.25,0.25) to (1,1);
\node at (-1,-1.5) {\tiny $1$};
\node at (1,-1.5) {\tiny $1$};
\node at (1,1.45) {\tiny $1$};
\node at (-1,1.45) {\tiny $1$};
\end{tikzpicture}
};
\endxy  = q^{\frac{1}{n}-1} \Bigg( 
\xy
(0,0)*{
\begin{tikzpicture}[scale=.3]
\draw [very thick, directed=.55] (-1, -2.5) to (-1,2.5);
\draw [very thick, directed=.55] (1, -2.5) to (1,2.5);
\node at (-1,3) {\tiny $1$};
\node at (1,3) {\tiny $1$};
\node at (-1,-3) {\tiny $1$};
\node at (1,-3) {\tiny $1$};
\end{tikzpicture}
};
\endxy- q
\xy
(0,0)*{
\begin{tikzpicture}[scale=.3]
\draw [very thick, directed=.55] (0,-1) to (0,.75);
\draw [very thick, directed=.65] (0,.75) to [out=30,in=270] (1,2.5);
\draw [very thick, directed=.65] (0,.75) to [out=150,in=270] (-1,2.5); 
\draw [very thick, directed=.65] (-1, -2.5) to [out=90, in=210] (0,-1);
\draw [very thick, directed=.65] (1,-2.5) to [out =90, in=330] (0,-1);
\node at (-1,-3) {\tiny $1$};
\node at (1,-3) {\tiny $1$};
\node at (0.75,0) {\tiny $2$};
\node at (-1,3) {\tiny $1$};
\node at (1,3) {\tiny $1$};
\end{tikzpicture}
};
\endxy
\Bigg).
\end{equation*}
This apparently minor difference has huge impact, in particular, integrally. 
For example, while exterior webs fit perfectly \ochanged{with} tilting modules (see e.g. \cite{El-ladders-clasps} and \cite{AnStTu-cellular-tilting}), symmetric webs do not have any such interpretation known 
to the authors.

More details on minimal and integral exterior and red-green webs can 
be found in \cite{La-gln-webs}.
\end{Remark}

\ochanged{Following \cite{TuVaWe-super-howe}, we} colour the thick edges of symmetric webs by \ochanged{red} {\color{red}(symmetric = red)} and those of exterior webs by \ochanged{green} {\color{green}(exterior = green)} \ochanged{(since the first symmetric power is the same as the first exterior power, we use black for edges labelled $1$)}. These webs are related by the following relation 
\begin{gather}\label{Eq:RedGreen}
\xy
(0,0)*{
\begin{tikzpicture}[scale=.3]
\draw [very thick, directed=.55, red] (0,-1) to (0,.75);
\draw [very thick, directed=.65] (0,.75) to [out=30,in=270] (1,2.5);
\draw [very thick, directed=.65] (0,.75) to [out=150,in=270] (-1,2.5); 
\draw [very thick, directed=.65] (-1, -2.5) to [out=90, in=210] (0,-1);
\draw [very thick, directed=.65] (1,-2.5) to [out =90, in=330] (0,-1);
\node at (-1.25,-3) {\tiny $1$};
\node at (1.25,-3) {\tiny $1$};
\node at (0.75,0) {\tiny $2$};
\node at (-1.25,3) {\tiny $1$};
\node at (1.25,3) {\tiny $1$};
\end{tikzpicture}
};
\endxy 
=  - 
\xy
(0,0)*{
\begin{tikzpicture}[scale=.3]
\draw [very thick, directed=.55, green] (0,-1) to (0,.75);
\draw [very thick, directed=.65] (0,.75) to [out=30,in=270] (1,2.5);
\draw [very thick, directed=.65] (0,.75) to [out=150,in=270] (-1,2.5); 
\draw [very thick, directed=.65] (-1, -2.5) to [out=90, in=210] (0,-1);
\draw [very thick, directed=.65] (1,-2.5) to [out =90, in=330] (0,-1);
\node at (-1.25,-3) {\tiny $1$};
\node at (1.25,-3) {\tiny $1$};
\node at (0.75,0) {\tiny $2$};
\node at (-1.25,3) {\tiny $1$};
\node at (1.25,3) {\tiny $1$};
\end{tikzpicture}
};
\endxy + [2] \ 
\xy
(0,0)*{
\begin{tikzpicture}[scale=.3]
\draw [very thick, directed=.55] (-1, -2.5) to (-1,2.5);
\draw [very thick, directed=.55] (1, -2.5) to (1,2.5);
\node at (-1,3) {\tiny $1$};
\node at (1,3) {\tiny $1$};
\node at (-1,-3) {\tiny $1$};
\node at (1,-3) {\tiny $1$};
\end{tikzpicture}
};
\endxy.
\end{gather}
\changed{This is the \textit{red-green dumbbell relation}. The (purely symmetric) \textit{dumbbell relation} is implicit. Namely, by definition, one annihilates the green clasp on $n+1$ strands, and then 
one uses \cite[Corollary 2.13]{TuVaWe-super-howe} to replace all green edges by red edges, which gives a purely symmetric relation. (Note that this relation is only imposed on endomorphisms of $1_\uparrow^{\otimes (n+1)}$ only.)}

\begin{Remark}
\changed{In exterior webs, there is a relation that annihilates a green edge of too thick label:
\begin{gather*}
\xy
(0,0)*{
\begin{tikzpicture}[scale=.3]
\draw [very thick, directed=.55, green] (0,-1) to (0,.75);
\node at (0.75,0) {\tiny $a$};
\end{tikzpicture}
};
\endxy 
=
0,
\end{gather*}
where $a>n$. This is needed to annihilate the so-called antisymmetrizer on $n$ strands, a certain expression in crossings on $n$ strands. For example, for $n=1$
\begin{gather*}
0=
\xy
(0,0)*{\begin{tikzpicture}[scale=.3]
\draw [very thick, ->] (-1,-1) to (-1,1);
\draw [very thick, ->] (1,-1) to (1,1);
\node at (-1,-1.5) {\tiny $1$};
\node at (1,-1.5) {\tiny $1$};
\node at (1,1.45) {\tiny $1$};
\node at (-1,1.45) {\tiny $1$};
\end{tikzpicture}};
\endxy
-
\xy
(0,0)*{\begin{tikzpicture}[scale=.3]
\draw [very thick, ->] (-1,-1) to (1,1);
\draw [very thick] (1,-1) to (0.25,-0.25);
\draw [very thick, ->] (-0.25,0.25) to (-1,1);
\node at (-1,-1.5) {\tiny $1$};
\node at (1,-1.5) {\tiny $1$};
\node at (1,1.45) {\tiny $1$};
\node at (-1,1.45) {\tiny $1$};
\end{tikzpicture}};
\endxy
=
\xy
(0,0)*{\begin{tikzpicture}[scale=.3]
\draw [very thick, ->] (-1,-1) to (-1,1);
\draw [very thick, ->] (1,-1) to (1,1);
\node at (-1,-1.5) {\tiny $1$};
\node at (1,-1.5) {\tiny $1$};
\node at (1,1.45) {\tiny $1$};
\node at (-1,1.45) {\tiny $1$};
\end{tikzpicture}};
\endxy
-\left(
\xy
(0,0)*{\begin{tikzpicture}[scale=.3]
\draw [very thick, ->] (-1,-1) to (-1,1);
\draw [very thick, ->] (1,-1) to (1,1);
\node at (-1,-1.5) {\tiny $1$};
\node at (1,-1.5) {\tiny $1$};
\node at (1,1.45) {\tiny $1$};
\node at (-1,1.45) {\tiny $1$};
\end{tikzpicture}};
\endxy
-
\xy
(0,0)*{
\begin{tikzpicture}[scale=.3]
\draw [very thick, directed=.55,green] (0,-1) to (0,.75);
\draw [very thick, directed=.65] (0,.75) to [out=30,in=270] (1,2.5);
\draw [very thick, directed=.65] (0,.75) to [out=150,in=270] (-1,2.5); 
\draw [very thick, directed=.65] (-1, -2.5) to [out=90, in=210] (0,-1);
\draw [very thick, directed=.65] (1,-2.5) to [out =90, in=330] (0,-1);
\node at (-1,-3) {\tiny $1$};
\node at (1,-3) {\tiny $1$};
\node at (0.75,0) {\tiny $2$};
\node at (-1,3) {\tiny $1$};
\node at (1,3) {\tiny $1$};
\end{tikzpicture}};
\endxy
\right),
\end{gather*}
where we set $q=1$ for simplicity, and scaled by $2$. Thus, the exterior $n=2$ edge needs to be zero. This is not true for symmetric webs, but by construction, the dumbbell relation implies the vanishing of the antisymmetrizer (which, via explosion, then has further consequences for the symmetric calculus on arbitrary thick strands).}

\changed{For the reader familiar with \cite{TuVaWe-super-howe}, the category of upwards symmetric webs in this paper then corresponds to the red webs only subcategory in that paper (where the dumbbell relation holds automatically).}
\end{Remark}

\subsection{A basis for thin webs}\label{S:Basis}

In \cite{El-ladders-clasps}, Elias constructed a basis for exterior $\mathfrak{sl}_{n}$-webs whose elements are called double light ladders. We will be considering only the basis of $\text{End} _{\mathbf{EWeb}_{\uparrow, \downarrow}(\mathfrak{sl}_n)}(1_\uparrow^{ \otimes k}), \ k \geq 1$. 
(The notation means exterior webs, see above.) We modify the \ochanged{latter} by performing an algorithm leading to a basis of $\text{End} _{\mathbf{SWeb}_{\uparrow}(\mathfrak{gl}_n)}(1_\uparrow^{\otimes k})$.

The three definitions from below will be used for both exterior and symmetric webs.

\begin{Definition} 
A web with $1$-labelled boundary edges is called a \textit{local dumbbell} if it contains only dumbbells and identity components.
\end{Definition}

\begin{Lemma}\label{L:Remove}
Let $X$ be a local dumbbell. Then $X$ can be written as a linear combination of webs containing only 2-dumbbells (dumbbells whose thickest label is 2) and identity components. 
\end{Lemma}

\begin{proof}
This follows from the $\mathfrak{gl}_n$-projector recursion in \autoref{JW-recursion}. 
\end{proof}

\begin{Definition}
We call the linear combination from \ochanged{\autoref{L:Remove}} above \textit{the thin expression} of the local dumbbell $X$. 
\end{Definition}

\begin{Definition}
We say that a local dumbbell $X$ is \textit{smaller or equal} to the local dumbbell $Y$, we write $ X \leq Y$, if the total number of the 2-dumbbells appearing in the thin expression of $X$ is smaller or equal to the total number of the 2-dumbbells appearing in the thin expression of $Y$. 
\end{Definition}

Now, let $B_1,\dots,B_l$ be the basis elements of $\text{End} _{\mathbf{EWeb}_{\uparrow,\downarrow}(\mathfrak{sl}_n)}(1^{\otimes k})$ from \cite{El-ladders-clasps}. We perform the following two steps:\\
\medskip
$\mathbf{Step  \ 1.}$ In each of the webs $B_1,\dots,B_l$ we change every cup and cap of the form 
\begin{equation*}
\xy
(0,0)*{
\begin{tikzpicture} [scale=1]
\draw[very thick,rdirected=-.95, ] (0,1) to [out=270,in=180] (.25,.5) to [out=0,in=270] (.5,1);
\node at (.25,.25) {};
\node at (0,1.25) {\tiny $n-k$};
\node at (.5,1.25) {\tiny $k$};
\end{tikzpicture}
}
\endxy 
\quad , \quad 
\xy
(0,0)*{
\begin{tikzpicture} [scale=1]
\draw[very thick,rdirected=-.95, ] (0,0) to [out=90,in=180] (.25,.5) to [out=0,in=90] (.5,0);
\node at (.25,.75) {};
\node at (0,-.25) {\tiny $n-k$};
\node at (.5,-.25) {\tiny $k$};
\end{tikzpicture}
}
\endxy   \quad  
\end{equation*}
to splits and merges of the form
\begin{equation*}
\xy
(0,0)*{
\begin{tikzpicture}[scale=.3]
\draw [very thick, directed=.55,green] (0,-1) to (0,.75);
\draw [very thick, directed=.65,green] (0,.75) to [out=30,in=270] (1,2.5);
\draw [very thick, directed=.65,green] (0,.75) to [out=150,in=270] (-1,2.5); 
\node at (0, -1.5) {\tiny $n$};
\node at (-1,3) {\tiny $n-k$};
\node at (1,3) {\tiny $k$};
\end{tikzpicture}
};
\endxy \quad , \quad
\xy
(0,0)*{
\begin{tikzpicture}[scale=.3]
\draw [very thick, directed=.55,green] (0, .75) to (0,2.5);
\draw [very thick, directed=.45,green] (1,-1) to [out=90,in=330] (0,.75);
\draw [very thick, directed=.45,green] (-1,-1) to [out=90,in=210] (0,.75); 
\node at (0, 3) {\tiny $n$};
\node at (-1,-1.5) {\tiny $n-k$};
\node at (1,-1.5) {\tiny $k$};
\end{tikzpicture}
};
\endxy,
\end{equation*}
where the $n$-labelled edge has no self-intersections. Similarly we define the operation for the mirrored webs. In order for this operation to be well-defined we take the construction of webs as morphisms in a monoidal category given by \ochanged{generators-relations}. Now, whenever this would lead to crossings of webs, we pull the $n$-labelled edges under and outside to the right like below:
\[
\xy
(0,0)*{
\begin{tikzpicture}[scale=.3]
\draw [very thick,green] (1,-1) to [out=90,in=0] (0,.5);
\draw [very thick, rdirected=-.9,green] (-1,-1) to [out=90,in=180] (0,.5); 
\draw [very thick, directed=.45,green] (-3,-1) to [out=90,in=210] (1,3); 
\draw [very thick, directed=.45,green] (5,-1) to [out=90,in=330] (1,3);
\draw [very thick, directed=.55,green] (1,3) to (1,6);
\node at (-1,-1.5) {\tiny $n-k$};
\node at (1,-1.5) {\tiny $k$};
\end{tikzpicture}
};
\endxy\quad \mapsto \quad 
\xy
(0,0)*{
\begin{tikzpicture}[scale=.3]
\draw [very thick, directed=.75,green, ](0, .75) to (4,6);
\draw [very thick, directed=.45,green] (1,-1) to [out=90,in=330] (0,.75);
\draw [very thick, directed=.45,green] (-1,-1) to [out=90,in=210] (0,.75); 
\draw [very thick, directed=.45,green] (-3,-1) to [out=90,in=210] (1,3); 
\draw [very thick, directed=.45,green,crossline] (5,-1) to [out=90,in=330] (1,3);
\draw [very thick, directed=.55,green] (1,3) to (1,6);
\node at (0, 3) {\tiny $n$};
\node at (2, 6.5) {\tiny $n$};
\node at (-1,-1.5) {\tiny $n-k$};
\node at (1,-1.5) {\tiny $k$};
\end{tikzpicture}
};
\endxy.
\]
Note that here we are making a choice. The reader might wonder why this is well-defined and what happens if for example we let
\[
\xy
(0,0)*{
\begin{tikzpicture}[scale=.3]
\draw [very thick,green] (1,-1) to [out=90,in=0] (0,.5);
\draw [very thick, rdirected=-.9,green] (-1,-1) to [out=90,in=180] (0,.5); 
\draw [very thick, directed=.45,green] (-3,-1) to [out=90,in=210] (1,3); 
\draw [very thick, directed=.45,green] (5,-1) to [out=90,in=330] (1,3);
\draw [very thick, directed=.55,green] (1,3) to (1,6);
\node at (-1,-1.5) {\tiny $n-k$};
\node at (1,-1.5) {\tiny $k$};
\end{tikzpicture}
};
\endxy\quad \mapsto \quad 
\xy
(0,0)*{
\begin{tikzpicture}[scale=.3]
\draw [very thick, directed=.75,green, ](0, .75) to (-3,3);
\draw [very thick, directed=.75,green, ](0,5) to (5,6);
\draw [very thick,green,] (-3,3) to [out=90,in=180] (0,5); 
\draw [very thick, directed=.45,green] (1,-1) to [out=90,in=330] (0,.75);
\draw [very thick, directed=.45,green] (-1,-1) to [out=90,in=210] (0,.75); 
\draw [very thick, directed=.45,green,crossline] (-3,-1) to [out=90,in=210] (1,3); 
\draw [very thick, directed=.45,green,crossline] (5,-1) to [out=90,in=330] (1,3);
\draw [very thick, directed=.55,green, crossline] (1,3) to (1,6);
\node at (0, 3) {\tiny $n$};
\node at (2, 6.5) {\tiny $n$};
\node at (-1,-1.5) {\tiny $n-k$};
\node at (1,-1.5) {\tiny $k$};
\end{tikzpicture}
};
\endxy.
\]
We point out that all such choices would lead to the same web due to the ambient isotopies such as zigzag relations, Reidemeister 2 and 3 moves, pitchfork relations, cup and cap slides etc. Moreover, if there are multiple caps, then we choose an arbitrary order of their end points noting that 
\[
\begin{tikzpicture}[anchorbase,scale=0.7]
\draw[very thick,directed=0.99,green,crossline]
(1,0) node[below, black]{\tiny $n$}to (0,1)node[above, black]{\tiny $n$};
\draw[very thick,directed=0.99,green,crossline]
(0,0) node[below, black]{\tiny $n$}to (1,1)node[above, black]{\tiny $n$};
\end{tikzpicture}
=
\begin{tikzpicture}[anchorbase,scale=0.7]
\draw[very thick,directed=0.99,green,crossline]
(1,0) node[below, black]{\tiny $n$}to (1,1)node[above, black]{\tiny $n$};
\draw[very thick,directed=0.99,green,crossline]
(0,0) node[below, black]{\tiny $n$}to (0,1)node[above, black]{\tiny $n$};
\end{tikzpicture}
=
\begin{tikzpicture}[anchorbase,scale=0.7]
\draw[very thick,directed=0.99,green,crossline]
(0,0) node[below, black]{\tiny $n$}to (1,1)node[above, black]{\tiny $n$};
\draw[very thick,directed=0.99,green,crossline]
(1,0) node[below, black]{\tiny $n$}to (0,1)node[above, black]{\tiny $n$};
\end{tikzpicture}, \]
modulo lower \ochanged{terms (see \cite{La-gln-webs}).}\\ 
$\mathbf{Step \  2.}$ We change the colours from green to red. 

\begin{Lemma} \label{ext thin basis} 
If we apply Step 1 from above to \ochanged{the} double ladder basis, then we obtain a basis for $\mathrm{End}_{\mathbf{EWeb}_\uparrow(\mathfrak{gl}_n)}(1^{\otimes k}).$
\end{Lemma}

\begin{proof} 
Let $B_1',\dots,B_l'$ be the $\mathfrak{gl}_n$-webs obtained from the $\mathfrak{sl}_n$-webs $B_1,\dots,B_l$ (double ladders) after Step 1. Note that the reason why we are performing Step 1 is because for $\mathfrak{sl}_n$-webs we do not have edges of labels $\geq n$ and exterior $\mathfrak{sl}_n$-webs are a quotient of exterior $\mathfrak{gl}_n$-webs by 
\[
\begin{tikzpicture}[anchorbase]
\draw[very thick,->] (0,1) to node[right]{$k$} (0,0);
\end{tikzpicture}
=
\begin{tikzpicture}[anchorbase]
\draw[very thick,->] (0,0) to node[right]{$n-k$} (0,1);
\end{tikzpicture}.
\]
Now, the webs $B_1',\dots,B_l'$ are linearly independent. Indeed, if they were linearly dependent, after passing to $\mathfrak{sl}_n$-webs (by imposing the relation above) we would get that $B_1,\dots,B_l$ are linearly dependent, contradicting the fact that they form a basis. It remains to show that $B_1',\dots,B_l'$ span $\mathrm{End}_{\mathbf{EWeb}_\uparrow(\mathfrak{gl}_n)}(1^{\otimes k}).$ This is a word for word \ochanged{adaptation} of the proof of \cite[Theorem 2.39]{El-ladders-clasps}. Indeed, Elias points out that \ochanged{the} proof would apply to any class of morphism satisfying \cite[Proposition 2.42]{El-ladders-clasps}. This finishes the proof
\ochanged{(see \cite{La-gln-webs} for more details).}
\end{proof}

\begin{Proposition}\label{thin basis} 
If we apply Step 1 and Step 2 to the double ladder basis we obtain a basis $\{\tilde{B}_{1},\dots,\tilde{B}_l \}$ for $\mathrm{End}_{\mathbf{SWeb}_{\uparrow}(\mathfrak{gl}_n)}(1_\uparrow^{\otimes k}).$
\end{Proposition}

\begin{proof} 
Let $\{B_1',\dots,B_l' \}$ be the basis for $\text{End}_{\mathbf{EWeb}_{\uparrow}(\mathfrak{gl}_n)}(1^{\otimes k})$ from the lemma above. 
We use the following strategy:\\
(a) First, we replace each $B_{i}'$ by a linear combination of webs containing only 2-dumbbells and identities locally. \\
(b) Second, we change the colours of the resulting webs.\\
(c) In these steps it is crucial that we have an order on webs given by counting 2-dumbbells, so both of the above will produce a triangular change-of-basis matrix.

We now give the details. By construction, each of the webs $B_{i}'$ is locally built from dumbbells and identities. Now we use the exterior version of $\mathfrak{gl}_n$-projector recursion from \autoref{JW-recursion} to write every $B_i'$ as a linear combination of webs built from only $2$-dumbbells and identities. Without loss of generality we assume that $B_1' \leq \dots\leq B_l'$. Note that there is a biggest web in the $\mathfrak{gl}_n$-projector formula (\autoref{JW-recursion}), namely the rightmost web. After performing all possible recursive steps, we will have one web which is the overall biggest in the thin expression of each $k$-dumbbell and consequently of each $B_i'$. Let $W_i$ be the biggest web appearing in the thin expression of $B_i'$. From the discussion above, the coefficient of $W_i$ in the thin expression of $B_i'$ is invertible. Inductively one concludes that the change-of-basis matrix from $B_i'$ to $W_i$ is triangular with invertible diagonal entries. This means that $W_1,\dots,W_l$ are linearly independent and consequently they form a basis for  $\text{End} _{\mathbf{EWeb}_{\uparrow}(\mathfrak{gl}_n)}(1^{\otimes k})$. Also, note that $W_1 \leq\dots\leq W_l$. After applying Step 2, i.e. changing the colours of the webs $W_1,\dots,W_l$, we denote the new webs by $\tilde{W}_1,\dots,\tilde{W}_l$. We claim that they are linearly independent. Indeed, using the relation
\[ 
\xy
(0,0)*{
\begin{tikzpicture}[scale=.3]
\draw [very thick, directed=.55, red] (0,-1) to (0,.75);
\draw [very thick, directed=.65] (0,.75) to [out=30,in=270] (1,2.5);
\draw [very thick, directed=.65] (0,.75) to [out=150,in=270] (-1,2.5); 
\draw [very thick, directed=.65] (-1, -2.5) to [out=90, in=210] (0,-1);
\draw [very thick, directed=.65] (1,-2.5) to [out =90, in=330] (0,-1);
\node at (-1.25,-3) {\tiny $1$};
\node at (1.25,-3) {\tiny $1$};
\node at (0.75,0) {\tiny $2$};
\node at (-1.25,3) {\tiny $1$};
\node at (1.25,3) {\tiny $1$};
\end{tikzpicture}
};
\endxy 
=  - 
\xy
(0,0)*{
\begin{tikzpicture}[scale=.3]
\draw [very thick, directed=.55, green] (0,-1) to (0,.75);
\draw [very thick, directed=.65] (0,.75) to [out=30,in=270] (1,2.5);
\draw [very thick, directed=.65] (0,.75) to [out=150,in=270] (-1,2.5); 
\draw [very thick, directed=.65] (-1, -2.5) to [out=90, in=210] (0,-1);
\draw [very thick, directed=.65] (1,-2.5) to [out =90, in=330] (0,-1);
\node at (-1.25,-3) {\tiny $1$};
\node at (1.25,-3) {\tiny $1$};
\node at (0.75,0) {\tiny $2$};
\node at (-1.25,3) {\tiny $1$};
\node at (1.25,3) {\tiny $1$};
\end{tikzpicture}
};
\endxy + [2] \ 
\xy
(0,0)*{
\begin{tikzpicture}[scale=.3]
\draw [very thick, directed=.55] (-1, -2.5) to (-1,2.5);
\draw [very thick, directed=.55] (1, -2.5) to (1,2.5);
\node at (-1,3) {\tiny $1$};
\node at (1,3) {\tiny $1$};
\node at (-1,-3) {\tiny $1$};
\node at (1,-3) {\tiny $1$};
\end{tikzpicture}
};
\endxy
\]
we can write each symmetric web $\tilde{W_i}$ as a linear combination of exterior webs. Recall that by definition, the number of 2-dumbbell components in a web $\tilde{W_i}$ is the same as the number of 2-dumbbell components in $W_i$. So, $\tilde{W_1} \leq \dots \leq \tilde{W}_l$. Now, from the above relation we can write each $\tilde{W_i}$ as a linear combination of exterior basis webs $W_j$, such that $W_j \leq W_i$. Moreover, again from the same relation we have that the coefficient of $W_i$ in this linear combination is $(-1)^a$, where $a$ is the number of 2-dumbbell components of $W_i$. This means that the change-of-basis matrix from symmetric to exterior webs is triangular with determinant $\pm 1$. This concludes the proof that $\tilde{W_1} , \dots , \tilde{W}_l$ are linearly independent. Since the algorithm sends exterior web generators to symmetric web \ochanged{generators}, then $\{\tilde{W_1}, \dots , \tilde{W}_l\}$ is a basis for $\text{End} _{\mathbf{SWeb}_{\uparrow}(\mathfrak{gl}_n)}(1^{\otimes k})$. Finally, note that by construction, $\tilde{W}_i$ is the biggest web appearing in the thin expression of $\tilde{B}_i$. Thus, like before, the change-of-basis matrix from $\tilde{W}_i$ to $\tilde{B}_i$ is triangular with non-zero diagonal entries, meaning that $\{\tilde{B}_1,\dots,\tilde{B}_l\}$ is a basis for $\text{End} _{\mathbf{SWeb}_{\uparrow}(\mathfrak{gl}_n)}(1^{\otimes k})$. 
\end{proof}

\begin{Example}
\changed{Let $n=k=3$. Then Elias' light ladder basis for $\text{End} _{\mathbf{EWeb}_{\uparrow,\downarrow}(\mathfrak{sl}_n)}(1^{\otimes k})$ is
\begin{gather*}
\xy
(0,0)*{\includegraphics[height=6cm]{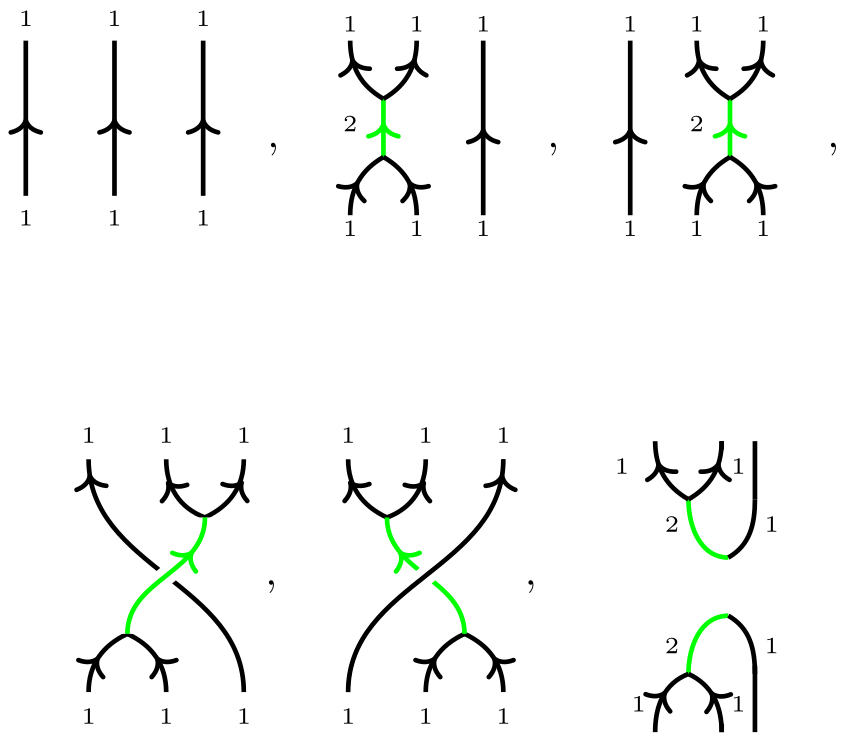}};
\endxy
\end{gather*}
given in our drawing convention (for the special linear groups, so-called tags included). Strictly speaking Elias would use so-called \textit{neutral ladders} instead of crossings, but crossings work equally well. Then
\begin{gather*}
\includegraphics[height=6cm]{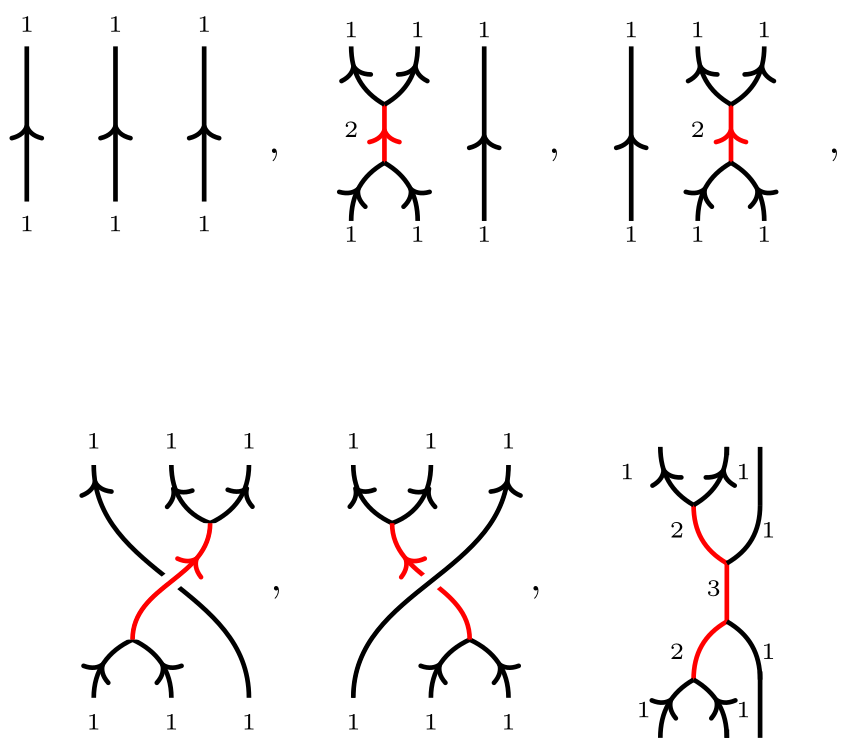}
\end{gather*}
is the symmetric web basis of $\text{End} _{\mathbf{SWeb}_{\uparrow}(\mathfrak{gl}_n)}(1^{\otimes k})$. More details can be found in \cite{La-gln-webs}.}
\end{Example}

\section{Equivalence Theorem}

We now turn to representation theory and we explain how webs describe the category of finite-dimensional modules over the quantum general linear group. 

\subsection{The quantum general linear group and its idempotented form}

We will use the notation $\epsilon_{i}= (0,\dots,1,\dots,0)  \in \mathbb{Z}^m$,  where 1 is in the $i$-th coordinate. For $i = 1,\dots,m-1$, let $ \alpha_{i} = \epsilon_{i} - \epsilon_{i+1} = (0,\dots,1,-1,\dots,0) \in \mathbb{Z}^m,$. Moreover, recall that the Euclidean inner product on $\mathbb{Z}^m$ is given by $(\epsilon_{i},\epsilon_{j})=\delta_{i,j}$.

\begin{Definition}\label{quantum gln} 
For $m \in \mathbb{Z} _{>1}$, the \textit{quantum general linear group} $\text{U}_{q}(\mathfrak{gl}_m)$ is the associative, unital $\Bbbk(q)$-algebra generated by $L_{i}$ and $L_{i}^{-1}$, for $i = 1,\dots,m$, and $E_{i}, F_{i}$, for $i=1,\dots,m-1$, subject to the following relations (for all admissible indices):
\begin{center}
$L_{i_{1}}L_{i_{2}} =L_{i_{2}}L_{i_{1}}, $
\medskip

$ L_{i}L_{i}^{-1}= L_{i}^{-1} L_{i}=1,$ 
\medskip

$ L_{i_1}E_{i_2}= q^{( \epsilon_{i_1}, \alpha_{i_2})}E_{i_2} L_{i_1}$, 

\medskip

$ L_{i_1}F_{i_2}= q^{-( \epsilon_{i_1}, \alpha_{i_2})}F_{i_2} L_{i_1}$, 

\medskip

$E_{i_1}F_{i_2}- F_{i_2}E_{i_1}= \delta_{i_1,i_2} \frac{ L_{i_1}L_{i_{1}+1}^{-1} - L_{i_1}^{-1}L_{i_{1}+1}}{q-q^{-1}}
$ , 
\medskip

$E_{i_1}^2E_{i_2}- [2]E_{i_1}E_{i_2}E_{i_1} + E_{i_2}E_{i_1}^2 = 0,   $ if $ |i_1- i_2|=1$, $E_{i_1}E_{i_2}=E_{i_2}E_{i_1} $, else, 

\medskip

$F_{i_1}^2F_{i_2}- [2]F_{i_1}F_{i_2}F_{i_1} + F_{i_2}F_{i_1}^2 = 0,   $ if $ |i_1- i_2|=1$, $F_{i_1}F_{i_2}=F_{i_2}F_{i_1} $, else.

\end{center}
\medskip
The last two relations are called the (quantum) \textit{Serre relations}.
\end{Definition}

\begin{Definition}\label{sln} 
For $m \in \mathbb{Z}_{>1}$ the \textit{quantum special linear algebra}  $\text{U}_{q}(\mathfrak{sl}_m)$ is the subalgebra of $\text{U}_{q}(\mathfrak{gl}_m)$ generated by $E_{i}, F_{i}, K_i=L_iL_{i+1}^{-1}$ and $K_i^{-1}=L_{i+1}L_i^{-1}$ for $i=1,\dots,m-1$.
\end{Definition}

\label{HA structure} We can endow $\text{U}_{q}(\mathfrak{gl}_m)$ with a Hopf algebra structure (which will be used throughout this paper), where the coproduct $\Delta$ is given by 

\begin{center}
\begin{gather*}
\Delta(E_{i}) = E_{i} \otimes L_{i}L_{i+1}^{-1} + 1 \otimes E_{i},   
\end{gather*}
$\Delta(F_{i})= F_{i} \otimes 1 +  L_{i}^{-1} L_{i+1} \otimes F_{i} $ and 

\bigskip

$\Delta ( L_{i}^{\pm 1}) =  L_{i}^{\pm 1}  \otimes  L_{i}^{ \pm 1} .$

\end{center}

\bigskip

\begin{flushleft}

The antipode and counit are given by 
\end{flushleft}

\bigskip

\begin{center}
$S(E_{i})= - E_{i}L_{i}^{-1}L_{i+1}$, $S(F_{i})=-L_{i}L_{i+1}^{-1}F_{i}$,  $S(L_{i}^{\pm 1})=L_{i}^{\mp 1} $\ochanged{,}
\medskip

$\epsilon(E_{i})= \epsilon(F_{i})=0$ and $\epsilon(L_{i}^{\pm 1})=1.$
\end{center}
The subalgebra $\text{U}_{q}(\mathfrak{sl}_m) $ inherits the Hopf algebra structure from $\text{U}_{q}(\mathfrak{gl}_m) $. Recall that 
\ochanged{it follows from the Hopf algebra structure that}
the existence of the trivial representation. It also allows us to extend actions to tensor products of representations.
Now, note that the antipode is invertible and the condition below is satisfied.
\begin{gather*}
S^2(h)=K_{2\rho}hK^{-1}_{2\rho},
\end{gather*}
for every $h \in \text{U}_{q}(\mathfrak{gl}_m)$, where $K_{2\rho}=K_1^{m-1}K_2^{m-2} \dots K_{m-1}$.

The following representation of $\text{U}_{q}(\mathfrak{gl}_m) $ is of particular interest for us.

\begin{Definition}\label{stand-rep}
The \textit{standard representation} $\Bbbk_q^m$ of $\text{U}_{q}(\mathfrak{gl}_m)$ is the $m$-dimensional $\Bbbk(q)$-vector space with basis $b_1,\dots,b_m$, where 

\begin{center}
$L_{i}b_i=qb_i, \ L_i^{-1}b_i=q^{-1}b_i, \ L_i^{\pm 1} b_j=b_j, \  j \neq i,$

\medskip

$E_{i-1}b_i=b_{i-1}, \ E_{i}b_{j}=0,$ if $i \neq j-1$,

\medskip

$F_{i}b_{i}=b_{i+1}, \ F_ib_j=0,$ if $ i \neq j$.
\end{center}
Note that this is the quantum version of the standard (vector) representation of $\text{GL}_m$.
\end{Definition}

We will later need \textit{Lusztig's idempotented form} of $\text{U}_{q}(\mathfrak{gl}_m)$. The idea is to add to $\text{U}_{q}(\mathfrak{gl}_m)$ an orthogonal system of idempotents in order to view this algebra as a category. Recall that (integral) $\mathfrak{gl}_m$-weights are tuples $\vec{k}=(k_1,\dots,k_m)$ in $\mathbb{Z}^m$. We define the \textit{divided powers}
\begin{center}
$F_{i}^{(j)}= \frac{F_{i}^j}{[j]!}$ and $E_{i}^{(j)}= \frac{E_{i}^j}{[j]!}$.
\end{center}

\begin{Definition}\label{integral idemp form} 
The \textit{category} $\dot{\text{U}}_{q}(\mathfrak{gl}_m)$ is the $\mathbb{Z}[q, q^{-1}]$-linear category with objects the $\mathfrak{gl}_m$-weights $\vec{k} =(k_1, \dots , k_m) \in \mathbb{Z}^m$. The identity morphism of $\vec{k}$ is denoted by $\mathbf{1}_{\vec{k}}$ and the morphism spaces are denoted by $\mathbf{1}_{\vec{l} } \dot{\text{U}}_{q}(\mathfrak{gl}_m) \mathbf{1} _{ \vec{k}} $. The morphisms are generated by 
$E_i^{(r)} \mathbf{1} _{ \vec{k}} \in \mathbf{1}_{\vec{k}+r\alpha_i } \dot{\text{U}}_{q}(\mathfrak{gl}_m) \mathbf{1} _{ \vec{k}} $ and $F_i^{(r)} \mathbf{1} _{ \vec{k}} \in \mathbf{1}_{\vec{k}-r\alpha_i } \dot{\text{U}}_{q}(\mathfrak{gl}_m) \mathbf{1} _{ \vec{k}} $, $r \in \mathbb{N}$ (here $ \alpha_{i} = (0,\dots,1,-1,\dots,0) \in \mathbb{Z}^m$ as before). Note that 
$\text{Hom}(\vec{k}, \vec{l})=0$ unless $\sum k_i =\sum l_i$.
These generating morphisms satisfy the following relations
\begin{gather*}
F_{i_1}^{(j_1)}F_{i_2}^{(j_2)}\mathbf{1}_{\vec{k}}=F_{i_2}^{(j_2)}F_{i_1}^{(j_1)}\mathbf{1}_{\vec{k}}, \ \text{if} \ |i_1-i_2|>1,
\end{gather*}
\begin{gather*}
F_{i_1}^{(j_1)}E_{i_2}^{(j_2)}\mathbf{1}_{\vec{k}}=E_{i_2}^{(j_2)}F_{i_1}^{(j_1)}\mathbf{1}_{\vec{k}}, \ \text{if} \ |i_1-i_2| > 0,
\end{gather*}

\begin{equation} \label{serre rel in dot}
F_{i_1} F_{i_2} F_{i_1} \mathbf{1}_{\vec{k}} =(F_{i_1}^{(2)}F_{i_2} + F_{i_2} F_{i_1}^{(2)}) \mathbf{1}_{\vec{k}}, \  \text{if} \ |i_1 -i_2|=1, 
\end{equation}

\begin{gather*}
F_{i}^{(j_1)}F_{i}^{(j_2)}\mathbf{1}_{\vec{k}}= \begin{bmatrix} j_1 + j_2 \\ j_1 \end{bmatrix}  F_{i}^{(j_1+j_2)}\mathbf{1}_{\vec{k}},
\end{gather*}

\begin{gather*}
E_{i}^{(j_2)}F_{i}^{(j_1)}\mathbf{1}_{\vec{k}}=\sum \limits_{j'}\begin{bmatrix} k_i-j_1-k_{i+1}+j_2 \\ j' \end{bmatrix} F_{i}^{(j_1-j')}E_{i}^{(j_2-j')}\mathbf{1}_{\vec{k}},
\end{gather*} 

\begin{gather*}
E_i^{(j)} \mathbf{1}_{\vec{k}}=\mathbf{1}_{\vec{k}+j\alpha_i}E_i^{(j)}, \  F_i^{(j)} \mathbf{1}_{\vec{k}}=\mathbf{1}_{\vec{k}-j\alpha_i}F_i^{(j)}
\end{gather*}
and their analogues for $E_{i}^{(j)}$ (when the weight is clear from the context we write  $E_{i}^{(j)}$ instead of  $E_{i}^{(j)}\mathbf{1}_{\vec{k}}$ in the first, third and fourth relation) Relation \autoref{serre rel in dot} is called the \textit{Serre relation}. 
\end{Definition}

\subsection{Quantum symmetric Howe functors}\label{Howe duality}

Let $\text{Sym}^k_q\Bbbk_q^n$ be the simple $\text{U}_{q}(\mathfrak{gl}_m)$-module of highest weight $(k,0,\dots,0).$ It has a basis given by
\begin{equation*}
\changed{\{y_{j_1} \cdot \dots  \cdot y_{j_k}=y_{j_1} \otimes \dots  \otimes y_{j_k}, \ \text{for} \ 1 \leq  j_1 \leq \dots \leq  j_k \leq n,\}}
\end{equation*}
where $\{y_1, \dots ,y_n\}$ is a basis of $\Bbbk_q^n$ (see \cite{BeZw-quantum-braided-sym-ext}) and its dimension is $\begin{pmatrix}
n +k-1 \\k
\end{pmatrix} $.
The following is a consequence of quantum symmetric Howe duality (Theorem 2.6 in \cite{RoTu-symmetric-howe}).

\begin{Corollary}\label{howe-func} 
(\cite[Theorem 2.6]{RoTu-symmetric-howe}.) There exists a functor
\[\Phi_{m,n} :\dot{\mathrm{U}} _{q}(\mathfrak{gl}_m) \rightarrow \mathrm{U}_{q}( \mathfrak{gl}_n)\text{-}\mathbf{fdMod}_S,
\]
which sends a $\mathfrak{gl}_m$-weight $\vec{k}=(k_1,\dots,k_m) \in \mathbb{Z}_{\geq0}^m$ to the $ \text{U}_{q}(\mathfrak{gl}_n)$-module $\text{Sym}_q^{k_1} \otimes \cdots \otimes \text{Sym}_q^{k_m }$ and morphisms $X \in \mathbf{1}_{\vec{l}} \dot{\mathrm{U}}_{q}(\mathfrak{gl}_m)\mathbf{1}_{\vec{k}} $ to $ f_{\vec{k}}^{\vec{l}} (X)$. The $\mathfrak{gl}_m$-weights with negative entries are sent to the zero module. This functor is moreover surjective on Hom-spaces as in equation \autoref{eq-symmor}:
\begin{equation} \label{eq-symmor}
f_{\vec{k}}^{\vec{l}}: \mathbf{1}_{\vec{l}} \dot{\text{U}}_{q}(\mathfrak{gl}_m)\mathbf{1}_{\vec{k}} \longrightarrow \text{Hom}_{\text{U}_{q}(\mathfrak{gl}_n)} (\text{Sym}_q^{k_1}  \otimes \cdots \otimes  \text{Sym}_q^{k_m}, \text{Sym}_q^{l_1} \otimes \cdots \otimes \text{Sym}_q^{l_m} )
\end{equation}
for any two $\vec{k}, \vec{l} \in \mathbb{Z}_{\geq 0}^m$ such that $\sum \limits_{i=0}^{m} k_{i} = \sum \limits _{i=0}^m l_{i}$, meaning that it is full.
\end{Corollary}

By dualising \cite[Theorem 2.6]{RoTu-symmetric-howe} we obtain \changed{the following. (More precisely, the dualization takes place in the source category by replacing weights with their negatives; under this identification, the contravariant nature of dualization aligns with the direction of the previous argument, so the same reasoning applies.)}

\begin{Corollary} 
There exists a functor 
\[\Phi_{m.n}^* : \dot{\mathrm{U}}_{q}(\mathfrak{gl}_m) \rightarrow \mathrm{U}_{q}( \mathfrak{gl}_n)\text{-}\mathbf{fdMod}_{S^*},
\]
which sends a $\mathfrak{gl}_m$-weight $\vec{k}=(-k_1,\dots,-k_m) \in \mathbb{Z}_{\leq 0}^m$ to the $ \text{U}_{q}(\mathfrak{gl}_n)$-module $(\text{Sym}_{q}^{k_m})^* \otimes \cdots \otimes (\text{Sym}_{q}^{k_1 })^*$ and morphisms $X \in \mathbf{1}_{\vec{l}} \dot{\mathrm{U}}_{q}(\mathfrak{gl}_m)\mathbf{1}_{\vec{k}} $ to $ (f_{\vec{k}}^{\vec{l}})^* (X)$. This functor is surjective on Hom-spaces analogous to those in \autoref{eq-symmor}.
\end{Corollary}

\changed{We refer to $\Phi_{m,n} $ and $\Phi_{m,n}^*$ as the \textit{q-symmetric Howe functors}.}

\subsection{The functor \texorpdfstring{$\Gamma_{\text{sym}}$}{GammaSym}} 

We start by recalling a ladder-type functor which has appeared several times in various forms in the literature. The main idea goes back to \cite{CaKaMo-webs-skew-howe}.

\begin{Lemma}(\cite[Lemma 2.17]{RoTu-symmetric-howe}.) \label{lem-ladderfunc}
For each $m \in \mathbb{Z}_{\geq 0}$, there exists a functor 
\[A_{m,n} \colon \dot{\mathrm{U}}_{q}(\mathfrak{gl}_m) \to \mathbf{SWeb}_\uparrow(\mathfrak{gl}_n) \]
which sends a $\mathfrak{gl}_m$-weight $\vec{k} \in \mathbb{Z}_{\geq 0}^m$ to the sequence obtained by removing all $0$'s, and sends all other objects $\vec{k}$ of $\dot{\mathrm{U}}_{q}(\mathfrak{gl}_m) $ to the zero object. This functor is determined on morphisms by the assignment
\begin{gather*}
A_{m,n}(F_i^{(j)} \mathbf{1}_{\vec{k}}) = 
\xy
(0,0)*{
\begin{tikzpicture}[scale=.3]
\draw [very thick, directed=.7] (-4,-2) to (-4,2);
\draw [very thick, directed=.7] (-1.5,-2) to (-1.5,0);
\draw [very thick] (-1.5,0) to (-1.5,2);
\draw [very thick, directed=.7] (1.5,-2) to (1.5,0);
\draw [very thick] (1.5,0) to (1.5,2);
\draw [very thick, directed=.55] (-1.5,0) to (1.5,0);
\draw [very thick, directed=.7] (4,-2) to (4,2);
\node at (-4,-2.5) {\tiny $k_1$};
\node at (-1.5,-2.5) {\tiny $k_i$};
\node at (1.5,-2.5) {\tiny $k_{i+1}$};
\node at (4,-2.5) {\tiny $k_m$};
\node at (-1.65,2.5) {\tiny $k_i{-}j$};
\node at (1.65,2.5) {\tiny $k_{i+1}{+}j$};
\node at (-4,2.5) {\tiny $k_1$};
\node at (4,2.5) {\tiny $k_m$};
\node at (0,0.75) {\tiny $j$};
\node at (-2.75,0) {$\cdots$};
\node at (2.75,0) {$\cdots$};
\end{tikzpicture}
};
\endxy 
\quad , \quad
A_{m,n}(E_i^{(j)} \mathbf{1}_ {\vec{k}}) = 
\xy
(0,0)*{
\begin{tikzpicture}[scale=.3]
\draw [very thick, directed=.7] (-4,-2) to (-4,2);
\draw [very thick, directed=.7] (-1.5,-2) to (-1.5,0);
\draw [very thick] (-1.5,0) to (-1.5,2);
\draw [very thick, directed=.7] (1.5,-2) to (1.5,0);
\draw [very thick] (1.5,0) to (1.5,2);
\draw [very thick, directed=.55] (1.5,0) to (-1.5,0);
\draw [very thick, directed=.7] (4,-2) to (4,2);
\node at (-4,-2.5) {\tiny $k_1$};
\node at (-1.5,-2.5) {\tiny $k_i$};
\node at (1.5,-2.5) {\tiny $k_{i+1}$};
\node at (4,-2.5) {\tiny $k_m$};
\node at (-1.65,2.5) {\tiny $k_i{+}j$};
\node at (1.65,2.5) {\tiny $k_{i+1}{-}j$};
\node at (-4,2.5) {\tiny $k_1$};
\node at (4,2.5) {\tiny $k_m$};
\node at (0,0.75) {\tiny $j$};
\node at (-2.75,0) {$\cdots$};
\node at (2.75,0) {$\cdots$};
\end{tikzpicture}
};
\endxy 
\end{gather*}
where we erase any zero labelled edges in the diagrams depicting the images.\qed
\end{Lemma}

Let us give some necessary $\text{U}_{q}(\mathfrak{gl}_n)$-intertwiners, using a thin version of the notations in \cite[Appendix B]{RoWa-sym-homology}.

\begin{Definition}\label{quantum proj incl ev coev} 
The \textit{thin projection, inclusion, evaluation and coevaluation maps} are defined as
\begin{equation*}
m_{1,1} : \text{Sym}_q^1 \otimes \text{Sym}_q^1 \rightarrow \text{Sym}_q^{2} 
\end{equation*}
\begin{equation*} b_i \otimes b_j  \mapsto \begin{cases}
b_i \otimes b_j, &\quad\text{if}  \quad i \leq j,\\
qb_j \otimes b_i, &\quad\text{if} \quad i>j\changed{,} \\
\end{cases}
\end{equation*}
\begin{equation*}
s_{1,1} :  \text{Sym}_q^{2} \rightarrow \text{Sym}_q^1 \otimes \text{Sym}_q^1 
\end{equation*}
\begin{equation*} b_i \otimes b_j  \mapsto \begin{cases}
q^{-1} b_i \otimes b_j +b_j \otimes b_i, &\quad\text{if}  \quad i < j,\\
[2]b_i \otimes b_i, &\quad\text{if} \quad i=j\changed{,} \\
\end{cases}
\end{equation*}
\begin{equation*}
ev_1^\text{left}: (\text{Sym}_q^1)^{*} \otimes \text{Sym}_q^1 \rightarrow \Bbbk_q
\end{equation*}
\begin{equation*}
b_i^* \otimes b_j \mapsto b_i^*(b_j)\changed{,}
\end{equation*}
\begin{equation*}
coev_1^\text{left}: \Bbbk_q \rightarrow \text{Sym}_q^1 \otimes (\text{Sym}_q^1)^*
\end{equation*}
\begin{equation*}
1 \mapsto   \sum_{i=1}^{n} b_i \otimes b_i^*.
\end{equation*}
\end{Definition}

\begin{Lemma}\label{thin maps are intertwiners} 
The maps $m_{1,1}, \ s_{1,1}, \  ev_1^\text{left}, \ coev_1^\text{left} $ are $\text{U}_{q}(\mathfrak{gl}_n)$-intertwiners.
\end{Lemma}

Recall now the $q$-symmetric Howe functors from \autoref{Howe duality}. We use them to define our main functor.

\begin{Definition}\label{mainfunctor} 
The functor 
\[ \Gamma_{\text{sym}}: \mathbf{SWeb}_{\uparrow,\downarrow}(\mathfrak{gl}_n) \rightarrow \text{U}_{q}(\mathfrak{gl}_n)\text{-}\mathrm{\mathbf{fdMod}}_{S,S^*} 
\]
is defined as follows.
\begin{itemize}
\item On objects $k_{\uparrow}$ and  $k_{\downarrow}$, $k \in \mathbb{Z}_{\geq 0}$ we set
\[ \Gamma_{\text{sym}}(k_{\uparrow} )= \text{Sym}_q^{k}   \quad \text{and} \quad
\Gamma_{\text{sym}}(k_{\downarrow} )= (\text{Sym}_q^{k})^* .
\]
By convention, we send the zero object to the zero module and the empty tuple to the trivial $\text{U}_{q}(\mathfrak{gl}_n)$-module.
\item On morphisms, we send the generating morphisms of $\mathbf{SWeb}_{\uparrow,\downarrow}(\mathfrak{gl}_n)$ to the following $\text{U}_{q}(\mathfrak{gl}_n)$-intertwiners and we extend monoidally
\[
\Gamma_\text{sym} \left(
\xy
(0,0)*{
\begin{tikzpicture}[scale=.3]
\draw [very thick, directed=.6] (0, .75) to (0,2.5);
\draw [very thick, directed=.6] (1,-1) to [out=90,in=330] (0,.75);
\draw [very thick,directed=.6] (-1,-1) to [out=90,in=210] (0,.75); 
\node at (0, 3) {\tiny $k{+}1$};
\node at (-1,-1.5) {\tiny $k$};
\node at (1,-1.5) {\tiny $1$};
\end{tikzpicture}
};
\endxy
\right)
=
\Phi_{2,n}(E_1^{(1)}\mathbf{1}_{(k,1)}), \
\Gamma_\text{sym} \left(
\xy
(0,0)*{
\begin{tikzpicture}[scale=.3]
\draw [very thick,directed=.6] (0,-1) to (0,.75);
\draw [very thick,directed=.7] (0,.75) to [out=30,in=270] (1,2.5);
\draw [very thick,directed=.7] (0,.75) to [out=150,in=270] (-1,2.5); 
\node at (0, -1.5) {\tiny $k{+}1$};
\node at (-1,3) {\tiny $k$};
\node at (1,3) {\tiny $1$};
\end{tikzpicture}
};
\endxy
\right)
=
\Phi_{2,n}(F_1^{(1)}\mathbf{1}_{(k+1,0)}), 
\]
\[
\Gamma_\text{sym} \left(
\xy
(0,0)*{
\begin{tikzpicture}[scale=.3]
\draw [very thick, rdirected=.6] (0, .75) to (0,2.5);
\draw [very thick, rdirected=.6] (1,-1) to [out=90,in=330] (0,.75);
\draw [very thick,rdirected=.6] (-1,-1) to [out=90,in=210] (0,.75); 
\node at (0, 3) {\tiny $k{+}1$};
\node at (-1,-1.5) {\tiny $k$};
\node at (1,-1.5) {\tiny $1$};
\end{tikzpicture}
};
\endxy
\right)
=
\Phi_{2,n}^*(F_1^{(1)}\mathbf{1}_{(-1,-k)}), \
\Gamma_\text{sym} \left(
\xy
(0,0)*{
\begin{tikzpicture}[scale=.3]
\draw [very thick,rdirected=.6] (0,-1) to (0,.75);
\draw [very thick,rdirected=.7] (0,.75) to [out=30,in=270] (1,2.5);
\draw [very thick,rdirected=.7] (0,.75) to [out=150,in=270] (-1,2.5); 
\node at (0, -1.5) {\tiny $k{+}1$};
\node at (-1,3) {\tiny $k$};
\node at (1,3) {\tiny $1$};
\end{tikzpicture}
};
\endxy
\right)
=
\Phi_{2,n}^*(E_1^{(1)}\mathbf{1}_{(0,-1-k)}), 
\]
\[
\Gamma_\text{sym} \left(
\xy
(0,0)*{
\begin{tikzpicture} [scale=1]
\draw[very thick,rdirected=-.9] (0,0) to [out=90,in=180] (.25,.5) to [out=0,in=90] (.5,0);
\node at (.25,.75) {};
\node at (0,-.15) {\tiny $1$};
\node at (.5,-.15) {\tiny $1$};
\end{tikzpicture}
}
\endxy
\right)
= ev_1^{\text{left}}, \
\Gamma_\text{sym} \left(
\xy
(0,0)*{
\begin{tikzpicture} [scale=1]
\draw[very thick,rdirected=-.9] (0,1) to [out=270,in=180] (.25,.5) to [out=0,in=270] (.5,1);
\node at (.25,.25) {};
\node at (0,1.15) {\tiny $1$};
\node at (.5,1.15) {\tiny $1$};
\end{tikzpicture}
}
\endxy
\right)
= coev_1^{\text{left}}.
\]
\end{itemize}
\end{Definition}

We now show that the functor $\Gamma_\text{sym}$ is actually well-defined. Here it will become clear why we constructed the diagrammatic category the way we did.

\begin{Lemma}\label{welldef of mf}
The functor 
\[\Gamma_\text{sym} :\mathbf{SWeb}_{\uparrow,\downarrow}(\mathfrak{gl}_n) \rightarrow \mathrm{U}_{q}(\mathfrak{gl}_n)\text{-}\mathrm{\mathbf{fdMod}}_{S,S^*}\] is well-defined.
\end{Lemma}

\begin{proof}
Let us show that the analogues of the relations from \autoref{symm-webs} hold in $\text{U}_{q}(\mathfrak{gl}_n)$-$\mathrm{\mathbf{fdMod}}$. The relations \autoref{(co)assoc} and \autoref{thinsqure-sw} hold in $\text{U}_{q}(\mathfrak{gl}_n)$-$\mathrm{\mathbf{fdMod}}$ as a consequence of \changed{$q$-Howe duality}, namely of \autoref{howe-func} since they come from relations in $\dot{\text{U}}_{q}(\mathfrak{gl}_m)$.
An easy calculation identifies the morphisms associated to the right-hand sides in \autoref{thin-braid} as 
\[ -q^{-1- \frac{1}{n}} \text{id}_1 \otimes \text{id}_1 + q^{-\frac{1}{n}} s_{1,1} \circ m_{1,1} \ \text{and} \ -q^{-1+ \frac{1}{n}} \text{id}_1 \otimes \text{id}_1 + q^{\frac{1}{n}} s_{1,1} \circ m_{1,1},
\]
respectively. These morphisms are invertible, and composing them with evaluation and coevaluation maps does not change that property, meaning that \autoref{inv-br} also holds. Now, a simple computation shows that \autoref{eq-snake} holds in $\text{U}_{q}(\mathfrak{gl}_n)$-$\mathrm{\mathbf{fdMod}}$. Indeed, from \autoref{quantum proj incl ev coev}, we have
\[ (ev_1^\text{left} \circ \text{id}_1) \circ (\text{id}_1 \circ coev_1^\text{left}) (b_j^* \otimes 1)=(ev_1^\text{left}  \circ \text{id}_1) (b_j^* \otimes \sum_{i=1}^{n} b_i \otimes b_i^*)= \sum_{i=1}^n b_j^*(b_i)\otimes b_i^* =1 \otimes b_j^*,\]
meaning that $(ev_1^\text{left}  \circ \text{id}_1) \circ (\text{id}_1 \circ coev_1^\text{left} )=\text{id}_1$, which is exactly the representation theoretical version of the zigzag relation from \autoref{eq-snake}. Showing relation \autoref{eq-vertexslide} is again a brute force computation which we will omit here, see \cite{La-gln-webs} for details.
Finally, relation \autoref{circ-ev} holds because of our definition of $ev_1^\text{left}$:
\[ ev_1^\text{left} \circ coev_1^\text{right} (1)=ev_1^\text{left}(\sum_{i=1}^n q^{-n+2i-1}b_i^* \otimes b_i)=q^{-n+1}+q^{-n+3}+\dots+q^{n-1}=[n].\]
The proof completes.
\end{proof} 

\subsection{Relating the two web categories}

Now we want to identify $\mathbf{SWeb}_{\uparrow}(\mathfrak{gl}_n) $ as a full subcategory of $\mathbf{SWeb}_{\uparrow,\downarrow}(\mathfrak{gl}_n) $. Recall \autoref{dualising webs}. In particular, from this proposition we can assume that the boundaries of a web in $\mathbf{SWeb}_{\uparrow, \downarrow}(\mathfrak{gl}_n)$ are upward pointing. Below we show that we have an even stronger connection between our two symmetric web categories.

\begin{Theorem}\label{up-vs-down}
There exists a fully faithful functor of braided monoidal categories
\[ \mathrm{G}: \mathbf{SWeb}_{\uparrow}(\mathfrak{gl}_n) \rightarrow \mathbf{SWeb}_{\uparrow, \downarrow}(\mathfrak{gl}_n)
\]
sending objects and morphism from $\mathbf{SWeb}_{\uparrow}(\mathfrak{gl}_n)$ to objects and morphisms of the same name in $\mathbf{SWeb}_{\uparrow, \downarrow}(\mathfrak{gl}_n)$.
\end{Theorem}

\begin{proof}
We give a line-to-line adaptation of \cite[Theorem 7.5.1]{BrDaKu-qwebs-typep}. First, it is clear that this functor is well-defined because the defining relations of $\mathbf{SWeb}_{\uparrow}(\mathfrak{gl}_n) $ hold in $\mathbf{SWeb}_{\uparrow,\downarrow}(\mathfrak{gl}_n) $. Also, they have the same monoidal braided structures by construction. \ochanged{It remains} to prove fully faithfulness. Let $k_{\uparrow}, \ l_{\uparrow}\in \text{Obj}(\mathbf{SWeb}_{\uparrow}(\mathfrak{gl}_n)) $  and $w$ a web in $\text{Hom}_{\mathbf{SWeb}_{\uparrow,\downarrow}(\mathfrak{gl}_n)}(k_{\uparrow},l_{\uparrow})$. 
We want to show that 
\[ \mathrm{Hom}_{\mathbf{SWeb}_{\uparrow}(\mathfrak{gl}_n)}(k_\uparrow, l_\uparrow) \cong \mathrm{Hom}_{\mathbf{SWeb}_{\uparrow,\downarrow}(\mathfrak{gl}_n)} (\mathrm{G}(k_\uparrow), \mathrm{G}(l_\uparrow)).
\]
Let us first show fullness. Note that the bottom and top of a web in $\mathbf{SWeb}_{\uparrow, \downarrow}(\mathfrak{gl}_n)$ are upward pointing but in the middle of it we can have both upward and downward pointing merges, splits, cups, caps and crossings. We use relations \autoref{eq-snake} and \autoref{circ-ev} to straighten zigzags and to remove circles. Further, we use \autoref{downard thick merge and split} to remove any downward pointing merges and splits. From \autoref{thick cu(a)p sliding} and \autoref{pitchfork} we can slide crossing, cups and caps through merges and splits putting the former below the \ochanged{latter}. Proceeding inductively on the total number of crossings, we can therefore write $w$ as a linear combination of diagrams of the form $k_{\uparrow}\rightarrow l_{\uparrow}$, where each diagram has an upper part and a lower part. The upper part consists only of upward pointing merges and splits, whereas the lower part consists only of cups, caps and crossings (both upward and downward pointing). The webs of the upper part are clearly morphisms in $\mathbf{SWeb}_{\uparrow}(\mathfrak{gl}_n) $. In the lower part, after removing any circles, the remaining webs are (compositions of) crossings or compositions as in the Reidemeister 1 moves. Note that in the bottom part we have a geometric braid, which can be turned upwards. Now we use the three types of Reidemeister moves, in particular, we use \autoref{thick Reidemeister} to write the compositions of crossings with cups and caps as scalar multiples of identities. In a nutshell, every web of form $k_{\uparrow}\rightarrow l_{\uparrow}$ in $\mathbf{SWeb}_{\uparrow,\downarrow}(\mathfrak{gl}_n) $ can be written as a linear combination of diagrams of the form $k_{\uparrow} \rightarrow l_{\uparrow}$ in $\mathbf{SWeb}_{\uparrow}(\mathfrak{gl}_n) $. This shows fullness of $\mathrm{G}$. Finally, consider the diagram
\[
\begin{gathered}
\xymatrix{
\mathbf{SWeb}_\uparrow(\mathfrak{gl}_n) \ar[rr]^{\Gamma_\text{sym}} \ar[dr]_{G} & & \text{U}_q(\mathfrak{gl}_n) \text{-}\mathbf{fdMod} \\
& \hspace*{0.2cm} \mathbf{SWeb}_{\uparrow, \downarrow} (\mathfrak{gl}_n) \ar[ur]_{\Gamma_\text{sym}} &
}
\end{gathered} .
\]
It commutes by definition. Moreover, the top functor is faithful from \autoref{Gamma sym is fully faithful}. This then implies that $G$ is faithful.
\end{proof}

\subsection{Proof that \texorpdfstring{$\Gamma_{\text{sym}}$}{GammaSym} is an equivalence of categories} \label{proof of equiv (sym generic)} 

We extend our functor by passing to the additive Karoubi envelopes of our categories. We abuse notation and continue denoting it by $\Gamma_{\text{sym}}$
\[
\Gamma_{\text{sym}}: \textbf{Kar}(\mathbf{SWeb}_{\uparrow,\downarrow}(\mathfrak{gl}_n) )\rightarrow \text{U}_{q}(\mathfrak{gl}_n)\text{-}\mathrm{\mathbf{fdMod}}.\]
We are finally ready to prove the following. 

\begin{Theorem}\label{main-thm} 
\textit{(Equivalence Theorem for symmetric webs)} The functor 
\[
\Gamma_{\mathrm{sym}}: \mathbf{SWeb}_{\uparrow,\downarrow}(\mathfrak{gl}_n) \rightarrow \mathrm{U}_{q}(\mathfrak{gl}_n)\text{-}\mathbf{fdMod}_{S,S^*}\]
is fully faithful and
\[
\Gamma_{\mathrm{sym}}: \mathbf{Kar}(\mathbf{SWeb}_{\uparrow,\downarrow}(\mathfrak{gl}_n) )\rightarrow \mathrm{U}_{q}(\mathfrak{gl}_n)\text{-}\mathrm{\mathbf{fdMod}}\]
is an equivalence of braided pivotal categories.
\end{Theorem}

\begin{proof} We have already proved that $\Gamma_{\mathrm{sym}}$ is well-defined in \autoref{welldef of mf}. Now we show that it is fully faithful and essentially surjective after passing to Karoubi envelopes. \\
\underline{Fully faithful.} We have to show that for $\vec{k}, \vec{l} \in \mathbb{Z}^m$, $m>0$, we have
\[ \mathrm{Hom}_{\mathbf{SWeb}_{\uparrow,\downarrow}(\mathfrak{gl}_n)}(\vec{k}, \vec{l}) \cong \mathrm{Hom}_{\mathrm{U}_{q}(\mathfrak{gl}_n)\text{-}\mathrm{\mathbf{fdMod}}}(\Gamma_\text{sym}(\vec{k}),\Gamma_\text{sym}(\vec{l})).
\]
From \autoref{dualising webs} and \autoref{up-vs-down}, it suffices to show this for upward pointing webs only and we have proved the \ochanged{latter} in \autoref{Gamma sym is fully faithful}. \\
\underline{Essentially surjective.} This follows from the definition of $\Gamma_\text{sym}$ and the fact that every irreducible $\mathrm{U}_{q}(\mathfrak{gl}_n)$-module appears in a suitable tensor product of symmetric powers of the standard representation and their duals. Indeed, let $V$ be an irreducible $\mathrm{U}_{q}(\mathfrak{gl}_n)$-module appearing as a direct summand in $\bigotimes_i \text{Sym}_q^{k_i} \otimes \bigotimes_j (\text{Sym}_q^{k_j})^*.$ Then $\Gamma_{\mathrm{sym}}(\bigotimes_i {k_i}_{\uparrow} \otimes \bigotimes_j {k_j}_{\downarrow})=\bigotimes_i \text{Sym}_q^{k_i} \otimes \bigotimes_j (\text{Sym}_q^{k_j})^*$. Since both sides are now idempotent complete, tensor products decompose into irreducible objects and such decompositions are unique up to permutation. Thus, there exists a direct  summand of $\bigotimes_i {k_i}_{\uparrow} \otimes \bigotimes_j {k_j}_{\downarrow}$ that is sent to $V$ by $\Gamma_{\mathrm{sym}}$. \\
This finishes the proof that $\Gamma_\text{sym}$ is an equivalence of categories. Finally, after adapting Section 3.1 of \cite{RoTu-symmetric-howe} and Section 6 of \cite{CaKaMo-webs-skew-howe} to our setting we have that $\Gamma_\text{sym}$ is an equivalence of pivotal braided categories.
\end{proof}

Now we give the connection between the three functors of this section. The following is a generalisation of Corollary 2.22 in \cite{RoTu-symmetric-howe}.

\begin{Lemma}
The diagram below commutes \label{comm of dot-mod}
\begin{equation*}
\begin{gathered}
\xymatrix{
\dot{\mathrm{U}}_{q}(\mathfrak{gl}_m) \ar[rr]^{\Phi_{m,n}} \ar[dr]_{A_{m,n}} & & \mathrm{U}_{q}(\mathfrak{gl}_n)\text{-} \mathrm{\mathbf{fdMod}} \\
& \hspace*{0.2cm} \mathbf{SWeb}_\uparrow(\mathfrak{gl}_n) \ar[ur]_{\Gamma_\text{sym}} &
}
\end{gathered}.
\end{equation*}
\end{Lemma}

The following lemma is one of the main ingredients for proving \autoref{main-thm}. Note that it concerns upward pointing webs only. 

\begin{Lemma}\label{Gamma sym is fully faithful} 
The functor \[
\Gamma_{\text{sym}}: \mathbf{SWeb}_{\uparrow}(\mathfrak{gl}_n)\rightarrow \mathrm{U}_{q}(\mathfrak{gl}_n)\text{-}\mathrm{\mathbf{fdMod}}_{S}\]
is fully faithful. 
\end{Lemma} 

\begin{proof}
We have to show that for tuples $k_\uparrow=(k_1, \dots, k_m)_\uparrow \in \mathbb{Z}_{\geq 0}^m$ and $ l_\uparrow= (l_1, \dots, l_{m'})_\uparrow$ $\in \mathbb{Z}_{\geq 0}^{m'}$ we have
\[ \mathrm{Hom}_{\mathbf{SWeb}_{\uparrow}(\mathfrak{gl}_n)}(k_\uparrow, l_\uparrow) \cong \mathrm{Hom}_{\mathrm{U}_{q}(\mathfrak{gl}_n)\text{-}\mathrm{\mathbf{fdMod}}}(\Gamma_\text{sym}(k_\uparrow),\Gamma_\text{sym}(l_\uparrow)).
\]
Note that we can assume that $m=m'$, otherwise we add zeros to the shorter tuple. Now, $\Gamma_\text{sym}$ by definition comes from a $q$-Howe functor $\Phi_{m,n}$, for $m \in \mathbb{Z}_{\geq 0}$ big enough.\\
\textit{Surjectivity:} In the case $\sum \limits_{i=1}^mk_i=\sum \limits_{i=1}^{m}l_i$ surjectivity follows from \autoref{comm of dot-mod} and the fact that $q$-symmetric Howe functor $\Phi_{m,n}$ is surjective on these Hom-spaces, whereas for $\sum \limits_{i=1}^mk_i \neq \sum \limits_{i=1}^{m}l_i$ we have zero on both sides.\\
\textit{Injectivity:} By surjectivity and finite-dimensionality of the involved spaces it suffices to show that
\begin{equation} \label{hom spaces thin}  \mathrm{dimHom}_{\mathbf{SWeb}_{\uparrow}(\mathfrak{gl}_n)}(k_\uparrow,l_\uparrow) =\mathrm{dimHom}_{\mathrm{U}_{q}(\mathfrak{gl}_n)\text{-}\mathrm{\mathbf{fdMod}}}(\Gamma_\text{sym}(k_\uparrow),\Gamma_\text{sym}(l_\uparrow))
\xy
(0,0)*{
\begin{tikzpicture}[scale=.3]
\end{tikzpicture}
};
\endxy.
\end{equation}
We prove this first for $k_\uparrow=(1,\dots,1)_\uparrow=l_\uparrow$. In this case, recall that in $\mathrm{U}_{q}(\mathfrak{gl}_n)\text{-}\mathrm{\mathbf{fdMod}}$ we have
\begin{center}
$
\text{Sym}_q^{1} \Bbbk_q^n \cong \bigwedge_{q}^{1}\Bbbk_q^n \cong \Bbbk_q^n. 
$
\end{center}
Further, we have 
\begin{equation} \label{eq:some} \text{dim}\mathrm{End}_{\mathrm{U}_{q}(\mathfrak{gl}_n)\text{-}\mathrm{\mathbf{fdMod}}}(\Bbbk_q^n \otimes \dots \otimes \Bbbk_q^n )=\text{dim}\mathrm{End}_{\mathrm{U}_{q}(\mathfrak{sl}_n)\text{-}\mathrm{\mathbf{fdMod}}}(\Bbbk_q^n \otimes \dots \otimes \Bbbk_q^n ), 
\end{equation}
since every irreducible $\mathrm{U}_{q}(\mathfrak{sl}_n)$-module comes from a restriction of an irreducible $\mathrm{U}_{q}(\mathfrak{gl}_n)$-module.

Now, from \autoref{thin basis} we have that
\begin{equation}\label{eq:some2}
\mathrm{dimEnd}_{\mathbf{SWeb}_\uparrow(\mathfrak{gl}_n)}(1 \otimes \dots \otimes1) = \mathrm{dimEnd} _{\mathbf{EWeb}_{\uparrow,\downarrow}(\mathfrak{sl}n)}(1 \otimes \dots\otimes 1). 
\end{equation}
Finally, there exists an equivalence of categories (\cite{CaKaMo-webs-skew-howe}, Theorem 3.3.1) 
\[\mathbf{EWeb}_{\uparrow,\downarrow}(\mathfrak{sl}_n) \rightarrow \mathrm{U}_{q}(\mathfrak{sl}_n)\text{-}\mathrm{\mathbf{fdMod}}.\]
Combining this equivalence with \autoref{eq:some} and \autoref{eq:some2} shows that the equality in \autoref{hom spaces thin} holds for $k_\uparrow=(1,\dots,1)_\uparrow=l_\uparrow$. For the general case, we proceed as in the proof of Theorem 1.10 in \cite{RoTu-symmetric-howe}. Let $k_\uparrow, l_\uparrow$ be such that $\sum \limits_{i=1}^mk_i=\sum \limits_{i=1}^{m}l_i$ and let $u,v \in \text{Hom}_{\mathbf{SWeb}_\uparrow}(k_\uparrow, l_\uparrow)$, $u \neq v$. Composing $u$ and $v$ with merges and splits we get 
\[ u'=
\xy
(0,0)*{
\begin{tikzpicture} [scale=1]
\draw[very thick] (0,1.25) rectangle (2,.75);
\draw[very thick] (.25,0) to (.25,.75);
\draw[very thick] (.25,1.25) to (.25,2);
\draw[very thick] (1.75,0) to (1.75,.75);
\draw[very thick] (1.75,1.25) to (1.75,2);
\node at (1,1.75) {$\cdots$};
\node at (1,.25) {$\cdots$};
\node at (1,1) {$u$};
\draw [very thick] (0.5,-.5) to [out=90,in=330] (0.25,0);
\draw [very thick] (0,-.5) to [out=90,in=210] (0.25,0);
\draw [very thick] (2,-.5) to [out=90,in=330] (1.75,0);
\draw [very thick] (1.5,-.5) to [out=90,in=210] (1.75,0);
\draw [very thick] (0.25,2) to [out=30,in=270] (0.5,2.5);
\draw [very thick] (0.25,2) to [out=150,in=270] (0,2.5);
\draw [very thick] (1.75,2) to [out=30,in=270] (2,2.5);
\draw [very thick] (1.75,2) to [out=150,in=270] (1.5,2.5);
\node at (0.25,2.65) {\tiny $\cdots$};
\node at (0.25,-.65) {\tiny $\cdots$};
\node at (1.75,2.65) {\tiny $\cdots$};
\node at (1.75,-.65) {\tiny $\cdots$};
\node at (0,2.65) {\tiny $1$};
\node at (0.5,2.65) {\tiny $1$};
\node at (0,-.65) {\tiny $1$};
\node at (0.5,-.65) {\tiny $1$};
\node at (1.5,2.65) {\tiny $1$};
\node at (2,2.65) {\tiny $1$};
\node at (1.5,-.65) {\tiny $1$};
\node at (2,-.65) {\tiny $1$};
\end{tikzpicture}
};
\endxy, \quad 
v'=
\xy
(0,0)*{
\begin{tikzpicture} [scale=1]
\draw[very thick] (0,1.25) rectangle (2,.75);
\draw[very thick] (.25,0) to (.25,.75);
\draw[very thick] (.25,1.25) to (.25,2);
\draw[very thick] (1.75,0) to (1.75,.75);
\draw[very thick] (1.75,1.25) to (1.75,2);
\node at (1,1.75) {$\cdots$};
\node at (1,.25) {$\cdots$};
\node at (1,1) {$v$};
\draw [very thick] (0.5,-.5) to [out=90,in=330] (0.25,0);
\draw [very thick] (0,-.5) to [out=90,in=210] (0.25,0);
\draw [very thick] (2,-.5) to [out=90,in=330] (1.75,0);
\draw [very thick] (1.5,-.5) to [out=90,in=210] (1.75,0);
\draw [very thick] (0.25,2) to [out=30,in=270] (0.5,2.5);
\draw [very thick] (0.25,2) to [out=150,in=270] (0,2.5);
\draw [very thick] (1.75,2) to [out=30,in=270] (2,2.5);
\draw [very thick] (1.75,2) to [out=150,in=270] (1.5,2.5);
\node at (0.25,2.65) {\tiny $\cdots$};
\node at (0.25,-.65) {\tiny $\cdots$};
\node at (1.75,2.65) {\tiny $\cdots$};
\node at (1.75,-.65) {\tiny $\cdots$};
\node at (0,2.65) {\tiny $1$};
\node at (0.5,2.65) {\tiny $1$};
\node at (0,-.65) {\tiny $1$};
\node at (0.5,-.65) {\tiny $1$};
\node at (1.5,2.65) {\tiny $1$};
\node at (2,2.65) {\tiny $1$};
\node at (1.5,-.65) {\tiny $1$};
\node at (2,-.65) {\tiny $1$};
\end{tikzpicture}
};
\endxy.
\]
Now, recall that the digon removal is invertible since quantum numbers are invertible (\autoref{invertibility of digon}) implying that $u' \neq v'$. Finally, from the above argument we have that $\Gamma_\text{sym}(u') \neq \Gamma_\text{sym}(v')$, showing faithfulness for the general case. 
\end{proof}

\section{Minimality Theorem}\label{minimality} 

In this section we show that we have used a ``minimal'' set of web generators and relations to prove the Equivalence Theorem. 

It is worth emphasizing that in contrast to e.g. \cite{CaKaMo-webs-skew-howe}, in our setting, the generators of the web category are as thin as possible and we use a shorter list of defining relations to prove the equivalence between our web category and the category of finite-dimensional $ \text{U}_{q}(\mathfrak{gl}_n)$-modules. The whole point of this section will be to make the notion of minimality precise. Roughly, if we agreed on a set of generators for our web categories, like merges and splits, then the ones we have chosen are as thin as possible, while the relations they satisfy cannot be derived from one another. \\

Throughout this section, all the notions and concepts apply to symmetric, exterior and green-red webs \changed{for $n\gg0$ (so that we do not need to worry about antisymmetrizer relations)}. We will write $\mathrm{Hom}_{\textbf{Web}} $ for the Hom-spaces of any of these web categories. Moreover, note that a web is effectively an equivalence class, up to the web relations in the corresponding web category. 

\begin{Definition} 
Let $w$ be a web. Its \textit{complexity} is the double minimum, in lexicographical order and running over all webs in its equivalence class, of its biggest label $k$ and its number of appearances $\#k$ of $k$. We denote it by $(k, \#k)$.

\changed{We also apply the notion of complexity to linear combinations of webs. By definition, the complexity of a linear combination is the maximal complexity of appearing webs with non-zero coefficient.}
\end{Definition}

We define a pairing 
\begin{equation}
\langle-,- \rangle: \mathrm{Hom}_{\textbf{Web}} (\emptyset, \vec{k}) \times \mathrm{Hom}_{\textbf{Web}} (\vec{k}, \emptyset) \rightarrow \mathrm{End}_{\textbf{Web}}(\emptyset)
\end{equation}
by taking a web $w \in \mathrm{Hom}_{\textbf{Web}} (\emptyset, \vec{k}) $ and a web $w' \in \mathrm{Hom}_{\textbf{Web}} (\vec{k}, \emptyset) $ and glueing them along their common boundary. We often say that $w'$ is a \textit{closure} of the web $w$ and \ochanged{vice versa}.

For a general web $w \in \mathrm{Hom}_{\textbf{Web}} (\vec{k}, \vec{l}) $, we first bend it as explained in \autoref{dualising webs} (note that there is a choice of bending involved, and we take the one indicated in the proof of \autoref{dualising webs}). This way we obtain a web $\overline{w} \in \mathrm{Hom}_{\textbf{Web}} (\emptyset, \vec{m})$ and we can pair it as above.

\begin{Example} 
For the web
\[  
w=\begin{tikzpicture}[anchorbase,scale=0.7]
\draw[very thick,crossline,directed=1] (0,0) to (1,1);
\draw[very thick,crossline,directed=1] (0,1) to (1,0);
\node at (0,-.3) {\tiny$1$};
\node at (1,-.3) {\tiny$1$};
\node at (0,1.3) {\tiny$1$};
\node at (1,1.3) {\tiny$1$};
\end{tikzpicture}
\]
we have a closed web
\[
\begin{tikzpicture}[anchorbase,scale=0.7]
\draw[very thick,crossline,directed=1] (0,0) to (1,1);
\draw[very thick,crossline,directed=1] (0,1) to (1,0);
\draw[very thick,] (0,0) to[out=225,in=270] (-1,1);
\draw[very thick,] (1,0) to[out=315,in=270] (2,1);
\draw[very thick,blue] (-1,1) to[out=90,in=180] (-.5,2);
\draw[very thick,blue] (-.5,2) to[out=0,in=90] (0,1);
\draw[very thick,blue] (1,1) to[out=90,in=180] (1.5,2);
\draw[very thick,blue] (2,1) to[out=90,in=0] (1.5,2);
\draw[very thick,blue] (-.5,2) to[out=90,in=180] (.5,3);
\draw[very thick,blue] (.5,3) to[in=90,out=0] (1.5,2);
\node at (-1.2,1) {\tiny$1$};
\node at (2.2,1) {\tiny$1$};
\node at (.5,3.2) {\tiny$2$};
\node at (0.2,1) {\tiny$1$};
\node at (.6,1) {\tiny$1$};
\end{tikzpicture},
\]
which is obtained by first bending $w$, and then composing it with a half theta web (blue).
\end{Example}

The reason why we are particularly interested in closed webs becomes evident from the lemma below.

\begin{Lemma} 
We have $\mathrm{End}(\emptyset) \cong \Bbbk(q^{\frac{1}{n}}).$
\end{Lemma}

\begin{proof} 
Note that the statement is equivalent to the claim that to any closed web we can assign a scalar, which is the result of the evaluation rules in the corresponding web category. From \autoref{main-thm}, it suffices to prove the statement in the representation theoretical side. The claim then follows by a version of Schur's lemma. Precisely, although our ground field might not be algebraically closed, the endomorphism space of $\text{Sym}_q^0= \Bbbk_q$ is trivial by \cite[Corollary 7.4]{AnPoWe-representation-qalgebras}.
\end{proof}

\begin{Definition} 
We refer to the scalar  $\langle w,w' \rangle$ as the \textit{evaluation} of the web $w$ with respect to $w'$.
\end{Definition}

\begin{Definition} 
Let $w_1$ and $w_2$ be webs. Two web relations $w_1=r_1(w_1), \ w_2=r_2(w_2)$ are called \textit{dependent} if $w_1$ and $w_2$ have the same boundary up to bending and $\langle w_1, w' \rangle= \langle w_2, w' \rangle$, for every $w'$. Otherwise we call them \textit{independent}. (Here the webs $w_1$ and $w_2$ are allowed to be linear combinations of webs.)
\end{Definition}

\begin{Remark}
\changed{Our notion of dependence is phrased in terms of closures because, by the universal construction (see e.g. \cite{BlHaMaVo-tqft-kauffman-bracket}; standard in TQFT-like settings), pairings with all closures detect morphisms. Thus, if two candidate relations have the same evaluations against all closures, then they are indistinguishable from the point of view of the quotient category. In particular, this is the criterion relevant for the independence argument below, where we separate relations by exhibiting closures with different evaluations.}
\end{Remark}

We are now ready to state and prove the Minimality Theorem for symmetric webs, which can be proved similarly for exterior and green-red webs.

\begin{Theorem}\label{Minimality for symmetric webs} 
(Minimality Theorem for symmetric webs) We have the following:
\begin{itemize}
\item[1.] The relations \autoref{(co)assoc}-\autoref{circ-ev} are sufficient to prove \autoref{main-thm}.
\item[2.] The relations \autoref{(co)assoc}-\autoref{circ-ev} are independent.
\item[3.] The generators of $\mathbf{SWeb}_{\uparrow,\downarrow}(\mathfrak{gl}_n)$ cannot be made thinner while still having a fully faithful functor to $\mathrm{U}_q(\mathfrak{gl}_n)\text{-}\mathbf{fdMod}_{S,S^*}$.
\end{itemize}
\end{Theorem}

\begin{proof} \leavevmode
\begin{itemize}
\item[1.] Throughout we have used only relations \autoref{(co)assoc}-\autoref{circ-ev} to show that $\Gamma_{\text{sym}}$ is an equivalence of categories.
\item[2.] We list the boundary points of each of the relations. The relations \autoref{(co)assoc}-\autoref{inv-br} have all four boundary points, \autoref{eq-snake} has two, \autoref{eq-vertexslide} has three and \autoref{circ-ev} has zero. This means that the relations \autoref{eq-snake}, \autoref{eq-vertexslide} and \autoref{circ-ev} are incomparable to each other and to other relations.  
To show that \autoref{(co)assoc} and \autoref{inv-br} are independent we bend the webs such that all the boundary points are at the top (see also \autoref{dualising webs}). Then these two relations have different boundary orientations, so they are independent. Similarly, relations \autoref{(co)assoc} and \autoref{thinsqure-sw} are independent. Finally, we need to compare \autoref{thinsqure-sw} and \autoref{inv-br}, however only for $k=l=1$. The bend versions are 
\[
\begin{tikzpicture}[anchorbase,scale=0.5]
\draw[very thick,rdirected=1] (0,0) to (0.25,0.5)
to (0.75,0.5) to (1,1)node[above]{\tiny$1$};
\draw[very thick,rdirected=1] (0,1)node[above]{\tiny$1$} to (0.25,0.5)
to (0.75,0.5) to (1,0);
\draw[very thick] (0,0) to[out=225,in=270] (-1,1)node[above]{\tiny$1$};
\draw[very thick] (1,0) to[out=315,in=270] (2,1)node[above]{\tiny$1$};
\end{tikzpicture} 
\quad \text{and}
\begin{tikzpicture}[anchorbase,scale=0.5]
\draw[very thick,crossline,rdirected=0] (0,0) to (1,1)node[above]{\tiny$1$};
\draw[very thick,crossline,rdirected=0.1] (0,1)node[above]{\tiny$1$} to (1,0);
\draw[very thick] (0,0) to[out=225,in=270] (-1,1)node[above]{\tiny$1$};
\draw[very thick] (1,0) to[out=315,in=270] (2,1)node[above]{\tiny$1$};
\end{tikzpicture}.
\]
We close these webs (using the same closure) as below, 
\[
\begin{tikzpicture}[anchorbase,scale=0.7]
\draw[very thick,rdirected=1] (0,0) to (0.25,0.5)
to (0.75,0.5) to (1,1);
\draw[very thick,rdirected=1] (0,1) to (0.25,0.5)
to (0.75,0.5) to (1,0);
\draw[very thick] (0,0) to[out=225,in=270] (-1,1);
\draw[very thick] (1,0) to[out=315,in=270] (2,1);
\draw[very thick] (-1,1) to[out=90,in=180] (-.5,2);
\draw[very thick] (-.5,2) to[out=0,in=90] (0,1);
\draw[very thick] (1,1) to[out=90,in=180] (1.5,2);
\draw[very thick] (2,1) to[out=90,in=0] (1.5,2);
\draw[very thick] (-.5,2) to[out=90,in=180] (.5,3);
\draw[very thick] (.5,3) to[in=90,out=0] (1.5,2);
\node at (-1.2,1) {\tiny$1$};
\node at (2.2,1) {\tiny$1$};
\node at (.5,3.2) {\tiny$2$};
\node at (0.2,1) {\tiny$1$};
\node at (.7,1) {\tiny$1$};
\end{tikzpicture} 
\quad\text{and}\quad
\begin{tikzpicture}[anchorbase,scale=0.7]
\draw[very thick,crossline,rdirected=0] (0,0) to (1,1);
\draw[very thick,crossline,rdirected=0] (0,1) to (1,0);
\draw[very thick] (0,0) to[out=225,in=270] (-1,1);
\draw[very thick] (1,0) to[out=315,in=270] (2,1);
\draw[very thick] (-1,1) to[out=90,in=180] (-.5,2);
\draw[very thick] (-.5,2) to[out=0,in=90] (0,1);
\draw[very thick] (1,1) to[out=90,in=180] (1.5,2);
\draw[very thick] (2,1) to[out=90,in=0] (1.5,2);
\draw[very thick] (-.5,2) to[out=90,in=180] (.5,3);
\draw[very thick] (.5,3) to[in=90,out=0] (1.5,2);
\node at (-1.2,1) {\tiny$1$};
\node at (2.2,1) {\tiny$1$};
\node at (.5,3.2) {\tiny$2$};
\node at (0.2,1) {\tiny$1$};
\node at (.7,1) {\tiny$1$};
\end{tikzpicture}.
\]

After evaluating, using the digon removal and the circle evaluation, we get $[2]^2\begin{bmatrix}
n+1\\
2
\end{bmatrix}$ for the first web. On the other hand, using \autoref{splits&merges comp crossings}, digon removal and the circle evaluation we get $q^{-\frac{1}{n} +1}[2]\begin{bmatrix}
n+1\\
2
\end{bmatrix}$ for the second web. These are different scalars meaning that the relations are independent.
\item[3.] We show the statement by induction on the complexity. Recall that the ``thinnest'' generators on $\mathbf{SWeb}_{\uparrow,\downarrow}(\mathfrak{gl}_n)$ are
\[
\xy
(0,0)*{
\begin{tikzpicture}[scale=.3]
\draw [very thick, directed=.55] (0, .75) to (0,2.5);
\draw [very thick, directed=.45] (1,-1) to [out=90,in=330] (0,.75);
\draw [very thick, directed=.45] (-1,-1) to [out=90,in=210] (0,.75); 
\node at (0, 3) {\tiny $2$};
\node at (-1,-1.5) {\tiny $1$};
\node at (1,-1.5) {\tiny $1$};
\end{tikzpicture}  
};
\endxy, \quad 
\xy
(0,0)*{
\begin{tikzpicture}[scale=.3]
\draw [very thick, directed=.55] (0,-1) to (0,.75);
\draw [very thick, directed=.65] (0,.75) to [out=30,in=270] (1,2.5);
\draw [very thick, directed=.65] (0,.75) to [out=150,in=270] (-1,2.5); 
\node at (0, -1.5) {\tiny $2$};
\node at (-1,3) {\tiny $1$};
\node at (1,3) {\tiny $1$};
\end{tikzpicture}
};
\endxy , \quad
\xy
(0,0)*{
\begin{tikzpicture}[scale=.3]
\draw [very thick, rdirected=.55] (0, .75) to (0,2.5);
\draw [very thick, rdirected=.45] (1,-1) to [out=90,in=330] (0,.75);
\draw [very thick, rdirected=.45] (-1,-1) to [out=90,in=210] (0,.75); 
\node at (0, 3) {\tiny $2$};
\node at (-1,-1.5) {\tiny $1$};
\node at (1,-1.5) {\tiny $1$};
\end{tikzpicture}
};
\endxy, \quad
\xy
(0,0)*{
\begin{tikzpicture}[scale=.3]
\draw [very thick, rdirected=.55] (0,-1) to (0,.75);
\draw [very thick, rdirected=.65] (0,.75) to [out=30,in=270] (1,2.5);
\draw [very thick, rdirected=.65] (0,.75) to [out=150,in=270] (-1,2.5); 
\node at (0, -1.5) {\tiny $2$};
\node at (-1,3) {\tiny $1$};
\node at (1,3) {\tiny $1$};
\end{tikzpicture}
};
\endxy, \quad
\xy
(0,0)*{
\begin{tikzpicture} [scale=1]
\draw[very thick,rdirected=-.95] (0,1) to [out=270,in=180] (.25,.5) to [out=0,in=270] (.5,1);
\node at (.25,.25) {};
\node at (0,1.25) {\tiny $1$};
\node at (.5,1.25) {\tiny $1$};
\end{tikzpicture}
}
\endxy ,\quad 
\xy
(0,0)*{
\begin{tikzpicture} [scale=1]
\draw[very thick,rdirected=-.95] (0,0) to [out=90,in=180] (.25,.5) to [out=0,in=90] (.5,0);
\node at (.25,.75) {};
\node at (0,-.25) {\tiny $1$};
\node at (.5,-.25) {\tiny $1$};
\end{tikzpicture}
}
\endxy . 
\]
\end{itemize}
Note that the corresponding Hom-spaces in $\mathrm{U}_q(\mathfrak{gl}_n)\text{-}\mathbf{fdMod}_{S,S^*}$ are the first non-trivial ones.

Let us inductively assume that we have already shown that the $(k, 1)$-merge and $(k, 1)$-split are needed as generators. On the representation theoretical side $(k + 1, 1)$ corresponds to $\text{Sym}^{k+1} \otimes \Bbbk_q^n$. The \ochanged{latter} has the summand $\text{Sym}^{k+2}$ which we have not seen in any shorter tensor product. In particular, none of the inductively added maps can go through this summand, so we need to add at least two additional maps (corresponding to the splitting of $\text{Sym}^{k+2}$ from $\text{Sym}^{k+1} \otimes \Bbbk_q^n$) in order for the Hom-space to be of the right dimension. We add the $(k + 1, 1)$-merge and the $(k + 1, 1)$-split. The same arguments work in the dual setting.
\end{proof}

\section{Some other properties of webs}

We now point out that webs can be adjusted to work integrally, meaning over $\mathbb{Z}[q^\frac{1}{n}, q^{-\frac{1}{n}}]$. By specialisation, this could lead to understanding something about the representation theory in positive characteristic.

\begin{Definition}\label{free int up-webs} 
The \textit{free symmetric upward pointing category of integral} $\mathfrak{gl}_n$-webs $\mathbf{fSWeb}_{\uparrow}^{\text{int}}(\mathfrak{gl}_n)$ is the $\mathbb{Z}[q^\frac{1}{n}, q^{-\frac{1}{n}}]$-linear category monoidally generated by \\
$\bullet $ objects $\{k_{\uparrow} \ | \ k\in \mathbb{Z}_{\geq 0} \}$,\\
$\bullet$ morphisms 
\[
\xy
(0,0)*{
\begin{tikzpicture}[scale=.3]
\draw [very thick, directed=.55] (0, .75) to (0,2.5);
\draw [very thick, directed=.45] (1,-1) to [out=90,in=330] (0,.75);
\draw [very thick, directed=.45] (-1,-1) to [out=90,in=210] (0,.75); 
\node at (0, 3) {\tiny $k{+}l$};
\node at (-1,-1.5) {\tiny $k$};
\node at (1,-1.5) {\tiny $l$};
\end{tikzpicture}  
};
\endxy : k_{\uparrow} \otimes l_{\uparrow} \rightarrow (k+l)_{\uparrow} , \quad
\xy
(0,0)*{
\begin{tikzpicture}[scale=.3]
\draw [very thick, directed=.55] (0,-1) to (0,.75);
\draw [very thick, directed=.65] (0,.75) to [out=30,in=270] (1,2.5);
\draw [very thick, directed=.65] (0,.75) to [out=150,in=270] (-1,2.5); 
\node at (0, -1.5) {\tiny $k{+}l$};
\node at (-1,3) {\tiny $k$};
\node at (1,3) {\tiny $l$};
\end{tikzpicture}
};
\endxy :(k+l)_{\uparrow} \rightarrow k_{\uparrow} \otimes l_{\uparrow} ,
\]
\[
\xy
(0,0)*{
\begin{tikzpicture}[scale=.3]
\draw [very thick, ->] (-1,-1) to (1,1);
\draw [very thick] (1,-1) to (0.25,-0.25);
\draw [very thick, ->] (-0.25,0.25) to (-1,1);
\node at (-1,-1.5) {\tiny $k$};
\node at (1,-1.5) {\tiny $l$};
\node at (1,1.45) {\tiny $k$};
\node at (-1,1.45) {\tiny $l$};
\end{tikzpicture}
};
\endxy :k_{\uparrow} \otimes l_{\uparrow} \rightarrow l_{\uparrow} \otimes k_{\uparrow} , \quad 
\xy
(0,0)*{
\begin{tikzpicture}[scale=.3]
\draw [very thick, ->] (1,-1) to (-1,1);
\draw [very thick] (-1,-1) to (-0.25,-0.25);
\draw [very thick, ->] (0.25,0.25) to (1,1);
\node at (-1,-1.5) {\tiny $k$};
\node at (1,-1.5) {\tiny $l$};
\node at (1,1.45) {\tiny $k$};
\node at (-1,1.45) {\tiny $l$};
\end{tikzpicture}
};
\endxy: k_{\uparrow} \otimes l_{\uparrow} \rightarrow l_{\uparrow} \otimes k_{\uparrow}.
\]
\end{Definition}

\begin{Definition}\label{int up-webs} 
The \textit{symmetric upward pointing category of integral} $\mathfrak{gl}_n$-webs $\mathbf{SWeb}_{\uparrow}^{\text{int}}(\mathfrak{gl}_n)$ is the quotient of $\mathbf{fSWeb}_{\uparrow}^{\text{int}}(\mathfrak{gl}_n)$ by the following relations:\\
\begin{flushleft}
Associativity and coassociativity
\end{flushleft} 
\begin{equation}\label{int (co)assoc}
\xy,
(0,0)*{
\begin{tikzpicture}[scale=.3]
\draw [very thick, directed=.45] (0,.75) to [out=90,in=220] (1,2.5);
\draw [very thick, directed=.45] (1,-1) to [out=90,in=330] (0,.75);
\draw [very thick, directed=.45] (-1,-1) to [out=90,in=210] (0,.75);
\draw [very thick, directed=.45] (3,-1) to [out=90,in=330] (1,2.5);
\draw [very thick, directed=.45] (1,2.5) to (1,4.25);
\node at (-1,-1.5) {\tiny $h$};
\node at (1,-1.5) {\tiny $m$};
\node at (-1.375,1.5) {\tiny $h{+}m$};
\node at (3,-1.5) {\tiny $l$};
\node at (1,4.75) {\tiny $h{+}m{+}l$};
\end{tikzpicture}
};
\endxy=\xy
(0,0)*{
\begin{tikzpicture}[scale=.3]
\draw [very thick, directed=.45] (0,.75) to [out=90,in=340] (-1,2.5);
\draw [very thick, directed=.45] (-1,-1) to [out=90,in=210] (0,.75);
\draw [very thick, directed=.45] (1,-1) to [out=90,in=330] (0,.75);
\draw [very thick, directed=.45] (-3,-1) to [out=90,in=220] (-1,2.5);
\draw [very thick, directed=.45] (-1,2.5) to (-1,4.25);
\node at (1,-1.5) {\tiny $l$};
\node at (-1,-1.5) {\tiny $m$};
\node at (1.25,1.5) {\tiny $m{+}l$};
\node at (-3,-1.5) {\tiny $h$};
\node at (-1,4.75) {\tiny $h{+}m{+}l$};
\end{tikzpicture}
};
\endxy
\quad\text{and}\quad
\xy
(0,0)*{\rotatebox{180}{
\begin{tikzpicture}[scale=.3]
\draw [very thick, rdirected=.55] (0,.75) to [out=90,in=220] (1,2.5);
\draw [very thick, rdirected=.55] (1,-1) to [out=90,in=330] (0,.75);
\draw [very thick, rdirected=.55] (-1,-1) to [out=90,in=210] (0,.75);
\draw [very thick, rdirected=.55] (3,-1) to [out=90,in=330] (1,2.5);
\draw [very thick, rdirected=.55] (1,2.5) to (1,4.25);
\node at (-1,-1.5) {\rotatebox{180}{\tiny $l$}};
\node at (1,-1.5) {\rotatebox{180}{\tiny $m$}};
\node at (-1.325,1.5) {\rotatebox{180}{\tiny $m\! +\! l$}};
\node at (3,-1.5) {\rotatebox{180}{\tiny $h$}};
\node at (1,4.75) {\rotatebox{180}{\tiny $h\! +\! m\! +\! l$}};
\end{tikzpicture}
}};
\endxy=
\xy
(0,0)*{\reflectbox{\rotatebox{180}{
\begin{tikzpicture}[scale=.3]
\draw [very thick, rdirected=.55] (0,.75) to [out=90,in=220] (1,2.5);
\draw [very thick, rdirected=.55] (1,-1) to [out=90,in=330] (0,.75);
\draw [very thick, rdirected=.55] (-1,-1) to [out=90,in=210] (0,.75);
\draw [very thick, rdirected=.55] (3,-1) to [out=90,in=330] (1,2.5);
\draw [very thick, rdirected=.55] (1,2.5) to (1,4.25);
\node at (-1,-1.5) {\reflectbox{\rotatebox{180}{\tiny $h$}}};
\node at (1,-1.5) {\reflectbox{\rotatebox{180}{\tiny $m$}}};
\node at (-1.325,1.5) {\reflectbox{\rotatebox{180}{\tiny $h\! +\! m$}}};
\node at (3,-1.5) {\reflectbox{\rotatebox{180}{\tiny $l$}}};
\node at (1,4.75) {\reflectbox{\rotatebox{180}{\tiny $h\! +\! m\! +\! l$}}};
\end{tikzpicture}
}}}; 
\endxy.
\end{equation}
Digon removal 
\begin{equation}\label{int digon}
\xy
(0,0)*{
\begin{tikzpicture}[scale=.3]
\draw [very thick, directed=.55] (0,.75) to (0,2.5);
\draw [very thick, directed=.55] (0,-2.75) to [out=30,in=330] (0,.75);
\draw [very thick, directed=.55] (0,-2.75) to [out=150,in=210] (0,.75);
\draw [very thick, directed=.55] (0,-4.5) to (0,-2.75);
\node at (0,-5) {\tiny $k{+}l$};
\node at (0,3) {\tiny $k{+}l$};
\node at (-2,-1) {\tiny $k$};
\node at (2,-1) {\tiny $l$};
\end{tikzpicture}
};
\endxy=\begin{bmatrix} k+l \\ l\end{bmatrix}\!\xy
(0,0)*{
\begin{tikzpicture}[scale=.3]
\draw [very thick, directed=.55] (0,-4.5) to (0,2.5);
\node at (0,-5) {\tiny $k{+}l$};
\node at (0,3) {\tiny $k{+}l$};
\end{tikzpicture}
}; 
\endxy.  
\end{equation}
Dumbbell-crossing relations 
\begin{equation} \label{int dum-cross}
\xy
(-5,0)*{
\begin{tikzpicture}[scale=.4]
\draw [very thick, directed=.6] (0,-1) to (0,.75);
\draw [very thick, directed=.75] (0,.75) to [out=30,in=270] (1,2.5);
\draw [very thick, directed =.75] (0,.75) to [out=150,in=270] (-1,2.5); 
\draw [very thick, directed=-.5] (1,-2.75) to [out=90,in=330] (0,-1);
\draw [very thick, directed=-.5] (-1,-2.75) to [out=90,in=210] (0,-1);
\node at (-1,3) {\tiny $r$};
\node at (1,3) {\tiny $s$};
\node at (-1,-3.25) {\tiny $k$};
\node at (1,-3.25) {\tiny $l$};
\end{tikzpicture}
};
\endxy =  \sum_{i-j=k-r}q^{-(k-i)(l-j)+ \frac{ij}{n}} 
\xy
(0,0)*{
\begin{tikzpicture}[scale=.4]
\draw [very thick,] (1,-1) to (-1,1);
\draw [very thick,crossline] (-1,-1) to (1,1);
\node at (-.7,-1.2) {\tiny $i$};
\node at (.7,-1.2) {\tiny $j$};
\node at (1,-3) {\tiny $l$};
\node at (-1,-3) {\tiny $k$};
\draw [very thick,directed=.9,directed=.2] (-1,-2.5) to (-1,2.5);
\draw [very thick,directed=.9,directed=.2] (1,-2.5) to (1,2.5);
\node at (-1,3) {\tiny $r$};
\node at (1,3) {\tiny $s$};
\end{tikzpicture}
};
\endxy,
\end{equation}
\[
\xy
(-5,0)*{
\begin{tikzpicture}[scale=.4]
\draw [very thick, directed=.6] (0,-1) to (0,.75);
\draw [very thick, directed=.75] (0,.75) to [out=30,in=270] (1,2.5);
\draw [very thick, directed =.75] (0,.75) to [out=150,in=270] (-1,2.5); 
\draw [very thick, directed=-.5] (1,-2.75) to [out=90,in=330] (0,-1);
\draw [very thick, directed=-.5] (-1,-2.75) to [out=90,in=210] (0,-1);
\node at (-1,3) {\tiny $r$};
\node at (1,3) {\tiny $s$};
\node at (-1,-3.25) {\tiny $k$};
\node at (1,-3.25) {\tiny $l$};
\end{tikzpicture}
};
\endxy =  \sum_{i-j=k-r}q^{(k-i)(l-j)- \frac{ij}{n}} 
\xy
(0,0)*{
\begin{tikzpicture}[scale=.4]
\draw [very thick,] (-1,-1) to (1,1);
\draw [very thick,crossline] (1,-1) to (-1,1);
\node at (-.7,-1.2) {\tiny $i$};
\node at (.7,-1.2) {\tiny $j$};
\node at (1,-3) {\tiny $l$};
\node at (-1,-3) {\tiny $k$};
\draw [very thick,directed=.9,directed=.2] (-1,-2.5) to (-1,2.5);
\draw [very thick,directed=.9,directed=.2] (1,-2.5) to (1,2.5);
\node at (-1,3) {\tiny $r$};
\node at (1,3) {\tiny $s$};
\end{tikzpicture}
};
\endxy.
\]
\end{Definition}

The relations in \autoref{int dum-cross} are motivated by \cite[(4.23)]{BrEnAiEtOs-semisimple-tilting}. Note that the crossings can be omitted as generators, as we will show in \autoref{cross are redundant gens}. We only use them above to give a shorter description of 
the construction. Note also that we do \textit{not} assume that the crossings are invertible. This will be a later consequence. We will actually prove that $\mathbf{SWeb}_{\uparrow}^{\text{int}}(\mathfrak{gl}_n)$ is a braided category.\\

Now, since we are working over $\mathbb{Z}[q^\frac{1}{n}, q^{-\frac{1}{n}}]$, the digon removal relation \autoref{int digon} is \textit{not} an invertible operation. Consequently we cannot use the idea of exploding edges (from the previous chapters) since this requires the inversion of quantum binomials. This is the reason why in \autoref{int up-webs} our web generators are thick to begin with.

We point out here that \ochanged{many properties of generic webs also hold integrally}. One of these is the following.

\begin{Lemma}\label{int sq switch}  
Thick square switches hold in $\mathbf{SWeb}_{\uparrow}^{\text{int}}(\mathfrak{gl}_n)$. Namely,
\begin{equation*}\
\xy
(0,0)*{
\begin{tikzpicture}[scale=.3]
\draw [very thick, directed=.55] (-2,-4) to (-2,-2);
\draw [very thick, directed=1] (-2,-2) to (-2,0.25);
\draw [very thick, directed=.55] (2,-4) to (2,-2);
\draw [very thick, directed=1] (2,-2) to (2,0.25);
\draw [very thick, directed=.55] (-2,-2) to (2,-2);
\draw [very thick] (-2,0.25) to (-2,2);
\draw [very thick, directed=.55] (-2,2) to (-2,4);
\draw [very thick] (2,0.25) to (2,2);
\draw [very thick, directed=.55] (2,2) to (2,4);
\draw [very thick, rdirected=.55] (-2,2) to (2,2);
\node at (-2,-4.5) {\tiny $k$};
\node at (2,-4.5) {\tiny $l$};
\node at (-2,4.5) {\tiny $k{-}h{+}g$};
\node at (2,4.5) {\tiny $l{+}h{-}g$};
\node at (-3.5,0) {\tiny $k{-}h$};
\node at (3.5,0) {\tiny $l{+}h$};
\node at (0,-1.25) {\tiny $h$};
\node at (0,2.75) {\tiny $g$};
\end{tikzpicture}
};
\endxy = \sum_{i} \begin{bmatrix} k-l+g -h\\ i \end{bmatrix}
\xy
(0,0)*{
\begin{tikzpicture}[scale=.3]
\draw [very thick, directed=.55] (-2,-4) to (-2,-2);
\draw [very thick, directed=1] (-2,-2) to (-2,0.25);
\draw [very thick, directed=.55] (2,-4) to (2,-2);
\draw [very thick, directed=1] (2,-2) to (2,0.25);
\draw [very thick, rdirected=.55] (-2,-2) to (2,-2);
\draw [very thick] (-2,0.25) to (-2,2);
\draw [very thick, directed=.55] (-2,2) to (-2,4);
\draw [very thick] (2,0.25) to (2,2);
\draw [very thick, directed=.55] (2,2) to (2,4);
\draw [very thick, directed=.55] (-2,2) to (2,2);
\node at (-2,-4.5) {\tiny $k$};
\node at (2,-4.5) {\tiny $l$};
\node at (-2,4.5) {\tiny $k{-}h{+}g$};
\node at (2,4.5) {\tiny $l{-}g{+}h$};
\node at (-4,0) {\tiny $k{+}g{-}i$};
\node at (4,0) {\tiny $l{-}g{+}i$};
\node at (0,-1.25) {\tiny $g{-}i$};
\node at (0,2.75) {\tiny $h{-}i$};
\end{tikzpicture}
};
\endxy, \quad
\end{equation*}
\begin{equation*}
\xy
(0,0)*{
\begin{tikzpicture}[scale=.3]
\draw [very thick, directed=.55] (-2,-4) to (-2,-2);
\draw [very thick, directed=1] (-2,-2) to (-2,0.25);
\draw [very thick, directed=.55] (2,-4) to (2,-2);
\draw [very thick, directed=1] (2,-2) to (2,0.25);
\draw [very thick, rdirected=.55] (-2,-2) to (2,-2);
\draw [very thick] (-2,0.25) to (-2,2);
\draw [very thick, directed=.55] (-2,2) to (-2,4);
\draw [very thick] (2,0.25) to (2,2);
\draw [very thick, directed=.55] (2,2) to (2,4);
\draw [very thick, directed=.55] (-2,2) to (2,2);
\node at (-2,-4.5) {\tiny $k$};
\node at (2,-4.5) {\tiny $l$};
\node at (-2,4.5) {\tiny $k{+}h{-}g$};
\node at (2,4.5) {\tiny $l{-}h{+}g$};
\node at (-4,0) {\tiny $k{+}h$};
\node at (3.75,0) {\tiny $l{-}h$};
\node at (0,-1.25) {\tiny $h$};
\node at (0,2.75) {\tiny $g$};
\end{tikzpicture}
};
\endxy= \sum_{i} \begin{bmatrix} l-k+g-h \\ i \end{bmatrix}
\xy
(0,0)*{
\begin{tikzpicture}[scale=.3]
\draw [very thick, directed=.55] (-2,-4) to (-2,-2);
\draw [very thick, directed=1] (-2,-2) to (-2,0.25);
\draw [very thick, directed=.55] (2,-4) to (2,-2);
\draw [very thick, directed=1] (2,-2) to (2,0.25);
\draw [very thick, directed=.55] (-2,-2) to (2,-2);
\draw [very thick] (-2,0.25) to (-2,2);
\draw [very thick, directed=.55] (-2,2) to (-2,4);
\draw [very thick] (2,0.25) to (2,2);
\draw [very thick, directed=.55] (2,2) to (2,4);
\draw [very thick, rdirected=.55] (-2,2) to (2,2);
\node at (-2,-4.5) {\tiny $k$};
\node at (2,-4.5) {\tiny $l$};
\node at (-2,4.5) {\tiny $k{+}h{-}g$};
\node at (2,4.5) {\tiny $l{-}h{+}g$};
\node at (-4,0) {\tiny $k{-}g{+}i$};
\node at (4,0) {\tiny $l{+}g{-}i$};
\node at (0,-1.25) {\tiny $g{-}i$};
\node at (0,2.75) {\tiny $h{-}i$};
\end{tikzpicture}
};
\endxy.
\end{equation*}
\end{Lemma}

\begin{proof}
These can be derived from relations \autoref{int (co)assoc}-\autoref{int dum-cross} and $q$-combinatorics. Indeed, we use relation \autoref{int dum-cross} to write the LHS of the first relation as
\[
\sum_{t=\text{max}(0,g-l)}^{\text{min}(g,h)}q^{-t(l-g+t)+ \frac{(h-t)(g-t)}{n}} 
\xy
(0,0)*{
\begin{tikzpicture}[scale=.5]
\draw [very thick,] (1,-1) to (-1,1);
\draw [very thick,crossline] (-1,-1) to (1,1);
\node at (-.2,-1.1) {\tiny $h{-}t$};
\node at (-.2,1.1) {\tiny $g{-}t$};
\node at (-1.3,0) {\tiny $t$};
\node at (1,-4.5) {\tiny $l$};
\node at (-2,-4.5) {\tiny $k$};
\draw [very thick,directed=.9,directed=.2] (-1,-2.5) to (-1,2.5);
\draw [very thick,directed=.9,directed=.2] (1,-4) to (1,4);
\draw[very thick] (-1,2.5) to [in=0, out=90] (-2,3.5) to [in=90,out=180] (-3,2.5);
\draw[very thick] (-1,-2.5) to [in=0, out=270] (-2,-3.5) to [in=270,out=180] (-3,-2.5);
\draw [very thick,] (-3,-2.5) to (-3,2.5);
\draw [very thick,] (-2,3.5) to (-2,4);
\draw [very thick,] (-2,-3.5) to (-2,-4);
\node at (-2,4.5) {\tiny $k{-}h+g$};
\node at (1,4.5) {\tiny $l{+}h{-}g$};
\end{tikzpicture}
};
\endxy
\]
\[ \stackrel{\autoref{int (co)assoc} \ \text{and} \ \autoref{int digon}}{=}
\sum_{t=\text{max}(0,g-l)}^{\text{min}(g,h)}q^{-t(l-g+t)+ \frac{(g-t)(h-t)}{n}} \begin{bmatrix} k-h+t \\ t \end{bmatrix}\!
\xy
(0,0)*{
\begin{tikzpicture}[scale=.5]
\draw [very thick,] (1,-1) to (-1,1);
\draw [very thick,crossline] (-1,-1) to (1,1);
\node at (-.2,-1.1) {\tiny $h{-}t$};
\node at (-.2,1.1) {\tiny $g{-}t$};
\node at (1,-3) {\tiny $l$};
\node at (-1,-3) {\tiny $k$};
\draw [very thick,directed=.9,directed=.2] (-1,-2.5) to (-1,2.5);
\draw [very thick,directed=.9,directed=.2] (1,-2.5) to (1,2.5);
\node at (-1.2,3) {\tiny $k{-}h{+}g$};
\node at (1.5,3) {\tiny $l{+}h{-}g$};
\end{tikzpicture}
};
\endxy.
\]
Similarly, the RHS is equal to 
\[
\sum_{i=0}^{\text{min}(g,h)} \begin{bmatrix} k-l+g -h\\ i \end{bmatrix} \sum_{s=0}^{\text{min}(g,h)-i} q^{-s(k-h+i+s)+ \frac{(g-i-s)(h-i-s)}{n}} 
\xy
(0,0)*{
\begin{tikzpicture}[scale=.5]
\draw [very thick,] (1,-1) to (-1,1);
\draw [very thick,crossline] (-1,-1) to (1,1);
\node at (-.1,-1.1) {\tiny $h{-}i{-}s$};
\node at (-.1,1.1) {\tiny $g{-}i{-}s$};
\node at (2,-4.5) {\tiny $l$};
\node at (-1,-4.5) {\tiny $k$};
\node at (1.3,0) {\tiny $s$};
\draw [very thick,directed=.9,directed=.2] (-1,-4) to (-1,4);
\draw [very thick,directed=.9,directed=.2] (1,-2.5) to (1,2.5);
\draw[very thick] (1,2.5) to [in=180, out=90] (2,3.5) to [in=90,out=0] (3,2.5);
\draw[very thick] (1,-2.5) to [in=180, out=270] (2,-3.5) to [in=270,out=0] (3,-2.5);
\draw [very thick,] (3,-2.5) to (3,2.5);
\draw [very thick,] (2,-3.5) to (2,-4);
\draw [very thick,] (2,3.5) to (2,4);
\node at (-1.2,4.5) {\tiny $k{-}h{+}g$};
\node at (2,4.5) {\tiny $l{+}h{-}g$};
\end{tikzpicture}
};
\endxy
\]
\[= \sum_{i=0}^{\text{min}(g,h)} \begin{bmatrix} k-l+g -h\\ i \end{bmatrix}
\sum_{s=0}^{\text{min}(g,h)-i}q^{-s(k-h+i+s)+ \frac{(g-i-s)(h-i-s)}{n}}  \begin{bmatrix} l-g+i+s \\ s \end{bmatrix}\!
\xy
(0,0)*{
\begin{tikzpicture}[scale=.5]
\draw [very thick,] (1,-1) to (-1,1);
\draw [very thick,crossline] (-1,-1) to (1,1);
\node at (-.1,-1.1) {\tiny $h{-}i{-}s$};
\node at (-.1,1.1) {\tiny $g{-}i{-}s$};
\node at (1,-3) {\tiny $l$};
\node at (-1,-3) {\tiny $k$};
\draw [very thick,directed=.9,directed=.2] (-1,-2.5) to (-1,2.5);
\draw [very thick,directed=.9,directed=.2] (1,-2.5) to (1,2.5);
\node at (-1.2,3) {\tiny $k{-}h{+}g$};
\node at (1.5,3) {\tiny $l{+}h{-}g$};
\end{tikzpicture}
};
\endxy
\]

\[ \stackrel{\text{\cite[Lemma 3.1.6]{La-gln-webs}}}{=} 
\sum_{i=\text{max}(0,g-l)}^{\text{min}(g,h)} \begin{bmatrix} k-h+i+s\\ i+s \end{bmatrix}
q^{-(i+s)(l-g+i+s)+ \frac{(g-i-s)(h-i-s)}{n}}  
\xy
(0,0)*{
\begin{tikzpicture}[scale=.5]
\draw [very thick,] (1,-1) to (-1,1);
\draw [very thick,crossline] (-1,-1) to (1,1);
\node at (-.1,-1.1) {\tiny $h{-}i{-}s$};
\node at (-.1,1.1) {\tiny $g{-}i{-}s$};
\node at (1,-3) {\tiny $l$};
\node at (-1,-3) {\tiny $k$};
\draw [very thick,directed=.9,directed=.2] (-1,-2.5) to (-1,2.5);
\draw [very thick,directed=.9,directed=.2] (1,-2.5) to (1,2.5);
\node at (-1.2,3) {\tiny $k{-}h{+}g$};
\node at (1.5,3) {\tiny $l{+}h{-}g$};
\end{tikzpicture}
};
\endxy
\]
\[ \stackrel{\text{put} \ i+s=t}{=}
\sum_{t=\text{max}(0,g-l)}^{\text{min}(g,h)}q^{-t(l-g+t)+ \frac{(g-t)(h-t)}{n}} \begin{bmatrix} k-h+t \\ t \end{bmatrix}\!
\xy
(0,0)*{
\begin{tikzpicture}[scale=.5]
\draw [very thick,] (1,-1) to (-1,1);
\draw [very thick,crossline] (-1,-1) to (1,1);
\node at (-.2,-1.1) {\tiny $h{-}t$};
\node at (-.2,1.1) {\tiny $g{-}t$};
\node at (1,-3) {\tiny $l$};
\node at (-1,-3) {\tiny $k$};
\draw [very thick,directed=.9,directed=.2] (-1,-2.5) to (-1,2.5);
\draw [very thick,directed=.9,directed=.2] (1,-2.5) to (1,2.5);
\node at (-1.2,3) {\tiny $k{-}h{+}g$};
\node at (1.5,3) {\tiny $l{+}h{-}g$};
\end{tikzpicture}
};
\endxy.
\]
The other square switch relation is proved analogously. 
\end{proof}

\begin{Lemma}\label{cross are redundant gens} 
The crossings from \autoref{int up-webs} satisfy the following identities
\[\raisebox{-.05cm}{\xy
(0,0)*{
\begin{tikzpicture}[scale=.3]
\draw [very thick, ->,] (-1,-1) to (1,1);
\draw [very thick, ->,] (-0.25,0.25) to (-1,1);
\draw [very thick,] (0.25,-0.25) to (1,-1);
\node at (-1,-1.5) {\tiny $k$};
\node at (1,-1.5) {\tiny $l$};
\end{tikzpicture}
};
\endxy} =
q^{-\frac{kl}{n}}\!\! \ \sum_{t=0}^{\text{min}(k,l)} (-q)^{-t}
\xy
(0,0)*{
\begin{tikzpicture}[scale=.3]
\draw [very thick,] (-2,-4) to (-2,-2);
\draw [very thick,directed=1] (-2,-2) to (-2,0.25);
\draw [very thick,] (2,-4) to (2,-2);
\draw [very thick,directed=1] (2,-2) to (2,0.25);
\draw [very thick, rdirected=.55,] (-2,-2) to (2,-2);
\draw [very thick,] (-2,0.25) to (-2,2);
\draw [very thick,] (-2,2) to (-2,4);
\draw [very thick,] (2,0.25) to (2,2);
\draw [very thick,] (2,2) to (2,4);
\draw [very thick, directed=.55,] (-2,2) to (2,2);
\node at (-2,-4.5) {\tiny $k$};
\node at (2,-4.5) {\tiny $l$};
\node at (-2,4.5) {\tiny $l$};
\node at (2,4.5) {\tiny $k$};
\node at (-4,0) {\tiny $k{+}l{-}t$};
\node at (2.75,0) {\tiny $t$};
\node at (0,-1.25) {\tiny $l-t$};
\node at (0,2.75) {\tiny $k-t$};
\end{tikzpicture}
};
\endxy,
\]
\[
\raisebox{-.05cm}{\xy
(0,0)*{
\begin{tikzpicture}[scale=.3]
\draw [very thick, ->,] (1,-1) to (-1,1);
\draw [very thick, ->,] (0.25,0.25) to (1,1);
\draw [very thick,] (-0.25,-0.25) to (-1,-1);
\node at (-1,-1.5) {\tiny $l$};
\node at (1,-1.5) {\tiny $k$};
\end{tikzpicture}
};
\endxy}
=
q^{\frac{kl}{n}}\!\!\sum_{t=0}^{\text{min}(k,l)}(-q)^{t}
\xy
(0,0)*{
\begin{tikzpicture}[scale=.3]
\draw [very thick,] (-2,-4) to (-2,-2);
\draw [very thick,directed=1] (-2,-2) to (-2,0.25);
\draw [very thick,] (2,-4) to (2,-2);
\draw [very thick,directed=1] (2,-2) to (2,0.25);
\draw [very thick, directed=.55,] (2,-2) to (-2,-2);
\draw [very thick,] (-2,0.25) to (-2,2);
\draw [very thick,] (-2,2) to (-2,4);
\draw [very thick,] (2,0.25) to (2,2);
\draw [very thick,] (2,2) to (2,4);
\draw [very thick, rdirected=.55,] (2,2) to (-2,2);
\node at (-2,-4.5) {\tiny $l$};
\node at (2,-4.5) {\tiny $k$};
\node at (-2,4.5) {\tiny $k$};
\node at (2,4.5) {\tiny $l$};
\node at (-4,0) {\tiny $l{+}k{-}t$};
\node at (2.75,0) {\tiny $t$};
\node at (0,-1.25) {\tiny $k-t$};
\node at (0,2.75) {\tiny $l-t$};
\end{tikzpicture}
};
\endxy,
\]
meaning that they are redundant as generators of $\mathbf{SWeb}_{\uparrow}^{\text{int}}(\mathfrak{gl}_n)$.
\end{Lemma}

The proof of \autoref{cross are redundant gens} is similar to that of \autoref{sq switch implies diggon-cross}. We prove the latter in the end of the section.

\begin{Lemma}\label{sq switch implies diggon-cross}
We can deduce relations \autoref{int dum-cross} from \autoref{int (co)assoc}, \autoref{int digon} and the thick square switches from \autoref{int sq switch}.
\end{Lemma}

\begin{proof}
We show that \autoref{int (co)assoc}, \autoref{int digon} and thick square switches from \autoref{int sq switch} imply the first equality in \autoref{int dum-cross}. The second equality can be proved similarly. At this point, we give the quantum version of \cite[Appendix A]{BrEnAiEtOs-semisimple-tilting}. First, note that up to symmetry, it suffices to prove that the following holds in $\mathbf{SWeb}_{\uparrow}(\mathfrak{gl}_n)$, for $r\geq 0$; 
\[
\xy
(-5,0)*{
\begin{tikzpicture}[scale=.4]
\draw [very thick, directed=.6] (0,-1) to (0,.75);
\draw [very thick, directed=.75] (0,.75) to [out=30,in=270] (1,2.5);
\draw [very thick, directed =.75] (0,.75) to [out=150,in=270] (-1,2.5); 
\draw [very thick, directed=-.5] (1,-2.75) to [out=90,in=330] (0,-1);
\draw [very thick, directed=-.5] (-1,-2.75) to [out=90,in=210] (0,-1);
\node at (-1,3) {\tiny $l$};
\node at (1,3) {\tiny $k+r$};
\node at (-1,-3.25) {\tiny $k$};
\node at (1,-3.25) {\tiny $l+r$};
\end{tikzpicture}
};
\endxy =  \sum_{i-j=k-l}q^{-(k-i)(l+r-j)+ \frac{ij}{n}} 
\xy
(0,0)*{
\begin{tikzpicture}[scale=.4]
\draw [very thick] (-1,-1) to (1,1);
\draw [very thick] (1,-1) to (0.25,-0.25);
\draw [very thick] (-0.25,0.25) to (-1,1);
\node at (-.7,-1.2) {\tiny $i$};
\node at (.7,-1.2) {\tiny $j$};
\node at (1,-3) {\tiny $l+r$};
\node at (-1,-3) {\tiny $k$};
\draw [very thick,directed=.9,directed=.2] (-1,-2.5) to (-1,2.5);
\draw [very thick,directed=.9,directed=.2] (1,-2.5) to (1,2.5);
\node at (-1,3) {\tiny $l$};
\node at (1,3) {\tiny $k+r$};
\end{tikzpicture}
};
\endxy.
\]
For simplicity, we relabel some of the edges above and the final version of relation \autoref{int dum-cross} that we are aiming to prove is 
\begin{equation*}
\xy
(-5,0)*{
\begin{tikzpicture}[scale=.4]
\draw [very thick, directed=.6] (0,-1) to (0,.75);
\draw [very thick, directed=.75] (0,.75) to [out=30,in=270] (1,2.5);
\draw [very thick, directed =.75] (0,.75) to [out=150,in=270] (-1,2.5); 
\draw [very thick, directed=-.5] (1,-2.75) to [out=90,in=330] (0,-1);
\draw [very thick, directed=-.5] (-1,-2.75) to [out=90,in=210] (0,-1);
\node at (-1,3) {\tiny $l$};
\node at (1,3) {\tiny $k+r$};
\node at (-1,-3.25) {\tiny $k$};
\node at (1,-3.25) {\tiny $l+r$};
\end{tikzpicture}
};
\endxy =  \sum_{s=0}^{\text{min}(k,l)} q^{-s(s+r)+ \frac{(k-s)(l-s)}{n}} 
\xy
(0,0)*{
\begin{tikzpicture}[scale=.4]
\draw [very thick] (-1,-1) to (1,1);
\draw [very thick] (1,-1) to (0.25,-0.25);
\draw [very thick] (-0.25,0.25) to (-1,1);
\node at (-1.5,-.3) {\tiny $s$};
\node at (1,-3) {\tiny $l+r$};
\node at (-1,-3) {\tiny $k$};
\draw [very thick,directed=.9,directed=.2] (-1,-2.5) to (-1,2.5);
\draw [very thick,directed=.9,directed=.2] (1,-2.5) to (1,2.5);
\node at (-1,3) {\tiny $l$};
\node at (1,3) {\tiny $k+r$};
\end{tikzpicture}
};
\endxy. 
\end{equation*}
Using \autoref{cross are redundant gens} (which is proved independently) for the  $(k-s,l-s)$ braiding, we have that the right-hand side of the above relation is equal to 
\[
\sum_{s=0}^{\text{min}(k,l)} q^{-s(s+r)}\sum_{t=0}^{\text{min}(k,l)-s}(-q)^{-t} 
\xy
(0,0)*{
\begin{tikzpicture}[scale=.25]
\draw [very thick, directed=.55] (-2,-4) to (-2,-2);
\draw [very thick, directed=1] (-2,-2) to (-2,0.25);
\draw [very thick, directed=.55] (2,-4) to (2,-2);
\draw [very thick, directed=1] (2,-2) to (2,0.25);
\draw [very thick, directed=.55] (-2,-2) to (2,-2);
\draw [very thick] (-2,0.25) to (-2,2);
\draw [very thick, directed=.55] (-2,2) to (-2,4);
\draw [very thick] (2,0.25) to (2,2);
\draw [very thick, directed=.55] (2,2) to (2,4);
\draw [very thick, rdirected=.55] (-2,2) to (2,2);
\draw [very thick, ] (6,4) to [out=90,in=0] (4,6);
\draw [very thick, ] (4,6) to [out=180,in=90] (2,4);
\draw [very thick, ] (6,-4) to [out=270,in=0] (4,-6);
\draw [very thick, ] (4,-6) to [out=180,in=270] (2,-4);
\draw [very thick,]  (6,-4) to (6,4);
\draw [very thick, ] (-2,4) to [out=90,in=0] (-4,6);
\draw [very thick, ] (-4,6) to [out=180,in=90] (-6,4);
\draw [very thick, ] (-2,-4) to [out=270,in=0] (-4,-6);
\draw [very thick, ] (-4,-6) to [out=180,in=270] (-6,-4);
\draw [very thick, ] (-6,-4) to (-6,4);
\node at (-3,0) {\tiny $t$};
\node at (0,-1.2) {\tiny $k{-}s{-}t$};
\node at (0,3) {\tiny $l{-}s{-t}$};
\draw [very thick,directed=.5]  (4,6) to (4,8);
\node at (4,8.5) {\tiny $k+r$};
\draw [very thick,rdirected=.5]  (4,-6) to (4,-8);
\node at (4,-8.5) {\tiny $l+r$};
\draw [very thick,directed=.5]  (-4,6) to (-4,8);
\node at (-4,8.5) {\tiny $l$};
\draw [very thick,rdirected=.5]  (-4,-6) to (-4,-8);
\node at (-4,-8.5) {\tiny $k$};
\node at (-6.5,0) {\tiny $s$};
\node at (7.5,0) {\tiny $r+s$};
\end{tikzpicture}
};
\endxy 
\]
\[ \stackrel{\autoref{int (co)assoc},\autoref{int digon}}{=}
\sum_{s=0}^{\text{min}(k,l)} q^{-s(s+r)}\sum_{t=0}^{\text{min}(k,l)-s}(-q)^{-t} \begin{bmatrix} s+t \\ s\end{bmatrix}\!
\xy
(0,0)*{
\begin{tikzpicture}[scale=.25]
\draw [very thick, directed=.55] (-2,-4) to (-2,-2);
\draw [very thick, directed=1] (-2,-2) to (-2,0.25);
\draw [very thick, directed=.55] (2,-4) to (2,-2);
\draw [very thick, directed=1] (2,-2) to (2,0.25);
\draw [very thick, directed=.55] (-2,-2) to (2,-2);
\draw [very thick] (-2,0.25) to (-2,2);
\draw [very thick, directed=.55] (-2,2) to (-2,4);
\draw [very thick] (2,0.25) to (2,2);
\draw [very thick, directed=.55] (2,2) to (2,4);
\draw [very thick, rdirected=.55] (-2,2) to (2,2);
\draw [very thick, ] (6,4) to [out=90,in=0] (4,6);
\draw [very thick, ] (4,6) to [out=180,in=90] (2,4);
\draw [very thick, ] (6,-4) to [out=270,in=0] (4,-6);
\draw [very thick, ] (4,-6) to [out=180,in=270] (2,-4);
\draw [very thick,]  (6,-4) to (6,4);
\draw [very thick,]  (-2,-8) to (-2,8);
\node at (-3.5,0) {\tiny $s+t$};
\node at (0,-1.2) {\tiny $k{-}s{-}t$};
\node at (0,3) {\tiny $l{-}s{-}t$};
\draw [very thick,directed=.5]  (4,6) to (4,8);
\node at (4,8.5) {\tiny $k+r$};
\draw [very thick,rdirected=.5]  (4,-6) to (4,-8);
\node at (4,-8.5) {\tiny $l+r$};
\node at (-2,8.5) {\tiny $l$};
\node at (-2,-8.5) {\tiny $k$};
\node at (7.5,0) {\tiny $r+s$};
\end{tikzpicture}
};
\endxy 
\]
\[ \stackrel{\text{put} \ s+t=u}{=}
\sum_{s=0}^{\text{min}(k,l)} q^{-s(s+r)}\sum_{u=s}^{\text{min}(k,l)}(-q)^{s-u} \begin{bmatrix} u\\ s\end{bmatrix}\!
\xy
(0,0)*{
\begin{tikzpicture}[scale=.25]
\draw [very thick, directed=.55] (-2,-4) to (-2,-2);
\draw [very thick, directed=1] (-2,-2) to (-2,0.25);
\draw [very thick, directed=.55] (2,-4) to (2,-2);
\draw [very thick, directed=1] (2,-2) to (2,0.25);
\draw [very thick, rdirected=.55] (-2,-2) to (2,-2);
\draw [very thick] (-2,0.25) to (-2,2);
\draw [very thick, directed=.55] (-2,2) to (-2,4);
\draw [very thick] (2,0.25) to (2,2);
\draw [very thick, directed=.55] (2,2) to (2,4);
\draw [very thick, directed=.55] (-2,2) to (2,2);
\draw [very thick, ] (-2,4) to [out=90,in=0] (-4,6);
\draw [very thick, ] (-4,6) to [out=180,in=90] (-6,4);
\draw [very thick, ] (-2,-4) to [out=270,in=0] (-4,-6);
\draw [very thick, ] (-4,-6) to [out=180,in=270] (-6,-4);
\draw [very thick, ] (-6,-4) to (-6,4);
\draw [very thick, ] (2,-8) to (2,8);
\node at (0,-1.25) {\tiny $l-s$};
\node at (0,2.9) {\tiny $k-s$};
\node at (2,8.5) {\tiny $k+r$};
\node at (2,-8.5) {\tiny $l+r$};
\draw [very thick,directed=.5]  (-4,6) to (-4,8);
\node at (-4,8.5) {\tiny $l$};
\draw [very thick,rdirected=.5]  (-4,-6) to (-4,-8);
\node at (-4,-8.5) {\tiny $k$};
\node at (-6.5,0) {\tiny $u$};
\node at (3.75,0) {\tiny $r+s$};
\end{tikzpicture}
};
\endxy 
\]
\[  \stackrel{\text{square switch} }{=}
\sum_{s=0}^{\text{min}(k,l)} q^{-s(s+r)}\sum_{u=s}^{\text{min}(k,l)}(-q)^{s-u} \begin{bmatrix} u\\ s\end{bmatrix}\! \sum_{t=u-s}^{\text{min}(k,l)-s} \begin{bmatrix} u+r\\ t\end{bmatrix}\! 
\xy
(0,0)*{
\begin{tikzpicture}[scale=.25]
\draw [very thick, directed=.55] (-2,-4) to (-2,-2);
\draw [very thick, directed=1] (-2,-2) to (-2,0.25);
\draw [very thick, directed=.55] (2,-4) to (2,-2);
\draw [very thick, directed=1] (2,-2) to (2,0.25);
\draw [very thick, directed=.55] (-2,-2) to (2,-2);
\draw [very thick] (-2,0.25) to (-2,2);
\draw [very thick, directed=.55] (-2,2) to (-2,4);
\draw [very thick] (2,0.25) to (2,2);
\draw [very thick, directed=.55] (2,2) to (2,4);
\draw [very thick, rdirected=.55] (-2,2) to (2,2);
\draw [very thick, ] (-2,4) to [out=90,in=0] (-4,6);
\draw [very thick, ] (-4,6) to [out=180,in=90] (-6,4);
\draw [very thick, ] (-2,-4) to [out=270,in=0] (-4,-6);
\draw [very thick, ] (-4,-6) to [out=180,in=270] (-6,-4);
\draw [very thick, ] (-6,-4) to (-6,4);
\draw [very thick, ] (2,-8) to (2,8);
\node at (0,-1.2) {\tiny $k{-}s{-}t$};
\node at (0,3) {\tiny $l{-}s{-}t$};
\node at (2,8.5) {\tiny $k+r$};
\node at (2,-8.5) {\tiny $l+r$};
\draw [very thick,directed=.5]  (-4,6) to (-4,8);
\node at (-4,8.5) {\tiny $l$};
\draw [very thick,rdirected=.5]  (-4,-6) to (-4,-8);
\node at (-4,-8.5) {\tiny $k$};
\node at (-6.5,0) {\tiny $u$};
\end{tikzpicture}
};
\endxy 
\]
\[ \stackrel{\autoref{int (co)assoc},\autoref{int digon} }{=}
\scalebox{0.9}{$\sum_{s=0}^{\text{min}(k,l)} q^{-s(s+r)}\sum_{u=s}^{\text{min}(k,l)}(-q)^{s-u} \begin{bmatrix} u\\ s\end{bmatrix}\! \sum_{t=u-s}^{\text{min}(k,l)-s} \begin{bmatrix} u+r\\ t\end{bmatrix}\! \begin{bmatrix} s+t\\ u\end{bmatrix}\! 
\xy
(0,0)*{
\begin{tikzpicture}[scale=.3]
\draw [very thick, directed=.55] (-2,-4) to (-2,-2);
\draw [very thick, directed=1] (-2,-2) to (-2,0.25);
\draw [very thick, directed=.55] (2,-4) to (2,-2);
\draw [very thick, directed=1] (2,-2) to (2,0.25);
\draw [very thick, directed=.55] (-2,-2) to (2,-2);
\draw [very thick] (-2,0.25) to (-2,2);
\draw [very thick, directed=.55] (-2,2) to (-2,4);
\draw [very thick] (2,0.25) to (2,2);
\draw [very thick, directed=.55] (2,2) to (2,4);
\draw [very thick, rdirected=.55] (-2,2) to (2,2);
\node at (-2,-4.5) {\tiny $k$};
\node at (2,-4.5) {\tiny $l+r$};
\node at (-2,4.5) {\tiny $l$};
\node at (2,4.5) {\tiny $k+r$};
\node at (-3.5,0) {\tiny $s+t$};
\node at (0,-1.2) {\tiny $k{-}s{-}t$};
\node at (0,2.8) {\tiny $l{-}s{-t}$};
\end{tikzpicture}
};
\endxy$}
\]
\[ \stackrel{ \text{put} \ v=s+t}{=}
\sum_{s=0}^{\text{min}(k,l)} q^{-s(s+r)}\sum_{u=s}^{\text{min}(k,l)}(-1)^{s+u}q^{s-u} \begin{bmatrix} u\\ s\end{bmatrix}\! \sum_{v=u}^{\text{min}(k,l)} \begin{bmatrix} u+r\\ v-s\end{bmatrix}\! \begin{bmatrix} v \\ u\end{bmatrix}\! 
\xy
(0,0)*{
\begin{tikzpicture}[scale=.3]
\draw [very thick, directed=.55] (-2,-4) to (-2,-2);
\draw [very thick, directed=1] (-2,-2) to (-2,0.25);
\draw [very thick, directed=.55] (2,-4) to (2,-2);
\draw [very thick, directed=1] (2,-2) to (2,0.25);
\draw [very thick, directed=.55] (-2,-2) to (2,-2);
\draw [very thick] (-2,0.25) to (-2,2);
\draw [very thick, directed=.55] (-2,2) to (-2,4);
\draw [very thick] (2,0.25) to (2,2);
\draw [very thick, directed=.55] (2,2) to (2,4);
\draw [very thick, rdirected=.55] (-2,2) to (2,2);
\node at (-2,-4.5) {\tiny $k$};
\node at (2,-4.5) {\tiny $l+r$};
\node at (-2,4.5) {\tiny $l$};
\node at (2,4.5) {\tiny $k+r$};
\node at (-3.25,0) {\tiny $v$};
\node at (0,-1.25) {\tiny $k-v$};
\node at (0,2.75) {\tiny $l-v$};
\end{tikzpicture}
};
\endxy
\]
\[= \sum_{s=0}^{\text{min}(k,l)} \sum_{u=s}^{\text{min}(k,l)}\sum_{v=u}^{\text{min}(k,l)} (-1)^{s+u}q^{-s(s+r)}q^{s-u} \begin{bmatrix} u\\ s\end{bmatrix}\!  \begin{bmatrix} u+r\\ v-s\end{bmatrix}\! \begin{bmatrix} v \\ u\end{bmatrix}\! 
\xy
(0,0)*{
\begin{tikzpicture}[scale=.3]
\draw [very thick, directed=.55] (-2,-4) to (-2,-2);
\draw [very thick, directed=1] (-2,-2) to (-2,0.25);
\draw [very thick, directed=.55] (2,-4) to (2,-2);
\draw [very thick, directed=1] (2,-2) to (2,0.25);
\draw [very thick, directed=.55] (-2,-2) to (2,-2);
\draw [very thick] (-2,0.25) to (-2,2);
\draw [very thick, directed=.55] (-2,2) to (-2,4);
\draw [very thick] (2,0.25) to (2,2);
\draw [very thick, directed=.55] (2,2) to (2,4);
\draw [very thick, rdirected=.55] (-2,2) to (2,2);
\node at (-2,-4.5) {\tiny $k$};
\node at (2,-4.5) {\tiny $l+r$};
\node at (-2,4.5) {\tiny $l$};
\node at (2,4.5) {\tiny $k+r$};
\node at (-3.25,0) {\tiny $v$};
\node at (0,-1.25) {\tiny $k-v$};
\node at (0,2.75) {\tiny $l-v$};
\end{tikzpicture}
};
\endxy.
\]
Switching the order of summations, while keeping in mind that $(-1)^{s+u}=(-1)^v(-1)^s(-1)^{u-v}$, the last sum is equal to

\[
\sum_{v=0}^{\text{min}(k,l)} (-1)^v \sum_{s=0}^{v}(-1)^s \sum_{u=s}^{v} (-1)^{u-v}q^{-s(s+r)+s-u} \begin{bmatrix} v\\ s\end{bmatrix}\!  \begin{bmatrix} u+r\\ u-s,v-u\end{bmatrix}\! 
\xy
(0,0)*{
\begin{tikzpicture}[scale=.3]
\draw [very thick, directed=.55] (-2,-4) to (-2,-2);
\draw [very thick, directed=1] (-2,-2) to (-2,0.25);
\draw [very thick, directed=.55] (2,-4) to (2,-2);
\draw [very thick, directed=1] (2,-2) to (2,0.25);
\draw [very thick, directed=.55] (-2,-2) to (2,-2);
\draw [very thick] (-2,0.25) to (-2,2);
\draw [very thick, directed=.55] (-2,2) to (-2,4);
\draw [very thick] (2,0.25) to (2,2);
\draw [very thick, directed=.55] (2,2) to (2,4);
\draw [very thick, rdirected=.55] (-2,2) to (2,2);
\node at (-2,-4.5) {\tiny $k$};
\node at (2,-4.5) {\tiny $l+r$};
\node at (-2,4.5) {\tiny $l$};
\node at (2,4.5) {\tiny $k+r$};
\node at (-3.25,0) {\tiny $v$};
\node at (0,-1.25) {\tiny $k-v$};
\node at (0,2.75) {\tiny $l-v$};
\end{tikzpicture}
};
\endxy
\]
where we used the quantum trinomial coefficient. Further, ignoring the web, we write the last sum as 
\[
\sum_{v=0}^{\text{min}(k,l)} (-1)^v \sum_{s=0}^{v}(-1)^s q^{vs-s}\begin{bmatrix} v\\ s\end{bmatrix}\! \sum_{u=s}^{v} (-1)^{u-v}q^{-s(s+r)+s-u-vs+s}  \begin{bmatrix} u+r\\ u-s,v-u\end{bmatrix}\!
\]
and from \cite[Lemma \ 3.1.5]{La-gln-webs}, this is equal to 
\[
\sum_{v=0}^{\text{min}(k,l)} (-1)^v \delta_{v,0} \sum_{u=s}^{v} (-1)^{u-v}q^{-s(s+r+v-2)-u}  \begin{bmatrix} u+r\\ u-s,v-u\end{bmatrix}\!.
\]
This sum is nontrivial if and only if $v=0$. In that case also $u=s=0$ and 
\[
\sum_{u=s}^{v} (-1)^{u-v}q^{-s(s+r+v-2)-u}  \begin{bmatrix} u+r\\ u-s,v-u\end{bmatrix}\!=1,
\]
meaning that the right-hand side of the initial relation is equal to 
\[
\xy
(-5,0)*{
\begin{tikzpicture}[scale=.4]
\draw [very thick, directed=.6] (0,-1) to (0,.75);
\draw [very thick, directed=.75] (0,.75) to [out=30,in=270] (1,2.5);
\draw [very thick, directed =.75] (0,.75) to [out=150,in=270] (-1,2.5); 
\draw [very thick, directed=-.5] (1,-2.75) to [out=90,in=330] (0,-1);
\draw [very thick, directed=-.5] (-1,-2.75) to [out=90,in=210] (0,-1);
\node at (-1,3) {\tiny $l$};
\node at (1,3) {\tiny $k+r$};
\node at (-1,-3.25) {\tiny $k$};
\node at (1,-3.25) {\tiny $l+r$};
\end{tikzpicture}
};
\endxy
.
\]
This finishes the proof that \autoref{int dum-cross} can be deduced from \autoref{int (co)assoc}, \autoref{int digon} and thick square switches from \autoref{int sq switch}. 
\end{proof}

We stress that \autoref{int sq switch} and \autoref{sq switch implies diggon-cross} imply that the dumbbell-crossing and the thick square switch are equivalent relations in a certain sense.
Let us also \ochanged{stress that} if the underlying field of $\mathbf{SWeb}_{\uparrow}^{\text{int}}(\mathfrak{gl}_n)$ \ochanged{were} $\Bbbk(q^{\frac{1}{n}})$ then we could use thinner merges and splits as generators. Then we would be able to define thick merges, splits and braidings via explosion of strings as in the previous sections. Combining this with \autoref{int sq switch} and \autoref{sq switch implies diggon-cross}  we have the following.

\begin{Corollary}\label{int and complex webs are equivalent over c}
If we work over $\Bbbk(q^{\frac{1}{n}})$ (that is, any field where the quantum binomials are invertible), then $\mathbf{SWeb}_{\uparrow}^{\text{int}}(\mathfrak{gl}_n)$ is equivalent to $\mathbf{SWeb}_{\uparrow}(\mathfrak{gl}_n)$.
\end{Corollary}

While avoiding any scalar inversion, one can show that \ochanged{many of the} web properties hold also integrally (\ochanged{see} \cite[Section \ 6.2]{La-gln-webs}). As a consequence we have the following.

\begin{Theorem} 
There exists a full functor 
\[ 
\Gamma_{\text{sym}}^{\text{int}}: \mathbf{SWeb}_{\uparrow,\downarrow}^{\text{int}}(\mathfrak{gl}_n) \rightarrow \text{U}_{q}(\mathfrak{gl}_n) \text{-} \mathrm{\mathbf{fdMod}}_{S,S^*}
\]
defined similarly to $\Gamma_{\text{sym}}$.
\end{Theorem}

\begin{Remark}
\changed{It is natural to ask how the integral category considered here relates to its
generic counterpart. In particular, one may wonder whether the morphism spaces
over the ground ring are free, so that after base change to the fraction field
one recovers the generic category and the integral one defines a full lattice
inside it. At present we do not know whether this freeness holds in general.
However, based on the behaviour in related diagrammatic settings, we expect the
integral category to provide such a lattice.}
\end{Remark}

\end{document}